\documentclass[final]{siamltex}

\def\d{\delta} 
 
\def\e{{\epsilon}}

\def\norm#1{\|#1\|} 
\def\wb{{\bar w}} 

\usepackage{algorithm}
\usepackage{algorithmic}
\usepackage{amsmath}
\usepackage{amssymb}
\usepackage{caption}
\usepackage{comment}
\usepackage{graphicx}
\usepackage{float}
\usepackage[latin1]{inputenc}
\usepackage[pdfborder={0 0 0}, colorlinks=true, linkcolor=blue, citecolor=red]{hyperref}\hypersetup{pdfborder=0 0 0}
\usepackage{mathrsfs}
\usepackage{multicol}
\usepackage{multirow}
\usepackage{showkeys,cite}
\usepackage{subfigure}
\usepackage{wrapfig}
\usepackage{ulem}

\usepackage{soul,xargs}
\usepackage[pdftex,dvipsnames]{xcolor}

\usepackage[colorinlistoftodos,prependcaption,textsize=tiny]{todonotes}
\newcommandx{\question}[2][1=]{\todo[linecolor=red,backgroundcolor=red!20,bordercolor=red,#1]{#2}}
\newcommandx{\change}[2][1=]{\todo[linecolor=blue,backgroundcolor=blue!25,bordercolor=blue,#1]{#2}}
\newcommandx{\info}[2][1=]{\todo[linecolor=OliveGreen,backgroundcolor=OliveGreen!25,bordercolor=OliveGreen,#1]{#2}}
\newcommandx{\improve}[2][1=]{\todo[linecolor=Plum,backgroundcolor=Plum!25,bordercolor=Plum,#1]{#2}}
\newcommandx{\answer}[2][1=]{\todo[linecolor=blue,backgroundcolor=White!25,bordercolor=Plum,#1]{#2}}
\newcommandx{\suggest}[2][1=]{\todo[linecolor=blue,backgroundcolor=White!25,bordercolor=red,#1]{#2}}

\graphicspath{{./plots/}}
\begin{document}

\title{A Multilevel Block Preconditioner for the HDG Trace System Applied to Incompressible Resistive MHD}

    \author{Sriramkrishnan Muralikrishnan\thanks{Department of Aerospace Engineering and Engineering Mechanics, The University of Texas at Austin, Austin, TX 78712, USA. Current address: Paul Scherrer Institut, 5232 Villigen, Switzerland.} \and Stephen Shannon\thanks{Oden Institute for Computational Engineering and Sciences, The University of Texas at Austin, Austin, TX 78712, USA.} \and Tan Bui-Thanh\thanks{Department of Aerospace Engineering and Engineering Mechanics, and Oden Institute for Computational Engineering and Sciences, The University of Texas at Austin, Austin, TX 78712, USA.} \and John N. Shadid\thanks{Computational Mathematics Department, Sandia National Laboratories, Albuquerque, NM87185 and Department of Mathematics and Statistics, University of New Mexico, Albuquerque, NM 87131. }}

\bibliographystyle{siam}
\newcommand{\TODO}[1]{ \fbox{\parbox{3in}{\bf TODO: #1}}}

\newcommand{\grbf}[1] {\mbox{\boldmath${#1}$\unboldmath}}
\newcommand{\gbf}[1] {\mathbf{#1}}

\newcommand{\beq} {\begin{equation}}
\newcommand{\eeq} {\end{equation}}
\newcommand{\bdm} {\begin{displaymath}}
\newcommand{\edm} {\end{displaymath}}
\newcommand{\bit}{\begin{itemize}}
\newcommand{\eit}{\end{itemize}}
\newcommand{\bde}{\begin{description}}
\newcommand{\ede}{\end{description}}
\newcommand{\bce}{\begin{center}}
\newcommand{\ece}{\end{center}}
\newcommand{\ben} {\begin{enumerate}}
\newcommand{\een} {\end{enumerate}}
\newcommand{\bea} {\begin{eqnarray}}
\newcommand{\eea} {\end{eqnarray}}
\newcommand{\barr} {\begin{array}}
\newcommand{\earr} {\end{array}}
\newcommand{\bean} {\begin{eqnarray*}}
\newcommand{\eean} {\end{eqnarray*}}
\newcommand{\edoc} {\end{document}}
\newcommand{\On}{\hspace{0.2cm} \mbox{on} \hspace{0.2cm}}
\newcommand{\In}{\hspace{0.2cm} \mbox{in} \hspace{0.2cm}}
\newcommand{\Minimize}{\hspace{0.2cm} \mbox{minimize} \hspace{0.2cm}}
\newcommand{\Maximize}{\hspace{0.2cm} \mbox{maximize} \hspace{0.2cm}}
\newcommand{\SubjectTo}{\hspace{0.2cm} \mbox{subject} \hspace{0.15cm}
            \mbox{to} \hspace{0.2cm}}

\def\etal{{\it et al.~}}

\newcommand{\point}[1] {\ensuremath{\boldsymbol{#1}}}
\newcommand{\vvec}[1] {\ensuremath{\boldsymbol{#1}}}
\renewcommand{\vec}[1] {\ensuremath{\boldsymbol{#1}}}
\newcommand{\cvec}[1] {\ensuremath{\boldsymbol{{#1}}}}
\newcommand{\ten}[1] {\ensuremath{\boldsymbol{#1}}}
\newcommand{\tenfour}[1] {\ensuremath{\boldsymbol{\mathsf{#1}}}}
\newcommand{\comp}[1]{\begin{pmatrix}{ #1 }\end{pmatrix}}
\newcommand{\vv}{\ensuremath{\vec{v}}}
\newcommand{\vr}{\left(\ensuremath{\vec{r}}\right)}
\newcommand{\att}{\left(t\right)}
\newcommand{\tvv}{\ensuremath{\tilde{\vec{v}}}}
\newcommand{\ww}{\ensuremath{\vec{w}}}
\newcommand{\GG}{\ensuremath{\ten{G}}}
\newcommand{\FF}{\ensuremath{\boldsymbol{\mathfrak{F}}}}
\newcommand{\tEE}{\ensuremath{\tilde{\vec{E}}}}
\newcommand{\nn}{\ensuremath{\vec{n}}}
\renewcommand{\SS}{\ensuremath{\ten{S}}}
\newcommand{\tSS}{\ensuremath{\tilde{\ten{S}}}}
\newcommand\tK{{\bf K}}
\newcommand{\cp}{\ensuremath{{c_p}}}
\newcommand{\cs}{\ensuremath{{c_s}}}
\newcommand{\tcs}{\ensuremath{{\tilde{c}_s}}}
\newcommand{\Vp}{\ensuremath{{V^e_p}}}
\newcommand{\Vs}{\ensuremath{{V^e_s}}}
\newcommand{\Ss}{\ensuremath{{S^e_s}}}
\newcommand{\nxnx}[1] {\ensuremath{\nn\times\left(\nn\times{#1}\right)}}

\newcommand{\dd}[2] {\ensuremath{\frac{\partial {#1}}{\partial {#2}}}}
\newcommand{\DD}[2] {\ensuremath{\frac{d {#1}}{d {#2}}}}
\newcommand{\Grad} {\ensuremath{\nabla}}
\newcommand{\Gradv}[1]{\ensuremath{\nabla_{#1}}}
\newcommand{\Div} {\ensuremath{\nabla\cdot}}
\newcommand{\Divv}[1] {\ensuremath{\nabla_{#1}\cdot}}
\newcommand{\Curl} {\ensuremath{\nabla\times}}
\newcommand{\Cten}{\ensuremath{\tenfour{C}}}

\newcommand{\eval}[2][\right]{\relax
  \ifx#1\right\relax \left.\fi#2#1\rvert}

\providecommand{\abs}[1]{\left\lvert#1\right\rvert}
\providecommand{\norm}[1]{\left\lVert#1\right\rVert}

\newcommand{\Nel} {\ensuremath{{N_T}}}
\newcommand{\Ned} {\ensuremath{{N_\text{ed}}}}
\newcommand{\De} {\ensuremath{{\mathsf{D}^e}}}
\newcommand{\Dep} {\ensuremath{{\mathsf{D}^{e'}}}}
\newcommand{\Dhat} {\ensuremath{\hat{\mathsf{D}}}}
\newcommand{\Dehat} {\ensuremath{{\hat{\mathsf{D}}^e}}}
\newcommand{\half} {\ensuremath{\frac{1}{2}}}
\newcommand{\IN} {\ensuremath{{\mathcal{I}_N}}}
\newcommand{\INe} {\ensuremath{{\mathcal{I}_{N_e}}}}

\newcommand{\Ha}{\ensuremath{\text{Ha}}}
\newcommand{\Rey}{\ensuremath{\textup{Re}}}
\newcommand{\Rm}{\ensuremath{\textup{Rm}}}

\newcommand{\bb}{\mb{b}}
\newcommand{\bbn}{{\mb{b}^n}}
\newcommand{\bbt}{{\mb{b}^t}}

%
\newcommand{\bh}{\widehat{b}}
\newcommand{\ch}{\widehat{c}}
\newcommand{\chn}{{\ch^{n}_h}}
\newcommand{\fh}{\widehat{f}}
\newcommand{\gh}{\widehat{g}}
\newcommand{\rh}{\widehat{r}}
\newcommand{\sh}{\widehat{s}}
\newcommand{\shU}{\underline{\sh}}
\newcommand{\vh}{\widehat{v}}
\newcommand{\vhn}{\widehat{v}^n}
\newcommand{\zh}{\widehat{z}}

%
%
\newcommand{\bbh}{\widehat{\bb}}
\newcommand{\bbhn}{{\widehat{\bb}^n}}
\newcommand{\bbht}{{\widehat{\bb}^t}}
\newcommand{\cbh}{{\widehat{\cb}}}
\newcommand{\cbht}{{\widehat{\cb}^t}}
\newcommand{\cbhtU}{{\underline{\widehat{\cb}}^t}}
\newcommand{\Fbh}{\widehat{\Fb}}
\newcommand{\FBh}{\widehat{\FB}}
\newcommand{\GBh}{\widehat{\GB}}
\newcommand{\Jbh}{\widehat{\Jb}}
\newcommand{\Jbht}{{\widehat{\Jb}^t}}
\newcommand{\qbh}{\widehat{\qb}}
\newcommand{\rbh}{\widehat{\rb}}
\newcommand{\rnh}{\widehat{r}^n}
\newcommand{\Ubh}{\widehat{\Ub}}
\newcommand{\vbh}{\widehat{\vb}}
\newcommand{\vbht}{{\widehat{\vb}^t}}
\newcommand{\vbhU}{\underline{\vbh}}
\newcommand{\wbh}{\widehat{\wb}}
\newcommand{\zbh}{\widehat{\zb}}
\newcommand{\zbht}{\widehat{\zb}^t}

\newcommand{\mc}[1]{\mathcal{#1}}
\newcommand{\mcC}[1]{\mathcal{#1}}
\newcommand{\mb}[1]{\boldsymbol{#1}}
\newcommand{\mbb}[1]{\mathbb{#1}}
\newcommand{\mi}[1]{\mathit{#1}}
\newcommand{\mf}[1]{\mathfrak{#1}}
\newcommand{\Oi}[1]{\int_{\mc{D}} #1 \, d\vec{x}}
\newcommand{\TOi}[1]{\int_{0}^T \int_{\mc{D}} #1 \, d\vec{x}dt}
\newcommand{\Ti}[1]{\int_{0}^T #1 \, dt}
\newcommand{\TDei}[1]{\int_{0}^T\int_{\De} #1 \, d\vec{x}dt}
\newcommand{\TDeid}[1]{\int_{0}^T\int_{\De,N} #1 \, d\vec{x}dt}
\newcommand{\Dei}[1]{\int_{\De} #1 \, d\vec{x}}
\newcommand{\DeI}[1]{\int_{D} #1 \, d\vec{x}}
\newcommand{\DeiM}[1]{\int_{\Dhat} #1 \, d\vec{r}}
\newcommand{\DeMd}[1]{\int_{\De,N_e} #1 \, d\vec{x}}
\newcommand{\DeiMd}[1]{\int_{\Dhat,N_e} #1 \, d\vec{r}}
\newcommand{\DeMdN}[1]{\int_{D,N} #1 \, d\vec{x}}
\newcommand{\DeiMdN}[1]{\int_{\Dhat,N} #1 \, d\vec{r}}
\newcommand{\DeiMdNC}[2]{\int_{\Dhat,N_{#2}} #1 \, d\vec{r}}
\newcommand{\pDei}[1]{\int_{\partial \De} #1 \, d\vec{x}}
\newcommand{\pDeiM}[1]{\int_{\partial \Dhat} #1 \, d\vec{r}}
\newcommand{\pDeiMd}[1]{\int_{\partial \Dhat,N_e} #1 \, d\vec{r}}
\newcommand{\pDeiMdN}[1]{\int_{\partial \Dhat,N} #1 \, d\vec{r}}
\newcommand{\pDeiMdNC}[2]{\int_{\partial \Dhat,N_{#2}} #1 \, d\vec{r}}
\newcommand{\pMortar}[2]{\int_{\mc{M}_{#2}} #1 \, d\vec{r}}
\newcommand{\pMortard}[2]{\int_{\mc{M}^{#2},N_{#2}} #1 \, d\vec{r}}
\newcommand{\pMortardd}[3]{\int_{\mc{M}^{#2},N_{#3}} #1 \, d\vec{r}}
\newcommand{\pMortardh}[3]{\int_{\mc{M}_{#2},N_{#3}} #1 \, d\vec{r}}
\newcommand{\pMortardhn}[4]{\int_{{#2}_{#3},N_{#4}} #1 \, d\vec{r}}
\newcommand{\TpDei}[1]{\int_{0}^T\int_{\partial \De} #1 \, d\vec{x}}
\newcommand{\TpDeid}[1]{\int_{0}^T\int_{\partial \De,N} #1 \, d\vec{x}}
\newcommand{\Deid}[1]{\int_{\De,N_e} #1 \, d\vec{x}}
\newcommand{\DeidN}[1]{\int_{\De,N} #1 \, d\vec{x}}
\newcommand{\pDeidN}[1]{\int_{\partial \De,N} #1 \, d\vec{x}}
\newcommand{\pDeid}[1]{\int_{\partial \De,N_e} #1 \, d\vec{x}}
\newcommand{\nor}[1]{\left\| #1 \right\|}
\newcommand{\snor}[1]{\left| #1 \right|}
\newcommand{\LRp}[1]{\left( #1 \right)}
\newcommand{\LRs}[1]{\left[ #1 \right]}
\newcommand{\LRa}[1]{\left< #1 \right>}
\newcommand{\LRc}[1]{\left\{ #1 \right\}}
\newcommand{\pp}[2]{\frac{\partial #1}{\partial #2}}
\newcommand{\ppcp}[1]{\frac{\partial #1}{\partial \vec{c_p}}}
\newcommand{\ppcpj}[1]{\frac{\partial #1}{\partial c_p}}
\newcommand{\ppcpmf}[1]{\frac{\partial #1}{\partial \vec{\mf{c_p}}}}
\newcommand{\ppm}[1]{\frac{\partial #1}{\partial \vec{m}}}
\newcommand{\ppho}[3]{\frac{\partial^{#3} #1}{{\partial #2}^{#3}}}
\newcommand{\ppmi}[1]{\frac{\partial #1}{\partial m_i}}
\newcommand{\ppmj}[1]{\frac{\partial #1}{\partial m_j}}
\newcommand{\ppcpm}[1]{\frac{\partial #1}{\partial \vec{c_p}^-}}
\newcommand{\ppcpp}[1]{\frac{\partial #1}{\partial \vec{c_p}^+}}
\newcommand{\ppcs}[1]{\frac{\partial #1}{\partial \vec{c_s}}}
\newcommand{\ppcsm}[1]{\frac{\partial #1}{\partial \vec{c_s}^-}}
\newcommand{\ppcsp}[1]{\frac{\partial #1}{\partial \vec{c_s}^+}}
\newcommand{\td}[2]{\frac{d #1}{d #2}}
\newcommand{\Rvert}[1]{\left. #1 \right|}
\newcommand{\bs}[1]{\mb{#1}}
\newcommand{\ipDed}[3]{\LRp{ #1, #2}_{\De,N}^{#3}}

\newcommand{\bigO}{\mathcal{O}}

\newcommand{\iOm}[1]{\int_{\Omega} #1 \, d\Omega}
\newcommand{\iGb}[1]{\int_{\Gamma_{b}} #1 \, ds}
\newcommand{\iGs}[1]{\int_{\Gamma_{s}} #1 \, ds}
\newcommand{\iGsy}[1]{\int_{\Gamma_{s}} #1 \, ds(\mb{y})}
\newcommand{\RE}{\mbox{{I}}\hspace{-0.5ex}\mbox{{R}}}
\newcommand{\iC}[1]{\int_{0}^{2\pi} #1 \, d\theta}
\newcommand{\iGr}[1]{\int_{\Gamma_{\infty}} #1 \, ds}

\newcommand{\jump}[1] {\ensuremath{[\![{#1}]\!]}}
\newcommand{\diff}[1] {\ensuremath{\LRs{#1}}}
\newcommand{\average}[1] {\ensuremath{\LRc{\!\!\LRc{#1}\!\!}}}

\newcounter{tablecell}

\newcommand*{\numbercell}{%
  \stepcounter{tablecell}%
  \thetablecell. 
}

\newtheorem{propo}{Proposition}
\newtheorem{coro}{Corollary}
\newtheorem{theo}{Theorem}
\newtheorem{lem}{Lemma}

\newtheorem{defi}{definition}

\newtheorem{rema}{remark}

\newcommand{\figlab}[1]{\label{fig:#1}}
\newcommand{\eqnlab}[1]{\label{eq:#1}}
\newcommand{\theolab}[1]{\label{theo:#1}}
\newcommand{\corolab}[1]{\label{coro:#1}}
\newcommand{\propolab}[1]{\label{propo:#1}}
\newcommand{\lemlab}[1]{\label{lem:#1}}
\newcommand{\defilab}[1]{\label{defi:#1}}
\newcommand{\remalab}[1]{\label{rema:#1}}
\newcommand{\tablab}[1]{\label{tab:#1}}

\newcommand{\figref}[1]{\ref{fig:#1}}
\newcommand{\theoref}[1]{\ref{theo:#1}}
\newcommand{\defiref}[1]{\ref{defi:#1}}
\newcommand{\remaref}[1]{\ref{rema:#1}}
\newcommand{\cororef}[1]{\ref{coro:#1}}
\newcommand{\proporef}[1]{\ref{propo:#1}}
\newcommand{\lemref}[1]{\ref{lem:#1}}
\newcommand{\eqnref}[1]{\eqref{eq:#1}}
\newcommand{\alglab}[1]{\label{alg:#1}}
\newcommand{\algref}[1]{\ref{alg:#1}}
\newcommand{\seclab}[1]{\label{sect:#1}}
\newcommand{\secref}[1]{\ref{sect:#1}}
\newcommand{\tabref}[1]{\ref{tab:#1}}

\renewcommand{\u}{u}
\newcommand{\uk}{\u^k}
\newcommand{\upk}{\u^{+k}}
\newcommand{\ukp}{\u^{k+1}}
\newcommand{\umkp}{\u^{-k+1}}
\newcommand{\umk}{\u^{-k}}
\renewcommand{\v}{v}
\newcommand{\vk}{{\v}^k}
\newcommand{\vkp}{{\v}^{k+1}}
\newcommand{\un}{\u_h}
\renewcommand{\e}{e}
\newcommand{\w}{w}
\renewcommand{\wb}{\mb{\w}}
\newcommand{\J}{J}
\newcommand{\Jb}{\mb{J}}
\newcommand{\q}{q}
\renewcommand{\r}{r}
\newcommand{\qh}{\hat{\q}}
\newcommand{\qs}{\q^*}
\newcommand{\qht}{\tilde{\qh}}
\newcommand{\qbar}{\overline{\q}}
\newcommand{\qb}{{\mb{\q}}}
\newcommand{\sig}{{\sigma}}
\newcommand{\thetah}{{\hat{\theta}}}
\newcommand{\thetas}{{\theta^*}}
\newcommand{\thetahn}{{\thetah_h}}
\newcommand{\sign}[1]{\text{ sgn}\!\LRp{#1}}
\newcommand{\sigh}{\hat{\sig}}
\newcommand{\sigs}{\sig^*}
\newcommand{\sigk}{\sig^{k}}
\newcommand{\sigpk}{\sig^{+k}}
\newcommand{\sigkp}{\sig^{k+1}}
\newcommand{\sigmkp}{\sig^{-k+1}}
\newcommand{\sigmk}{\sig^{-k}}
\newcommand{\sigb}{{\bs{\sig}}}
\newcommand{\sigbk}{\sigb^k}
\newcommand{\sigbkp}{\sigb^{k+1}}
\newcommand{\sigbn}{\sigb_h}
\newcommand{\sigbs}{\sigb^*}
\newcommand{\taub}{{\bs{\tau}}}
\newcommand{\taustar}{\tau_*}
\newcommand{\kappab}{{\bs{\kappa}}}
\newcommand{\gamb}{{\bs{\gamma}}}
\newcommand{\ut}{\tilde{\u}}
\newcommand{\sigbt}{\tilde{\sigb}}
\newcommand{\vt}{\tilde{\v}}
\newcommand{\taubt}{\tilde{\taub}}
\newcommand{\ub}{\mb{\u}}
\newcommand{\ubl}{\mb{\u}_{\lambda}}
\newcommand{\ubf}{\mb{\u}_{f}}
\newcommand{\ubp}{\mb{\u}^+}
\newcommand{\ubm}{\mb{\u}^-}
\newcommand{\vb}{\mb{\v}}
\newcommand{\um}{\u^-}
\newcommand{\up}{\u^+}
\newcommand{\ubk}{\ub^k}
\newcommand{\ubkp}{\ub^{k+1}}

\newcommand{\m}{{m}}

\newcommand{\ud}{{\dot{\u}}}
\newcommand{\vd}{{\dot{\v}}}
\newcommand{\sigbd}{{\dot{\sigb}}}
\newcommand{\taubd}{\dot{\taub}}

\renewcommand{\d}{{d}}
\newcommand{\R}{{\mathbb{R}}}
\newcommand{\kap}{\kappa}
\newcommand{\F}{\mb{F}}
\newcommand{\f}{f}
\newcommand{\g}{g}
\newcommand{\fb}{\mb{\f}}
\newcommand{\gb}{\mb{\g}}
\newcommand{\gD}{g_D}
\newcommand{\pOmega}{\partial \Omega}
\newcommand{\Omegah}{{\Omega_h}}
\newcommand{\pOmegah}{{\partial \Omegah}}
\newcommand{\K}{T}
\newcommand{\Kp}{\K^+}
\newcommand{\Km}{\K^-}
\newcommand{\pK}{{\partial \K}}
\newcommand{\pKin}{{\pK^{\text{in}}}}
\newcommand{\pKout}{{\pK^{\text{out}}}}
\newcommand{\Kj}{{\K_j}}
\newcommand{\Gh}{{\mc{E}_h}}
\newcommand{\Gho}{{\mc{E}_h^o}}
\newcommand{\Ghb}{{\mc{E}_h^{\partial}}}
\newcommand{\Poly}{{\mc{P}}}
\newcommand{\D}{{D}}
\newcommand{\p}{{p}}
\newcommand{\pres}{{q}}
\newcommand{\presl}{{q}_{\lambda}}
\newcommand{\presf}{{q}_{f}}
\newcommand{\Ts}{{\mb{T}}}
\newcommand{\Vb}{{\mb{V}}}
\newcommand{\V}{{V}}
\newcommand{\Lam}{{\Lambda}}
\newcommand{\Lamb}{\bs{\Lambda}}
\newcommand{\mub}{\bs{\mu}}
\newcommand{\lambdab}{\bs{\lambda}}
\newcommand{\Thn}{\mathcal{T}_h}
\newcommand{\Th}{{\Ts_h}}
\newcommand{\Vh}{{\V_h}}
\newcommand{\Vbh}{{\Vb_h}}
\newcommand{\Lamh}{{\Lam_h}}
\newcommand{\Lambh}{{\Lamb_h}}
\newcommand{\Lamho}{{\overline{\Lam}_h}}
\newcommand{\Lambht}{{\Lamb_h^t}}
\newcommand{\Ltwo}{{L^2\LRp{\Omega}}}
\newcommand{\Ltwoe}{{L^2\LRp{\e}}}
\newcommand{\Ltw}{{L^2\LRp{\K}}}
\newcommand{\VhK}{{\V_h\LRp{\K}}}
\newcommand{\ThK}{{\Ts_h\LRp{\K}}}
\newcommand{\VbhK}{{\Vb_h\LRp{\K}}}
\newcommand{\Lamhe}{{\Lam_h\LRp{\e}}}
\newcommand{\Lambhe}{{\Lamb_h\LRp{\e}}}
\renewcommand{\L}{\mc{L}}
\renewcommand{\D}{\mc{D}}
\newcommand{\T}{\mc{T}}
\newcommand{\Tt}{\tilde{\T}}
\newcommand{\A}{\mb{A}}
\newcommand{\B}{\mb{B}}
\renewcommand{\D}{\mb{D}}
\newcommand{\Am}{\A^-}
\newcommand{\Ap}{\A^+}
\newcommand{\Aa}{\snor{\A}}
\newcommand{\Rb}{\mb{R}}
\newcommand{\C}{\mb{C}}
\renewcommand{\S}{\mb{S}}
\newcommand{\Sm}{\S^{-}}
\newcommand{\Sp}{\S^{+}}
\newcommand{\Spm}{\S^{\pm}}

\newcommand{\Lte}{{L^2\LRp{\Gh}}}
\newcommand{\n}{\mb{n}}
\newcommand{\nm}{\n^-}
\newcommand{\np}{\n^+}
\newcommand{\uh}{{\hat{\u}}}
\newcommand{\uhk}{{\uh}^k}
\newcommand{\uhkp}{{\uh}^{k+1}}
\newcommand{\us}{{\u^*}}
\newcommand{\uhn}{{\uh_h}}
\newcommand{\ubh}{{\hat{\ub}}}
\newcommand{\ubs}{{\ub^*}}
\newcommand{\ubht}{{\tilde{\ubh}}}
\newcommand{\uht}{{\tilde{\uh}}}
\newcommand{\uhtn}{{\uht_h}}
\newcommand{\uhd}{{\dot{\uh}}}
\newcommand{\ubhk}{\ubh^k}
\newcommand{\ubhkp}{\ubh^{k+1}}
\newcommand{\sigbh}{{\hat{\sigb}}}
\newcommand{\Fh}{{\hat{\F}}}
\newcommand{\Fs}{\F^*}

\newcommand{\zb}{\mb{z}}
\renewcommand{\sb}{\mb{s}}
\newcommand{\rb}{\mb{r}}
\newcommand{\betab}{\bs{\beta}}
\newcommand{\veps}{{\varepsilon}}
\newcommand{\vepsb}{\bs{\veps}}
\newcommand{\Hb}{\mb{H}}
\newcommand{\Hbt}{\Hb^t}
\newcommand\E{\mathcal{E}}
\newcommand\Eh{\mathcal{E}_h}
\newcommand\Ehi{\mathcal{E}_{h,i}}
\newcommand\Ehiint{\mathcal{E}_{h,i}^\mathrm{int}}
\newcommand\EhG{\mathcal{E}_{h}^{\Gamma}}
\newcommand\EhiG{\mathcal{E}_{h,i}^{\Gamma}}
\newcommand\EhijG{\mathcal{E}_{h,i,j}^{\Gamma}}
\newcommand\EhjiG{\mathcal{E}_{h,j,i}^{\Gamma}}
\newcommand{\Eb}{\mb{E}}
\newcommand{\Ebt}{\Eb^t}
\newcommand{\hb}{\mb{h}}
\newcommand{\eb}{\mb{e}}
\newcommand{\Hbh}{\hat{\mb{H}}}
\newcommand{\Ebh}{\hat{\mb{E}}}
\newcommand{\Hbs}{{\mb{H}^*}}
\newcommand{\Ebs}{{\mb{E}^*}}
\newcommand{\Ebht}{\Ebh^t}
\newcommand{\Hbht}{\Hbh^t}
\newcommand{\Lb}{\mb{L}}
\newcommand{\Gb}{\mb{G}}
\newcommand{\Lbh}{\hat{\Lb}}
\newcommand{\Lbs}{\Lb^*}
\newcommand{\Ib}{\mb{I}}
\newcommand{\I}{I}
\newcommand{\Q}{Q}
\newcommand{\Sigb}{\bs{\Sigma}}
\newcommand{\Piu}{\Pi_\u}
\newcommand{\Pisigb}{\Pi_\sigb}
\newcommand{\Z}{Z}
\newcommand{\Y}{Y}
\newcommand{\rhou}{\rho\u}
\newcommand{\rhov}{\rho\v}
\renewcommand{\P}{P}
\newcommand{\h}{H}

\newcommand{\mygreen}[1]{{\color[rgb]{0,.65,0} #1}}
\newcommand{\mywhite}[1]{{\color[rgb]{1.0,1.0,1.0} #1}}
\newcommand{\myred}[1]{{\color[rgb]{0.65,0.0,0.0} #1}}
\newcommand{\myblue}[1]{{\color[rgb]{0.0,0.0,1.0} #1}}
\newcommand{\mygray}[1]{{\color[rgb]{.15,.35,.60} #1}}
\newcommand{\mythemec}[1]{{\color[RGB]{51,108,121} #1}}
\newcommand{\mydarkred}[1]{{\color[rgb]{.6,.1,.1} #1}}
\newcommand{\mydarkblue}[1]{{\color[rgb]{.1,.1,.9} #1}}
\newcommand{\mygreenback}[1]{{\color[rgb]{.75,1,.75} #1}}
\newcommand{\myorange}[1]{{\color[rgb]{.76,.39,.13} #1}}
\newcommand{\mygrass}[1]{{\color[rgb]{.19,.64,.13} #1}}
\newcommand{\mysierp}[1]{{\color[RGB]{209,28,209} #1}}

\newcommand{\tanbui}[2]{\textcolor{blue}{#1} \textcolor{red}{#2}}

\newcommand{\vel}{\bs{\vartheta}}
\newcommand{\vels}{\vel^*}
\newcommand{\vele}{\vel^e}
\newcommand{\velh}{\hat{\vel}}
\newcommand{\velk}{\vel^k}
\newcommand{\velkp}{\vel^{k+1}}
\newcommand{\vhk}{\hat\v^{k}}

\newcommand{\phis}{\phi^*}
\newcommand{\phie}{\phi^e}
\newcommand{\phih}{\hat{\phi}}
\newcommand{\phihk}{\phih^k}
\newcommand{\phik}{\phi^k}
\newcommand{\phikp}{\phi^{k+1}}

\newcommand{\algn}[1]{\begin{align} #1 \end{align}}
\newcommand{\algns}[1]{\begin{align*}#1\end{align*}}
\newcommand{\flalgns}[1]{\begin{flalign*}#1\end{flalign*}}
\newcommand{\subeqn}[1]{\begin{subequations} #1 \end{subequations}}
\newcommand{\mltln}[1]{\begin{multline} #1 \end{multline}}
\newcommand{\mltlns}[1]{\begin{multline*} #1 \end{multline*}}

%
\newcommand{\AB}{\mb{A}}
\newcommand{\BB}{\mb{B}}
\newcommand{\BBT}{\mb{B}^\top}
\newcommand{\CB}{\mb{C}}
\newcommand{\DB}{\mb{D}}
\newcommand{\FB}{\mb{F}}
\newcommand{\GB}{\mb{G}}
\newcommand{\IB}{\mb{I}}
\newcommand{\LB}{\mb{L}}
\newcommand{\NB}{\mb{N}}
\newcommand{\RB}{\mb{R}}
\newcommand{\SB}{\mb{S}}
\newcommand{\TB}{\mb{T}}
\newcommand{\uB}{\mb{u}}
\newcommand{\zeroB}{\mb{0}}
\newcommand{\zerob}{\mb{0}}
\newcommand{\db}{\mb{d}}

%
\newcommand{\ChM}{\widehat{C}_h}
\newcommand{\ChMD}{\widehat{C}_h^D}
\newcommand{\CbhM}{\widehat{\mb{C}}_h}
\newcommand{\CbhMi}{\widehat{\mb{C}}^i_h}
\newcommand{\CbhtM}{\widehat{\mb{C}}^t_h}
\newcommand{\CbhtMi}{\widehat{\mb{C}}^{t,i}_h}
\newcommand{\CbhtMD}{\widehat{\mb{C}}^{t,D}_h}
\newcommand{\FhM}{\widehat{F}_h}
\newcommand{\FhMi}{\widehat{F}_h^i}
\newcommand{\FhMN}{\widehat{F}_h^N}
\newcommand{\GBM}{\mb{G}_h}
\newcommand{\VbM}{\mb{V}_h}
\newcommand{\VBM}{\mb{V}_h}
\newcommand{\JBM}{\mb{H}_h}
\newcommand{\CBM}{\mb{C}_h}
\newcommand{\SBM}{\mb{S}_h}
\newcommand{\VhM}{\widehat{V}_h}
\newcommand{\VbhM}{\widehat{\mb{V}}_h}
\newcommand{\VbhMi}{\widehat{\mb{V}}^i_h}
\newcommand{\VbhtM}{\widehat{\mb{V}}^t_h}
\newcommand{\VbhtMi}{\widehat{\mb{V}}^{t,i}_h}
\newcommand{\QM}{Q_h}
\newcommand{\WM}{W_h}
\newcommand{\QBM}{\mb{Q}_h}
\newcommand{\HbM}{\mb{H}_h}
\newcommand{\CbM}{\mb{C}_h}
\newcommand{\SM}{S_h}
\newcommand{\ShM}{\widehat{S}_h}
\newcommand{\ShMi}{\widehat{S}_h^i}
\newcommand{\ShMD}{\widehat{S}_h^D}
\newcommand{\UBM}{\mb{U}_h}
\newcommand{\MuBM}{\mb{M}_h}
\newcommand{\MuBMo}{\mb{M}^o_h}
\newcommand{\LamBM}{\mb{\Lambda}^t_h}
\newcommand{\GamBM}{\mb{\Gamma}_h}
\newcommand{\UBhM}{\mb{\widehat{U}}_h}
\newcommand{\XBM}{\mb{X}_h}
\newcommand{\PkKscalar}{{\Poly_k(K)}}
\newcommand{\PkKvector}{{\LRs{\Poly_k(K)}^d}}
\newcommand{\PkKtensor}{{\LRs{\Poly_k(K)}^{d \times d}}}
\newcommand{\PkThscalar}{{\Poly_k(\K)}}
\newcommand{\PzpThscalar}{{\Poly_0(\pK)}}
\newcommand{\PkThvector}{{\LRs{\Poly_k(\K)}^d}}
\newcommand{\PkThtensor}{{\LRs{\Poly_k(\K)}^{d \times d}}}
\newcommand{\PkEhscalar}{{\Poly_k(\Gh)}}
\newcommand{\PkEhvector}{{\LRs{\Poly_k(\Gh)}^d}}
\newcommand{\PkEhtensor}{{\LRs{\Poly_k(\Gh)}^{d \times d}}}
\newcommand{\PkEhintscalar}{{\Poly_k(\Gho)}}
\newcommand{\PkEhintvector}{{\LRs{\Poly_k(\Gho)}^d}}
\newcommand{\PkEhintvectorD}{{\LRs{\Poly_k(\Gho\cup\GhD)}^d}}
\newcommand{\PkEhintvectorN}{{\LRs{\Poly_k(\Gho\cup\GhN)}^d}}
\newcommand{\PkEhinttensor}{{\LRs{\Poly_k(\Gho)}^{d \times d}}}
\newcommand{\PkEhintscalarD}{{\Poly_k(\Gho\cup\GhD)}}
\newcommand{\PkEhintscalarN}{{\Poly_k(\Gho\cup\GhN)}}
\newcommand{\PkpOmegaDscalar}{{\Poly_k(\pOmegaD)}}
\newcommand{\PkpOmegaNvector}{{{\Poly_k(\pOmegaN)}^d}}
\newcommand{\PkpOmegaDvector}{{{\Poly_k(\pOmegaD)}^d}}
\newcommand{\PkpOmegaNscalar}{{\Poly_k(\pOmegaN)}}

\newcommand{\wbar}{\overline{w}}
\newcommand{\presbar}{\overline{\pres}}

\newcommand{\taun}{\tau_n}
\newcommand{\taut}{\tau_t}
\newcommand{\betan}{\beta_n}
\newcommand{\betat}{\beta_t}
\newcommand{\nb}{\bs{n}}
\newcommand{\cb}{\bs{c}}
\newcommand{\tb}{\bs{t}}

\maketitle

\begin{abstract}
    We present a scalable block preconditioning strategy for the trace system coming from the high-order hybridized discontinuous Galerkin (HDG) discretization of incompressible resistive magnetohydrodynamics (MHD). We construct the block preconditioner with a least squares commutator (BFBT) approximation for the inverse of the Schur complement that segregates out the pressure unknowns of the trace system. The remaining velocity, magnetic field, and Lagrange multiplier unknowns form a coupled
    nodal unknown block (the upper block), for which a system algebraic multigrid (AMG) is used for the approximate inverse. The complexity of the MHD equations together with the algebraic nature of the statically condensed HDG trace system makes the choice of smoother in the system AMG part critical for the convergence and performance of the block preconditioner. Our numerical experiments show GMRES preconditioned by ILU(0) of overlap zero as a smoother inside system AMG performs best in terms of robustness, time per nonlinear iteration and memory requirements. With several transient test cases in 2D and 3D including the island coalescence problem at high Lundquist number we demonstrate the robustness and parallel scalability of the block preconditioner. Additionally for the upper block a preliminary study of an alternate nodal block system solver based on a multilevel approximate nested dissection is presented. On a 2D island coalescence problem the multilevel approximate nested dissection preconditioner shows better scalability with respect to mesh refinement than the system AMG, but is relatively less robust with respect to Lundquist number scaling. 
\end{abstract}

\pagestyle{myheadings} \thispagestyle{plain} \markboth{S.~Muralikrishnan, S.~Shannon, T.~Bui-Thanh, J.~N.~Shadid}{Block preconditioner for HDG applied to incompressible resistive MHD}





\section{Introduction}

Incompressible visco-resistive magnetohydrodynamics (MHD) equations play an important
role in modeling low Lundquist number liquid metal flows, high Lundquist number large-guide-field fusion plasmas and 
low flow-Mach-number compressible flows \cite{goedbloed2004principles,chacon2008optimal,shadid2016scalable}. 
The mathematical basis
for the continuum modeling of these systems is the solution of the governing partial differential equations (PDEs)
describing conservation of mass, momentum, and energy, augmented by the low-frequency Maxwell's equations.
The multiphysics phenomena produced by this strongly coupled nonlinear system includes, convective transport, momentum forces
evolving from pressure gradients and Lorentz forces from the magnetic field, Alfven wave propagation, and diffusion effects from the viscosity and resistivity of the plasma system. Additionally the coupling of the conservation equations to elliptic constraints from  the incompressiblity assumption and the solenoidal involution property of the magnetic field add significant complexity to the block structured interactions of the unknowns. Finally the resulting response of these systems
is also characterized by a wide range of spatial and temporal scales which makes efficient 
solution of these systems extremely challenging. In the context of time integration of these systems
some form of implicitness is required to deal with the infinite pressure wave speed and to remove one or more sources of stiffness arising from diffusion and fast-waves in the MHD system. 
The potential benefit of implicit methods however require 
robust, efficient and scalable nonlinear and linear iterative solvers/preconditioners which can handle the highly stiff algebraic systems generated from
incompressible MHD equations.

Over the past few years there has been significant progress in the area of fully implicit robust and scalable MHD simulations as evidenced in 
\cite{chacon2008optimal,cyr2013new,shadid2016scalable,phillips2016block,wathen2018preconditioners} among others. However, these simulations use low-order 
stabilized FEM, mixed FEM or finite volume methods for the spatial 
discretization and hence require larger numbers of unknowns for high
accuracy. They also have low computation to communication ratio which
is not preferable for achieving high performance (in terms of percentage of peak performance) 
in modern extreme scale computing architectures.

High-order spatial discretizations are attractive in this context as
they have high computation to communication ratio and can simulate 
the problem with a reduced number of unknowns compared to low-order
discretizations for similar or higher accuracy. Hybridized discontinuous 
Galerkin (HDG) methods introduced a decade ago are promising 
candidates for high-order spatial discretization due to the following
reasons. It combines the important advantages of DG methods
namely arbitrary high-order with compact stencil, flux upwinding on element boundaries and the
ability to handle complex geometries together with the smaller and 
sparser linear system which involves only trace unknowns, a 
characteristic of hybridized methods.

The first HDG method for incompressible resistive MHD equations is 
proposed in \cite{shannonCSRIProceedings,lee2019analysis} and very recently HDG methods for
ideal and resistive compressible MHD equations are introduced in \cite{ciucua2020implicit}. A current challenge in the
context of high-order HDG discretizations of MHD is the availability of robust, efficient and scalable iterative solvers/preconditioners
and that is the goal of our work. The availability of such a solver will enable large scale 3D high-order HDG simulations of MHD on realistic geometries and this will have large impact in many areas of plasma physics research including but not limited to nuclear 
fusion.

Towards that extent, we propose a block preconditioning strategy for trace systems \cite{muralikrishnan2019fast} arising from HDG discretizations of incompressible 
resistive MHD \cite{lee2019analysis}. Even though the concept of block preconditioning for saddle point systems has been studied by
various researchers \cite{benzi2005numerical,elman2006block,elman1999preconditioning,planas2013stabilized,elman2014finite,rudi2019global,cyr2013new,phillips2016block,wathen2018preconditioners,wathen2017preconditioners,wathen2020scalable}, it is mostly in the 
context of volume based low-order discretizations such as mixed FEM, 
stabilized FEM, finite volume and finite difference. 

In this development the linear systems coming from HDG discretization of incompressible resistive MHD \cite{lee2019analysis} pose several challenges which are as follows. First, it involves only trace unknowns which live 
on the skeleton of the mesh and moreover they all are not of the same nature. 
Second, since it is obtained after static 
condensation of the volume unknowns the linear system is mostly 
algebraic and it is non-trivial to identify the nature of different
blocks in the linear system as opposed to linear systems coming from stabilized and mixed FEM. In \cite{rhebergen2018preconditioning,southworth2020fixed} block preconditioners for HDG discretization of
Stokes and incompressible Navier Stokes equations are presented. 
However, the authors eliminated only velocity volume unknowns by static condensation and the final linear system contains velocity and pressure trace unknowns togther with the pressure volume unknowns. \emph{To
the best of our knowledge we are not aware of block preconditioners 
for HDG discretizations which completely eliminate all the volume 
unknowns and solve the linear system with only trace unknowns and in
that aspect our work is a first step towards it.}

In this effort we employ the BFBT framework of block preconditioning first introduced in \cite{elman1999preconditioning} for finite difference and finite element discretizations of incompressible Navier-Stokes equations. 
However, the novelty of our work lies in the careful choice of 
different components involved in the framework so that it leads to 
an effective preconditioner for HDG discretization of incompressible
resistive MHD. We demonstrate this numerically by means of several 2D
and 3D benchmark problems and a theoretical analysis of the current 
approach is left for future work. 

This paper is organized as follows. In section \secref{notation}, 
we introduce the notations which will be followed throughout the rest 
of the paper. We present the incompressible resistive MHD equations and the relevant non-dimensional parameters in 
section \secref{MHD_equations}. Then in section \secref{HDG_MHD}, we present an HDG scheme for the discretization of the MHD system and identify the block structure in it. We then proceed to introduce a block preconditioning strategy for the linear system in section \secref{block}. Section \secref{numerical_mhd} presents various 2D and 3D transient test cases to test the robustness and scalability of the
block preconditioner. Finally in section \secref{conclusions} we discuss the conclusions and directions for future research.

\section{Notation}
\seclab{notation}

In this section we introduce the common notations and conventions which will be followed in the rest of the paper. Let us partition $\Omega \in \R^\d$, an  open and bounded domain,  into
$\Nel$ nonoverlapping elements $\Kj, j = 1, \hdots, \Nel$ with
Lipschitz boundaries such that $\Omega_h := \cup_{j=1}^\Nel \Kj$ and
$\overline{\Omega} = \overline{\Omega}_h$. Here, $h$ is defined as $h
:= \max_{j\in \LRc{1,\hdots,\Nel}}diam\LRp{\Kj}$. We denote the
skeleton of the mesh by $\Gh := \cup_{j=1}^\Nel \partial \K_j$,
the set of all (uniquely defined) interfaces $\e$ between elements. We conventionally identify $\nm$ as the outward 
normal vector on the boundary $\pK$ of element $\K$ (also denoted as $\Km$) and $\np = -\nm$ as the outward normal vector of the boundary of a neighboring element (also denoted as $\Kp$). Furthermore, we use $\n$ to denote either $\nm$ or $\np$ in an expression that is valid for both cases, and this convention is also used for other quantities (restricted) on
a face $\e \in \Gh$.
   For the sake of convenience, we denote by $\Ghb$ the set of all boundary faces on $\pOmega$, by $\Gho := \Gh \setminus \Ghb$ the set of all interior faces, and $\pOmega_h := \LRc{\pK:\K \in \Omega_h}$.

For simplicity in writing we define $\LRp{\cdot,\cdot}_\K$ as the
$L^2$-inner product on a domain $\K \in \R^\d$ and
$\LRa{\cdot,\cdot}_\K$ as the $L^2$-inner product on a domain $\K$ if
$\K \in \R^{\d-1}$. We shall use $\nor{\cdot}_{\K} :=
\nor{\cdot}_{\Ltw}$ as the induced norm for both cases and the
particular value of $\K$ in a context will indicate which inner
product the norm is coming from. We also denote the $\veps$-weighted
norm of a function $\u$ as $\nor{\u}_{\veps, \K} :=
\nor{\sqrt{\veps}\u}_{\K}$ for any positive $\veps$. Moreover, we define
$\LRp{\u,\v}_{\Omega_h} := \sum_{\K\in \Omega_h}\LRp{\u,\v}_\K$ and
$\LRa{\u,\v}_\Gh := \sum_{\e\in \Gh}\LRa{\u,\v}_\e$ whose
induced (weighted) norms are clear, and hence their definitions are
omitted. We shall use
boldface lowercase letters, e.g. $\ub$, for vector-valued functions and boldface uppercase letters, e.g. $\mb{L}$, for matrices and tensors. The notations of inner products and norms are naturally 
extended for these cases in a component-wise manner. We use the terms ``skeletal unknowns'' and ``trace unknowns'' interchangeably and they both refer to the unknowns on the mesh skeleton.\ They are denoted 
with hats to differentiate it from the corresponding volume unknowns.

We define the gradient of a vector, the divergence of a matrix, and the outer product symbol $\otimes$ as:
\[
  \LRp{\nabla \ub}_{ij} = \pp{u_i}{x_j}, \quad
  \LRp{\nabla \cdot \Lb}_i = \nabla \cdot \Lb\LRp{i,:} = \sum_{j=1}^3\pp{\bs{L}_{ij}}{x_j}, \quad
  \LRp{\bs{a} \otimes \bs{b}}_{ij} = a_i b_j = \LRp{\bs{a}\bs{b}^T}_{ij}.
\]

We define $\Poly^\p\LRp{\K}$ as the space of polynomials of degree at
most $\p$ on a domain $\K$. Similarly, $\Poly^\p\LRp{\e}$ denotes the
polynomials of degree at most $\p$ on the mesh skeleton edge $\e$ and the extensions to vector- or matrix-valued polynomials $\LRs{\Poly^\p(\K)}^d$,
$\LRs{\Poly^\p(\K)}^{d\times d}$, $\LRs{\Poly^\p(\e)}^d$, etc, are straightforward. We define the ``jump'' operator for any quantity
 $\LRp{\cdot}$ as $\jump{\LRp{\cdot}} := \LRp{\cdot}^- +
 \LRp{\cdot}^+$. 

\section{Incompressible visco-resistive MHD system}
\seclab{MHD_equations}
The visco-resistive, incompressible MHD equations in non-dimensional form are given by
\subeqn{
  \label{eq:MHDnondim}
  \begin{alignat}{2}
    \label{eq:MHDnondim1}
    \frac{\partial \ub}{\partial t} + \ub \cdot \nabla \ub + \nabla \pres 
    - \frac{1}{\Rey} \Delta \ub - \kappa (\nabla \times \bb) \times \bb
    &= \fb, \qquad &&  \\
    \label{eq:MHDnondim2}
    \nabla \cdot \ub &= 0, \quad &&  \\
    \label{eq:MHDnondim3}
    \kappa \frac{\partial \bb}{\partial t} 
    + \frac{\kappa}{\Rm} \nabla \times (\nabla \times \bb) 
    - \kappa \nabla \times (\ub \times \bb) 
    + \nabla \r 
    &= \gb, &&  \\
    \label{eq:MHDnondim4}
    \nabla \cdot \bb &= 0. && 
  \end{alignat}
}
Here, $\ub$ is the fluid velocity field, $\bb$ is the magnetic field, $\pres$ is the 
mechanical pressure and $\r$ denotes the Lagrange 
multiplier which is used to enforce the solenoidality of the magnetic field constraint \cite{codina2006stabilized,shadid2016scalable}.
Thus allows the enforcement of the divergence-free constraint and at the same time it is just an auxiliary discrete variable because in essence we are solving for a variable that would be zero in the solution to the continuous problem. This 
can be seen by taking divergence of equation \eqref{eq:MHDnondim3} which gives $\Delta\r=0$ (assuming the magnetic source term $\gb$ to be divergence-free) and together with homogeneous Dirichlet boundary conditions 
gives $\r=0$. More details about this approach can be found in \cite{salah1999conservative,schotzau2004mixed,codina2006stabilized,shadid2016scalable,lee2019analysis}.

The non-dimensional parameters in \eqref{eq:MHDnondim} are as follows. $\Rey := \frac{\rho l_0 u_0}{\mu}$ is the fluid Reynolds number which measures the ratio of inertial forces to viscous forces, $\Rm := \frac{\mu_0 u_0 l_0}{\eta}$
is the magnetic Reynolds number which is the ratio of magnetic advection to magnetic diffusion, $\kappa := \frac{b_0^2}{\rho \mu_0 u_0^2}$ is the coupling
parameter and is the ratio of electromagnetic forces to inertial forces. Here, $\rho$ is the fluid density, $\mu$ is the viscosity, $\mu_0$ is the permeability of free space and $\eta$ the resistivity. The characteristic length, velocity and magnetic scales are given by  $l_0$, $u_0$ and $b_0$ respectively. The parameters $\kappa$, $\Rey$ and $\Rm$ are related by $\kappa = \frac{\Ha^2}{\Rey\Rm}$,
where $\Ha := \frac{b_0 l_0}{\sqrt{\mu \eta}}$ is the Hartmann number. We can also write $\kappa$ as $\kappa = \frac{u_A^2}{u_0^2}$,
where $u_A := \frac{b_0}{\sqrt{\rho \mu_0}}$ is the Alfv\'{e}n speed. If we choose the characteristic velocity as Alfv\'{e}n speed then the resulting magnetic Reynolds number is referred as Lundquist number $S$. We refer the readers to \cite{goedbloed2004principles,muller2013magnetofluiddynamics} for details on non-dimensional parameters in the MHD system.

We put \eqref{eq:MHDnondim} into first order form for discretizing with HDG and towards that end let us define the 
auxiliary variables $\LB$ and $\Jb$ which represents the velocity gradient and curl of magnetic field 
respectively. The first order system is given by
\subeqn{
  \label{eq:MHDnondim_1st_order}
  \begin{align}
    \label{eq:MHDnondim_1st_order_L}
    \Rey \LB - \nabla \ub &= \zerob, \\
    \label{eq:MHDnondim_1st_order_u}
    \frac{\partial \ub}{\partial t} + \nabla \cdot (\ub \otimes \ub) + \nabla \pres 
    - \nabla \cdot \LB - \kappa (\nabla \times \bb) \times \bb
    &= \fb, \\
    \label{eq:MHDnondim_1st_order_p}
    \nabla \cdot \ub &= 0, \\
    \label{eq:MHDnondim_1st_order_J}
    \frac{\Rm}{\kappa} \Jb - \nabla \times \bb &= \zerob, \\
    \label{eq:MHDnondim_1st_order_b}
    \kappa \frac{\partial \bb}{\partial t} + \nabla \r - \kappa \nabla \times (\ub \times \bb) 
    + \nabla \times \Jb &= \gb, \\
    \label{eq:MHDnondim_1st_order_r}
    \nabla \cdot \bb &= 0.
  \end{align}
}
We refer to $\Jb$ as the current density or simply the current and it should be understood in a non-dimensional sense
with the characteristic value defined by $J_0 = \frac{Rm}{\kappa} \frac{b_0}{\mu_0 l_0}$.

The MHD system \eqref{eq:MHDnondim_1st_order} is equipped with the following set of initial conditions
\begin{equation}
  \label{eq:ICs}
  \ub(t=0) = \ub_0, \qquad \bb(t=0) = \bb_0.
\end{equation}
We also need to specify boundary conditions for the fluid components, magnetic components and the Lagrange multiplier. 
Since it is not important for the current discussion we will defer this till section \secref{numerical_mhd} where
we specify these details for each numerical experiment separately.

We will use the Picard nonlinear solver in our study and hence we consider the HDG schemes posed on the (Picard) linearized version of \eqref{eq:MHDnondim_1st_order} about a prescribed velocity $\wb$ and a prescribed magnetic field $\db$
\cite{codina2006stabilized}:
\subeqn{
  \label{eq:MHDnondimlin}
  \begin{align}
    \label{eq:MHDnondimlin_L}
    \Rey \LB - \nabla \ub &= \zerob, \\
    \label{eq:MHDnondimlin_u}
    \frac{\partial \ub}{\partial t} + \nabla \cdot (\ub \otimes \wb) + \nabla \pres 
    - \nabla \cdot \LB - \kappa (\nabla \times \bb) \times \db
    &= \fb, \\
    \label{eq:MHDnondimlin_p}
    \nabla \cdot \ub &= 0, \\
    \label{eq:MHDnondimlin_J}
    \frac{\Rm}{\kappa} \Jb - \nabla \times \bb &= \zerob, \\
    \label{eq:MHDnondimlin_b}
    \kappa \frac{\partial \bb}{\partial t} + \nabla \r - \kappa \nabla \times (\ub \times \db) 
    + \nabla \times \Jb &= \gb, \\
    \label{eq:MHDnondimlin_r}
    \nabla \cdot \bb &= 0.
  \end{align}
}
Here $\wb$ is assumed to reside in $H(div,\Omega)$ and be divergence-free,
while $\db$ is assumed to reside in $H(div,\Omega) \cap H(curl,\Omega)$.

\section{HDG for incompressible MHD}
\seclab{HDG_MHD}

In this section we present the HDG scheme first proposed in \cite{lee2019analysis,stephen_thesis}.

  Find $( \LB, \ub, \pres, \Jb, \bb, \r, \ubh, \bbht, \rh, \rho )$ in 
  $\GBM \times \VbM \times \WM \times \HbM \times \CbM \times \SM \times 
    \VbhM \times \CbhtM \times \ShM \times \Poly^0(\pK)$
  such that
  the local equations
  \subeqn{
    \label{eq:MHD_HDG_2}
    \algn{
      \label{eq:MHD_HDG_2a}
      &\Rey \LRp{ \LB, \GB }_\K
      + \LRp{ \ub, \nabla \cdot \GB }_\K
      - \LRa{ \ubh, \GB \nb }_\pK
      = 
      0
      , \\
      \label{eq:MHD_HDG_2b}
      &\LRp{ \pp{ \ub }{ t }, \vb }_\K
      - \LRp{ \nabla \cdot \LB, \vb }_\K
      + \LRp{ \nabla \pres, \vb }_\K
      - \half \LRp{ \ub \otimes \wb, \nabla \vb }_\K
      + \half \LRp{ \nabla \ub, \vb \otimes \wb }_\K
      \notag \\
      & \qquad
      + \LRp{ \nabla \times \bb, \vb \times \kappa \db }_\K
      + \LRa{ \half \LRp{\wb \cdot \nb} \ubh + \SB_u \LRp{ \ub - \ubh }, \vb }_\pK
      \notag \\
      & \qquad
      - \LRs{1-\xi} \LRa{ \nb \times \LRp{ \bbt - \bbht }, \vb \times \kappa \db }_\pK
      = 
      \LRp{ \fb, \vb }_\K
      , \\
      \label{eq:MHD_HDG_2c}
      &- \LRp{ \ub, \nabla \w }_\K
      + \LRa{ \ubh \cdot \nb, \w - \wbar }_\pK
      + \LRa{ \presbar - \rho, \wbar }_\pK
      =
      0
      , \\
      \label{eq:MHD_HDG_2d}
      &\frac{\Rm}{\kappa} \LRp{ \Jb, \Hb }_\K
      - \LRp{ \bb, \nabla \times \Hb }_\K
      - \LRa{ \nb \times \bbht, \Hb }_\pK
      =
      0
      , \\
      \label{eq:MHD_HDG_2e}
      &\kappa \LRp{ \pp{ \bb }{ t }, \cb }_\K
      + \LRp{ \nabla \times \Jb, \cb }_\K
      - \LRp{ \r, \nabla \cdot \cb }_\K
      - \LRp{ \ub \times \kappa \db, \nabla \times \cb }_\K
      + \LRa{ \rh, \cb \cdot \nb }_\pK
      \notag \\
      & \qquad
      + \LRa{ \LRp{ \LRs{1-\xi} \ub + \xi \ubh } \times \kappa \db, \nb \times \cb }_\pK
      + \LRa{ \betat \LRp{ \bbt - \bbht }, \cb }_\pK
      =
      \LRp{ \gb, \cb }_\K
      , \\
      \label{eq:MHD_HDG_2f}
      &\LRp{ \nabla \cdot \bb, s }_\K 
      + \LRa{ \frac{1}{\betan} \LRp{ \r - \rh }, s }_\pK 
      =
      0
      , 
    }
    the conservation equations
    \algn{
      \label{eq:MHD_HDG_2g}
      &- \LRa{ \jump{-\LB \nb + \pres \nb + \half \LRp{ \wb \cdot \nb } \ub
      + \SB_u \LRp{ \ub - \ubh } + \kappa \db \times \LRp{ \nb \times \xi \bb }}, \vbh }_\e 
      = 
      0
      , \\ 
      \label{eq:MHD_HDG_2h}
      &- \LRa{ \jump{\nb \times \Jb + \betat \LRp{ \bbt - \bbht }
      - \nb \times \LRp{ \LRs{1-\xi} \ub \times \kappa \db }}, \cbht }_\e
      = 
      0
      , \\ 
      \label{eq:MHD_HDG_2i}
      &- \LRa{ \jump{\bb \cdot \nb + \frac{1}{\betan} \LRp{ \r - \rh }}, \sh }_\e
      = 
      0
      , 
    }
    and the additional constraint
    \algn{
      \label{eq:MHD_HDG_2j}
      &\LRa{ \ubh \cdot \nb, \psi }_\pK
      = 
      0
    }
  }
  hold
  for all $( \GB, \vb, \w, \Hb, \cb, s, \vbh, \cbht, \sh, \psi )$ in 
  $\GBM \times \VbM \times \WM \times \HbM \times \CbM \times \SM \times 
  \VbhM \times \CbhtM \times \ShM \times \Poly^0(\pK)$. In
  addition the pressure is subject to the constraint 
  \[
      (\pres,1)_\Omegah = 0.
  \]
In the above, the notation $\qbar$ is defined by $\qbar := \snor{ \pK }^{-1} \LRa{ q, 1 }_\pK$
as the $\pK$-wise average of $q$, and $\snor{ \pK }$ is the length of the perimeter of element $\K$.
The new unknowns $\rho$ which are sought in $\Poly^0(\pK)$ represent the $\pK$-wise or edge average pressure. Also, $\bbht$ represents the skeletal unknowns corresponding to the tangent magnetic field. 

Here, the volume spaces are defined as
\subeqn{
  \label{eq:MHD_volume_spaces}
  \algn{
      \GBM &:= \LRc{ \GB \in \LRs{ L^2(\Omega_h) }^{d \times d} \; : \; \GB|_\K \in \LRs{\Poly^\p\LRp{\K}}^{d \times d}, \forall \K \in \Omega_h}, \\
    \VbM &:= \LRc{ \vb \in \LRs{ L^2(\Omega_h) }^d            \; : \; \vb|_\K \in \LRs{\Poly^\p\LRp{\K}}^d, \forall \K \in \Omega_h }, \\
    \WM  &:= \LRc{ w   \in       L^2(\Omega_h)                \; : \;   w|_\K \in \Poly^\p\LRp{\K}, \forall \K \in \Omega_h }, \\
    \HbM &:= \LRc{ \Hb \in \LRs{ L^2(\Omega_h) }^{\tilde d}          \; : \; \Hb|_\K \in  \LRs{\Poly^\p\LRp{\K}}^{\tilde d}, \forall \K \in \Omega_h}, \\
    \CbM &:= \LRc{ \cb \in \LRs{ L^2(\Omega_h) }^d            \; : \; \cb|_\K \in  \LRs{\Poly^\p\LRp{\K}}^d, \forall \K \in \Omega_h}, \\
    \SM  &:= \LRc{ s   \in       L^2(\Omega_h)                \; : \;   s|_\K \in \Poly^\p\LRp{\K}, \forall \K \in \Omega_h},
  }
}
where $\tilde d$ takes the value of one for 2D and three for 3D.
We define the skeletal spaces as follows,
\algn{
    \VbhM &:= \LRc{ \vbh \in \LRs{ L^2(\Gh) }^d \; : \; \vbh|_e \in  \LRs{\Poly^\p\LRp{\e}}^d, \forall \e \in \Gh},\\
  \CbhtM  &:= \LRc{ \cbht \in \LRs{ L^2(\Gh) }^{d-1} \; : \; \cbht|_e \in \CbhtM(e) },\\ 
  \ShM  &:= \LRc{ \sh \in L^2(\Gh) \; : \; \sh|_e \in \Poly^\p\LRp{\e}, \forall \e \in \Gh}. 
}
Here $\CbhtM(e)$ is a vector valued polynomial space with no normal component, defined by
\algn{
  \CbhtM(e) 
  &= \LRc{ \sum\limits_{i=1}^{d-1} \tb^i \ch_{h,i} \; : \; \ch_{h,i} \in \Poly^\p\LRp{\e}, \forall \e \in \Gh },
}
where $\tb^i$ are tangent vectors to $e$. The values 
$\xi=\frac{1}{2}$, $\beta_n=\beta_t=1$ are chosen and the
stabilization $\SB_u$ is taken as
\algn{
  \label{eq:stab_tensor_MHD_u}
  \SB_u &:= \taut \TB + \taun \NB,
}
where $\NB:=\nb\otimes\nb$, $\TB:=-\nb\times(\nb\times .)=\IB-\NB$, $\taut = \frac{1}{2}\sqrt{4+(\wb\cdot\nb)^2}$ and $\taun = \frac{1}{2}\sqrt{8+(\wb\cdot\nb)^2}$.

The well-posedness of the local and global solvers of the scheme \eqref{eq:MHD_HDG_2} and also the error analysis 
is shown in \cite{stephen_thesis,lee2019analysis}. For this scheme the volume velocity and magnetic fields 
converge optimally as $\mc{O}(h^{\p+1})$, whereas all the other volume variables converge as $\mc{O}(h^{\p+1/2})$. The 
verification of this scheme for a number of prototypical MHD problems is shown in \cite{stephen_thesis,lee2019analysis}.

After discretizing the time derivative terms in equation \eqref{eq:MHD_HDG_2} by means of some time discretization scheme e.g., backward Euler,
for each nonlinear iteration of Picard we need to solve a linear system. To that extent, we express the volume unknowns 
in terms of the skeletal unknowns through the local solvers, then we use the conservation
conditions to generate the global linear system which has the following block form \cite{stephen_thesis}
\algn{
  \label{eq:MHD_HDG_2_matrix}
  \LRs{
    \begin{array}{cccc}
      A & -D^\top & E & G \\ 
      D & 0 & 0 & 0 \\
      F & 0 & C & J \\
      H & 0 & K & L
    \end{array}
  }
  \LRs{
    \begin{array}{c}
      \widehat{U} \\ 
      \rho \\ 
      \widehat{B}^t \\ 
      \widehat{R}
    \end{array}
  }
  =
  \LRs{
    \begin{array}{c}
      F_1 \\ 
      F_2 \\ 
      F_3 \\ 
      F_4
    \end{array}
  }.
}
We will use this block structure and develop a preconditioning strategy as shown in the next section.

\section{A block preconditioner for the linear system}
\seclab{block}
Before moving into the construction of the block preconditioner we first briefly
explain the need for it in this case. If we want to precondition the linear system
\eqref{eq:MHD_HDG_2_matrix} using multigrid or multilevel methods, we cannot apply it directly 
because of the difference in the nature of the trace unknowns. The unknowns $(\widehat{U},\widehat{B}^t,\widehat{R})$ are all nodal skeletal unknowns belonging to $\Poly^\p(\e)$, whereas
the edge average pressure $\rho$ is an element-wise constant and is independent
of the solution order $\p$. Thus with standard multigrid or multilevel methods, 
coarsening becomes an issue unless different strategies are employed for the different types of trace unknowns. This problem is also encountered in the linear systems arising from
mixed finite element methods and a strategy to tackle this issue is block
preconditioning. The idea is to identify and group blocks corresponding to different 
unknowns and use approximate block inverses for preconditioning \cite{elman2006block,elman1999preconditioning,elman2014finite,rudi2019global,cyr2013new,
phillips2016block,wathen2018preconditioners}. Here we use similar techniques to develop a preconditioner for the linear system \eqref{eq:MHD_HDG_2_matrix}.

Towards this goal we first rewrite equation \eqref{eq:MHD_HDG_2_matrix} into the saddle point form as follows
\algn{
  \label{eq:MHD_HDG_saddle}
  \LRs{
    \begin{array}{cccc}
      A & E & G & -D^\top \\ 
      F & C & J & 0 \\
      H & K & L & 0 \\
      D & 0 & 0 & 0
    \end{array}
  }
  \LRs{
    \begin{array}{c}
      \widehat{U} \\
      \widehat{B}^t \\ 
      \widehat{R} \\
        \rho
    \end{array}
  }
  =
  \LRs{
    \begin{array}{c}
      F_1 \\ 
      F_3 \\ 
      F_4 \\ 
      F_2
    \end{array}
  }.
}
Denoting the $3\times3$ block corresponding to the unknowns $(\widehat{U},\widehat{B}^t,\widehat{R})$ as $\mc{F}$ and $\LRs{D \quad 0 \quad 0}$ as $\mc{B}$ we can write 
the above matrix as
\algn{
  \label{eq:2x2_form}
  \LRs{
    \begin{array}{cc}
        \mc{F} & -\mc{B}^\top \\ 
           \mc{B} & 0 
    \end{array}
  }.
  }
For a $2\times2$ block matrix such as \eqref{eq:2x2_form} its block inverse (assuming $\mc{F}^{-1}$ and $\mc{S}^{-1}$ exists) can be written as \cite{vassilevski2008multilevel}
\algn{
  \label{eq:2x2_inverse}
  \LRs{
    \begin{array}{cc}
        \mc{F} & -\mc{B}^\top \\ 
           \mc{B} & 0 
    \end{array}
}^{-1}
=
\LRs{
    \begin{array}{cc}
        \mc{F}^{-1} & \mc{F}^{-1}\mc{B}^\top\mc{S}^{-1} \\ 
           0 &  \mc{S}^{-1}
    \end{array}
}
\LRs{
    \begin{array}{cc}
        I & 0 \\ 
           -\mc{B}\mc{F}^{-1} &  I
    \end{array}
},
}
where $\mc{S}:=\mc{B}\mc{F}^{-1}\mc{B}^\top$ is the Schur complement.  
Now, when we use the block upper triangular matrix of the inverse \eqref{eq:2x2_inverse}
as a right preconditioner for the saddle point matrix \eqref{eq:2x2_form} we get
\algn{
  \LRs{
    \begin{array}{cc}
        \mc{F} & -\mc{B}^\top \\ 
           \mc{B} & 0 
    \end{array}
  }
\LRs{
    \begin{array}{cc}
        \mc{F}^{-1} & \mc{F}^{-1}\mc{B}^\top\mc{S}^{-1} \\ 
           0 &  \mc{S}^{-1}
    \end{array}
}
=
\LRs{
    \begin{array}{cc}
        I & 0 \\ 
           \mc{B}\mc{F}^{-1} &  I
    \end{array}
}.
}
 All the eigenvalues of the preconditioned matrix have the value 1, and hence with a Krylov subspace method such as GMRES at most two iterations are needed to solve the system \cite{murphy2000note,elman2006block,elman1999preconditioning,elman2014finite}.

However, the problem with this ideal preconditioner is that we need inverses of 
$\mc{F}$ and $\mc{S}$ which are expensive to compute. Hence a natural idea is
to use approximations of these inverses in the construction of the preconditioner.
First, let us consider the approximation for the inverse of the Schur complement
matrix $\mc{S}$. We will follow the approach shown in \cite{elman2006block,elman1999preconditioning,elman2014finite} for incompressible Navier--Stokes equations to derive
an approximation for the inverse of the Schur complement.

To that extent, if we can find an approximate commutator $\tilde{\mc{F}}$ such that 
\algn{
    \label{eq:commutator}
    \mc{B}^\top \tilde{\mc{F}} \approx \mc{F} \mc{B}^\top,
}
then pre-multiplying by $\mc{F}^{-1}$ and post-multiplying by $\tilde{\mc{F}}^{-1}$ (assuming both inverses exists) on both sides of the above equation we get
\algn{
    \label{eq:commutator_mod}
    \mc{F}^{-1} \mc{B}^\top \approx \mc{B}^\top \tilde{\mc{F}}^{-1}.
}
Using \eqref{eq:commutator_mod} we can write an approximate inverse for the Schur complement as
\algn{
    \mc{S}^{-1} &= \LRp{\mc{B}\mc{F}^{-1} \mc{B}^\top}^{-1} \\
            &\approx \LRp{\mc{B}\mc{B}^\top \tilde{\mc{F}}^{-1}}^{-1} \\
    \label{eq:Schur_Approx}
 \tilde{\mc{S}}^{-1}           &= \tilde{\mc{F}}\LRp{\mc{B}\mc{B}^\top}^{-1}.
}

Now the only task remaining is to find the approximate commutator of $\mc{F}$ and we can use least squares minimization for that. 
Solving for the normal equations corresponding to \eqref{eq:commutator} gives \cite{elman2006block}
\algn{
    \tilde{\mc{F}} = \LRp{\mc{B}\mc{B}^\top}^{-1}\mc{B}\mc{F}\mc{B}^\top.
}
Substituting the above equation in the approximate inverse for the 
Schur complement \eqref{eq:Schur_Approx} we get
\algn{
    \label{eq:BFBT}
  \tilde{\mc{S}}^{-1} = \LRp{\mc{B}\mc{B}^\top}^{-1}\mc{B}\mc{F}\mc{B}^\top\LRp{\mc{B}\mc{B}^\top}^{-1}.
}
This is called BFBT approximation (because of the middle term sandwiched between two inverses) of the inverse Schur complement in the 
literature \cite{elman1999preconditioning,elman2006block,elman2014finite,rudi2019global} and
has been successfully used in the context of Stokes and incompressible Navier--Stokes
equations discretized with finite volume, finite difference and mixed FEM methods.
In those cases, $\mc{B}\mc{B}^\top$ corresponds to the Poisson operator and instead of 
taking inverses one typically use a geometric multigrid or AMG cycle \cite{elman1999preconditioning,elman2006block,elman2014finite}.\ Also, for
approximating $\mc{F}^{-1}$, AMG or geometric multigrid cycles are used \cite{elman1999preconditioning,elman2006block,elman2014finite}.\ Thus these choices provide a cheap approximation to the ideal preconditioner in equation \eqref{eq:2x2_inverse} for a certain class
of problems and discretizations. Scaled versions of the BFBT
approximation are proposed in \cite{elman2006block,may2008preconditioned,rudi2019global,elman2014finite} in order to
improve the robustness in certain cases.

The HDG linear system for skeletal unknowns in equation \eqref{eq:MHD_HDG_2_matrix} which comes after static condensation
of volume unknowns is mostly algebraic in nature. Thus it is more difficult to analyze than stabilized or mixed FEM linear systems
especially for complex systems like MHD. Moreover, our grouping of unknowns in equation \eqref{eq:MHD_HDG_saddle} to arrive
at the saddle point form is also algebraic and we can see that already in the $\mc{F}$ matrix which contains part of the fluid block (without pressure) and
the magnetic block. We base our choices for $\LRp{\mc{B}\mc{B}^\top}^{-1}$ and $\mc{F}^{-1}$ considering this aspect.

In the context of HDG discretization of the incompressible MHD system, $\mc{B}\mc{B}^\top$
corresponds to a matrix of size $\Nel\times\Nel$ which is much smaller in 
size compared to $\mc{F}$. It is independent of the solution order $\p$ and the
bandwidth is also small since it corresponds to the edge-average pressure which is piecewise 
constant per element. Since in high order HDG we try to minimize the number of elements while increasing the
order to increase the computation to communication ratio the $\mc{B}\mc{B}^\top$ matrix is small in this scenario. 
Thus as a first step we will use a parallel sparse 
direct solver \verb|SuperLU_DIST| \cite{li2003superlu_dist} for taking inverse of 
this matrix in the BFBT approximation. In our future work we will replace 
this part with preconditioned conjugate gradient solver of lenient tolerance 
to improve scalability. Some of the initial studies conducted in this direction shows promise. 

For approximating $\mc{F}^{-1}$, we use one v-cycle of a system projection AMG solver from the ML library \cite{gee2006ml} of Trilinos project \cite{heroux2005overview}. Similar to 
\cite{shadid2016scalable}, we order the unknowns in the nodal block such that the degrees-of-freedom within each node appear consecutively. This helps to preserve the coupling between different variables during coarsening. The 
ordering of unknowns within each node is $(\widehat{U},\widehat{B}^t,\widehat{R})$ i.e., 
velocity trace unknowns are ordered first, followed by the tangent magnetic 
field and then the Lagrange multiplier. Other orderings are also possible and it can affect the performance of the AMG cycle. Since it is beyond the current scope of this article we intend to compare them in our future work. We use non-smoothed and uncoupled aggregation with a sparse direct solver on the coarsest level as in \cite{shadid2016scalable}. More details about the aggregation strategy employed in AMG can be found in 
\cite{shadid2016scalable}. The choice of smoothers inside the AMG cycle is critical and the mixed parabolic-hyperbolic nature of MHD together with the algebraic nature of the HDG trace system requires us to use strong smoothers such as ILU(0) and GMRES instead of standard smoothers like Jacobi and Gauss-Seidel.

Another choice for $\mc{F}^{-1}$ is the multilevel preconditioner introduced in \cite{muralikrishnan2020multilevel}. Since the coarse
solver in the multilevel preconditioner is obtained by an approximation of the nested dissection direct solver it can also provide an effective approximation to $\mc{F}^{-1}$. We compare the results of the multilevel preconditioner to that of AMG for the approximation to $\mc{F}^{-1}$ in section \secref{amg_vs_ml}. In our numerical studies we observed that the performance of the scaled BFBT approximation is very similar to the non-scaled one \eqref{eq:BFBT} and hence we use that in all our
experiments. In summary we use a right preconditioner 
\algn{
    \label{eq:BFBT_AMG}
\LRs{
    \begin{array}{cc}
        \tilde{\mc{F}}^{-1} & \tilde{\mc{F}}^{-1}\mc{B}^\top\tilde{\mc{S}}^{-1} \\ 
           0 &  \tilde{\mc{S}}^{-1}
    \end{array}
},
}
where $\tilde{\mc{F}}^{-1}$ is one v-cycle of AMG or the multilevel algorithm (Algorithm 1) in \cite{muralikrishnan2020multilevel} and
the approximate Schur complement inverse $\tilde{\mc{S}}^{-1}$ is given by the BFBT approximation \eqref{eq:BFBT}. The $\LRp{\mc{B}\mc{B}^\top}^{-1}$
in the BFBT approximation is obtained by a linear system solve with parallel direct solver \verb|SuperLU_DIST|.

\section{Numerical results}
\seclab{numerical_mhd}
In this section we test the performance of the block preconditioner on some of the transient test cases in incompressible resistive MHD considered before \cite{shadid2016scalable,cyr2013new,phillips2016block}. In particular, we consider 2D and 3D versions of the island coalescence problem, hydromagnetic Kelvin-Helmholtz (HMKH) instability and hydromagnetic lid-driven cavity problems. We use quadrilateral elements in 2D and hexahedral elements in 3D. For time integration we use the backward Euler time stepping for all the 3D test cases and five stage fourth order diagonally implicit Runge-Kutta (DIRK) method of \cite{kennedy2016diagonally} for all the
2D test cases. For the nonlinear solver we employ the Picard iteration scheme with a stopping criterion based
on the discrete norm of the solution update vector given in \cite{shadid2016scalable} 
\[
    \sqrt{\frac{1}{N_u}\sum_{i=1}^{N_u}\LRs{\frac{|\Delta \chi_i|}{\varepsilon_r |\chi_i| + \varepsilon_a}}^2} < 1.
\]
Here, $\chi_i$ is any variable, $\Delta \chi_i$ is its corresponding update in the Picard iteration, $N_u$ is the total
number of unknowns and $\varepsilon_a$ and $\varepsilon_r$ are the absolute and relative 
tolerances which are set as $10^{-6}$ and $10^{-4}$ respectively. For the stopping tolerance of the linear solver, apart from
the HMKH test case, we set the value to be $10^{-6}$ multiplied by the norm of the right hand side of the Picard linear system.
For the HMKH test case, it turns out we need a stricter tolerance of $10^{-9}$ (not scaled by the norm of the right hand side) to 
make the Picard iterations converge. 

For all the parallel results, we implemented our algorithms on top of the deal.II FEM library \cite{BangerthHartmannKanschat2007,alzetta2018deal}. The weak scaling studies are conducted in the Knights Landing (KNL) nodes of
the Stampede2 supercomputer at the Texas Advanced Computing Center. Each node of KNL consists of 68 intel Xeon Phi 7250 1.4GHz processors with 4 threads per core. It has 96GB DDR4 RAM along with 16GB high speed MCDRAM which acts as L3 cache. We have used only MPI 
parallelism even though our deal.II code has task based parallelism using thread building blocks (TBB) in addition to MPI. The reasons
for this choice are: (i) to have memory locality and avoid memory contention which may
complicate the weak scaling studies (ii) the main focus of this study, which is the linear solver part, uses ML from Trilinos which does not have the threads support. In our future work we will use the second generation Trilinos library 
MueLu \cite{prokopenko2014muelu} instead of ML as it provides support for threads and GPU.

\subsection{Island coalescence}
\seclab{island}
Magnetic reconnection is a fundamental phenomenon by which a magnetic field changes its structure and is accompanied by conversion of magnetic field energy into plasma energy and transport \cite{shadid2016scalable}. Many physical phenomena which occurs
in outer space such as solar flares, coronal mass ejections involve magnetic reconnection as an important driving mechanism. 
It is also important in a laboratory scenario, especially in fusion experiments 
to understand and control plasma disruptions which
can lead to loss in plasma confinement and also damages to the machine. Since fusion reactors like tokamak are typically designed to handle only certain 
maximum number of these disruptions, understanding and controlling these phenomena is of significant interest to the fusion and in general plasma community. More details about magnetic reconnection can be found in \cite{biskamp1996magnetic,goedbloed2004principles,goedbloed2010advanced}.

While magnetic reconnection is important for its physical significance, it is also characterized by disparate spatial and temporal scales which serves as an ideal test bed for testing the robustness of our preconditioners/solvers.
In this section we consider a specific reconnection problem  which is the island coalescence studied in \cite{shadid2016scalable,chacon2008optimal}. It initially
consists of two islands embedded in a Harris current sheet as shown in Figure \figref{3D_initial_island}. A perturbation to the initial configuration and the combined magnetic field causes the center of the islands (referred as the o-points) to move towards each other and eventually coalesce
to form one island. When the reconnection happens, the islands form a `x' structure in the center of the domain as shown 
in Figures \figref{3D_x1} and \figref{J_p_t_150}. In what follows we will briefly describe the settings of this problem and then
evaluate the performance of our preconditioner in this case.

The domain is $[-1, 1]^d$, where $d$ is the dimension. The boundary conditions are periodic in the x-direction on the left and
right faces and also in the z-direction on the back and front faces. On
the top and bottom faces in the y-direction, for the magnetic part, perfect conducting boundary
conditions described by zero normal magnetic field $\bb \cdot \n = 0$ and zero tangential
electric field $\n \times \Eb = 0$ are applied. For the fluid part on the top and bottom faces, mirror boundary 
conditions described by zero normal velocity $\ub \cdot \n = 0$ and zero shear stress are applied. The Lagrange multiplier 
$r$ is set as zero on all the boundaries. We refer the readers to \cite{stephen_thesis} for an application of these boundary 
conditions in the context of HDG discretizations. 

The initial conditions consists of zero fluid velocity ($\ub^0 = {\bf 0}$), and the magnetic field given by
\[
  \bb^0 = \LRp{ \frac{          \sinh \LRp{ 2 \pi y } }{ \cosh\LRp{ 2 \pi y } + \epsilon \cos\LRp{ 2 \pi x } },
          \frac{ \epsilon \sin  \LRp{ 2 \pi x } }{ \cosh\LRp{ 2 \pi y } + \epsilon \cos\LRp{ 2 \pi x } }, 0 },
      \]
where in 2D the first two components of the field values are used. Here, $\epsilon$ refers to the width of the island and we choose it to be 0.2 as in \cite{shadid2016scalable}. In order for the initial configuration to be in
equilibrium a forcing of $\gb = \nabla \times \Jb^0$ is used where
\[
    \Jb^0 := \Jb( t = 0 ) = \LRp{0, 0, 
  - \frac{ 2 \pi \kappa \LRp{ 1 - \epsilon^2 } }{ \Rm \LRp{ \cosh\LRp{ 2 \pi y } + \epsilon \cos\LRp{ 2 \pi x } }^2 }}.
\]
For the momentum equation zero forcing ($\fb={\bf 0}$) is selected. To set the islands 
into motion in a reproducible manner rather than relying on the accumulation of 
round-off error an initial perturbation of
\[
  \delta \bb^0 = \sigma \LRp{ \frac{ \pi }{ 2 } \cos\LRp{ \pi x } \sin\LRp{ \frac{ \pi y }{ 2 } } \cos\LRp{ \pi z },
  - \pi \sin\LRp{ \pi x } \cos\LRp{ \frac{ \pi y }{2} } \cos\LRp{ \pi z }, 0},
\]
is used with the value of $\sigma=10^{-3}$ which sets the magnitude of perturbation \cite{stephen_thesis}. For 2D, the first two components are used without the $z-$term. 
As described in \cite{stephen_thesis} we choose the characteristic velocity as Alfv\'{e}n 
speed which gives $\kappa=1$. Also in all our numerical experiments we set the fluid Reynolds number and magnetic Reynolds 
number (which is Lundquist number in this case) equal to each other.


\begin{figure}[h!b!t!]
  \subfigure[$t = 0.1$]{
      \figlab{3D_initial_island}
    \includegraphics[trim=7.5cm 0cm 6.5cm 1.5cm,clip=true,width=0.48\columnwidth]{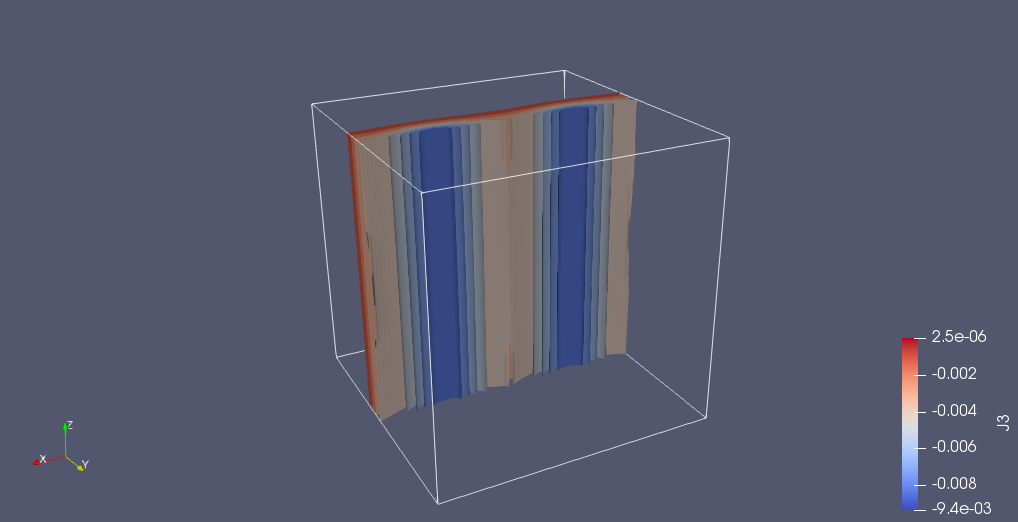}
  }
  \subfigure[$t = 1.1$]{
    \includegraphics[trim=7.5cm 0cm 6.5cm 1.5cm,clip=true,width=0.48\columnwidth]{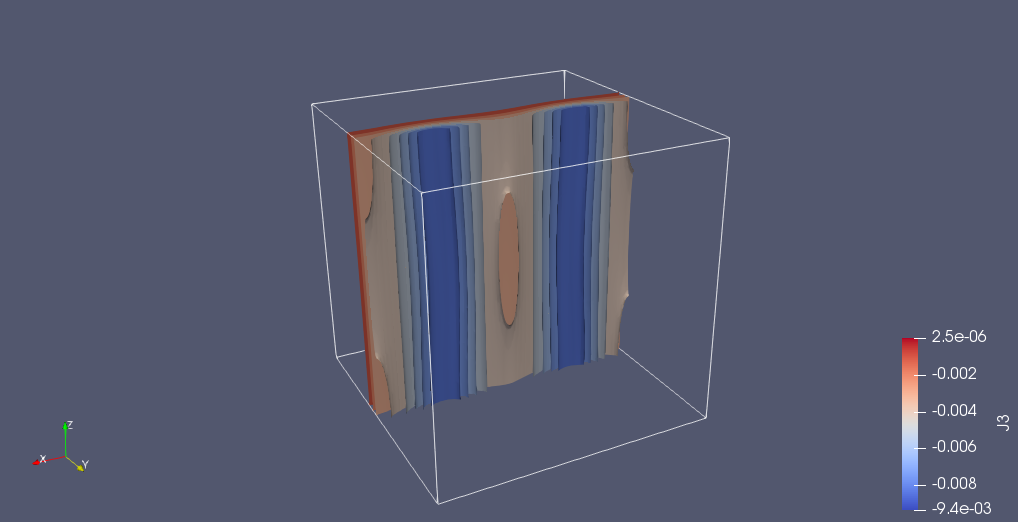}
  }
    \subfigure[$t = 2.1$]{
    \includegraphics[trim=7.5cm 0cm 6.5cm 1.5cm,clip=true,width=0.48\columnwidth]{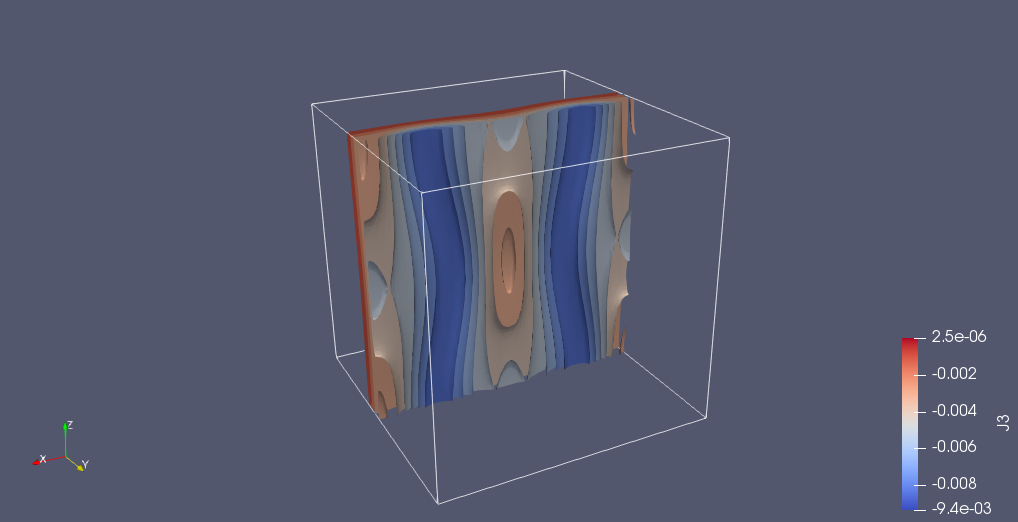}
  }
  \subfigure[$t = 3.1$]{
      \figlab{3D_x1}
    \includegraphics[trim=7.5cm 0cm 6.5cm 1.5cm,clip=true,width=0.48\columnwidth]{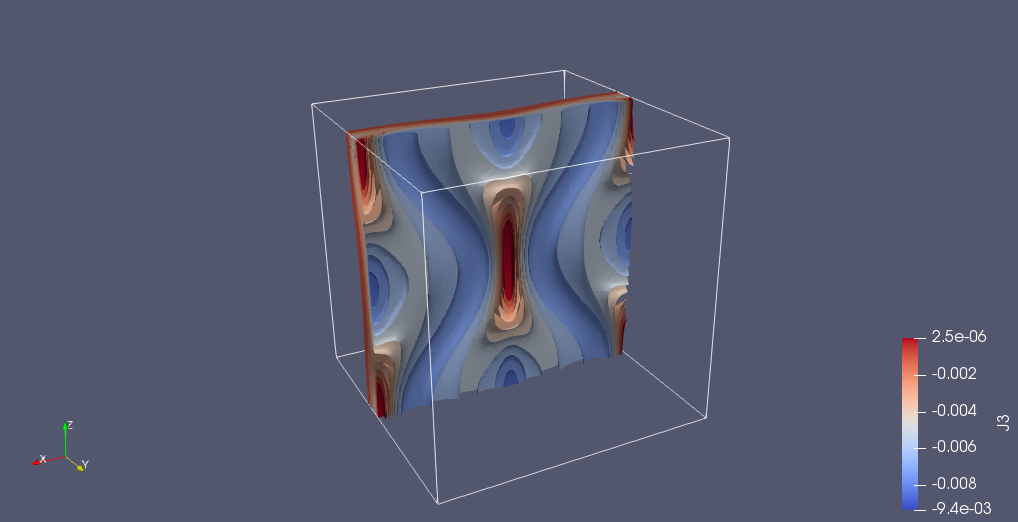}
  }

    \caption{Evolution of the z-component of the current ($J_z$) with time. The contours of $J_z$ show a highly kinked state at $t=3.1$ due to the perturbation in the $z-$direction.}
    \figlab{3D_island}
\end{figure}

Let us choose a Lundquist number $S=10^3$, and a $16\times12\times14$ clustered mesh of order $\p=5$. Figure \figref{3D_island}, shows the evolution of the $z-$component of current by taking a slice of the contour plot at mid-plane $z=0$.
As can be seen the perturbation in the $z-$direction causes the current contours to bend and results in a highly
kinked state as evidenced in Figure \figref{3D_x1}. Our results show good agreement (visually) with the results in \cite{chacon2008optimal,shadid2016scalable}.

Now we compare the performance of three preconditioners for this problem in a weak scaling sense i.e., we increase the problem size proportionally with the increase in number of processors so that the number of elements per processor remains constant. We take $4^3$, $8^3$, $16^3$ and $32^3$ elements and 2, 16, 128 and 1024 processors respectively so that we have 32 elements per core. A time stepsize of $0.1$ is selected and all the results are averaged over six time steps. The nonlinear Picard solver takes on average $3.2-4$ iterations in all the cases. 

First of the preconditioners is a one level domain decomposition method with an incomplete factorization 
sub-domain solver with zero fill-in (ILU(0)) and overlap of one\footnote{Here, overlap one means each processor will include its own set of rows and in addition it also includes the rows corresponding to its non-zero columns.} (denoted as DD, ILU(0) in Figure \figref{Island_weak_p4} and let us refer to this as DD with ILU(0)). The other two preconditioners are the block preconditioners given in equation \eqref{eq:BFBT_AMG} with BFBT 
approximation for $\tilde{\mc{S}}^{-1}$ and one AMG v-cycle for $\tilde{\mc{F}}^{-1}$. The difference between them is the smoothing
inside AMG, in one of 
them we use the ILU(0) smoother of overlap one (denoted as BFBT+AMG, ILU(0) in Figure \figref{Island_weak_p4} and let us refer to this as BFBT+AMG with ILU(0)), whereas in the other one we use the GMRES preconditioned by ILU(0) of overlap zero\footnote{Here, overlap zero means each processor will include only its own set of rows and thus no communication is needed.} as smoother (denoted as BFBT+AMG, GMRES in Figure \figref{Island_weak_p4} and let us refer to this as BFBT+AMG with GMRES). The reason for these non-standard choice of smoothers is, classical smoothers
like Jacobi, Chebyshev and Gauss--Seidel did not result in converging iterations. This is also observed 
in \cite{lin2017performance} for linear systems arising from stabilized FEM discretization of MHD. Thus it indicates that strong smoothers are needed for AMG cycles applied to linear systems coming from incompressible resistive MHD. 

We perform three pre- and three post-smoothing steps (in the case of GMRES smoother, these denote the number of inner iterations), whereas in the DD with ILU(0) preconditioner we perform three smoothing steps. The outer iterations are performed with GMRES
for the DD with ILU(0) and the BFBT+AMG with ILU(0) preconditioners. For the BFBT+AMG with GMRES, we use flexible GMRES (FGMRES) \cite{saad1993flexible} for 
the outer iterations due to the nonlinear nature of the inner iterations. In both cases similar to \cite{shadid2016scalable} we use non-restarted versions as it may result in degradation of convergence.

\begin{figure}[h!b!t!]
  \subfigure[]{
      \figlab{iter_island_p4}
    \includegraphics[trim=2cm 6cm 2cm 7cm,clip=true,width=0.48\columnwidth]{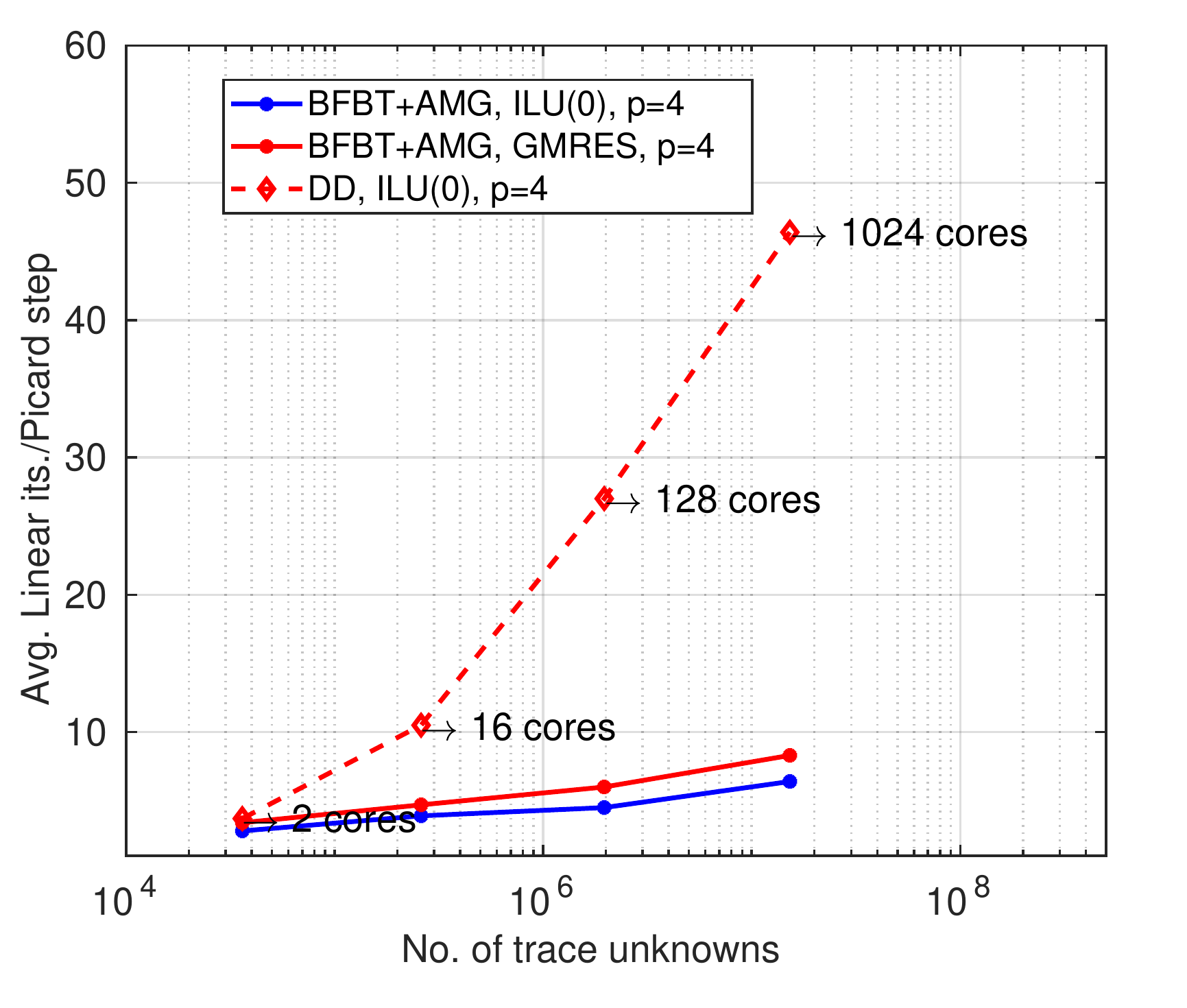}
  }
  \subfigure[]{
      \figlab{time_island_p4}
    \includegraphics[trim=1.5cm 6cm 2cm 7.0cm,clip=true,width=0.48\columnwidth]{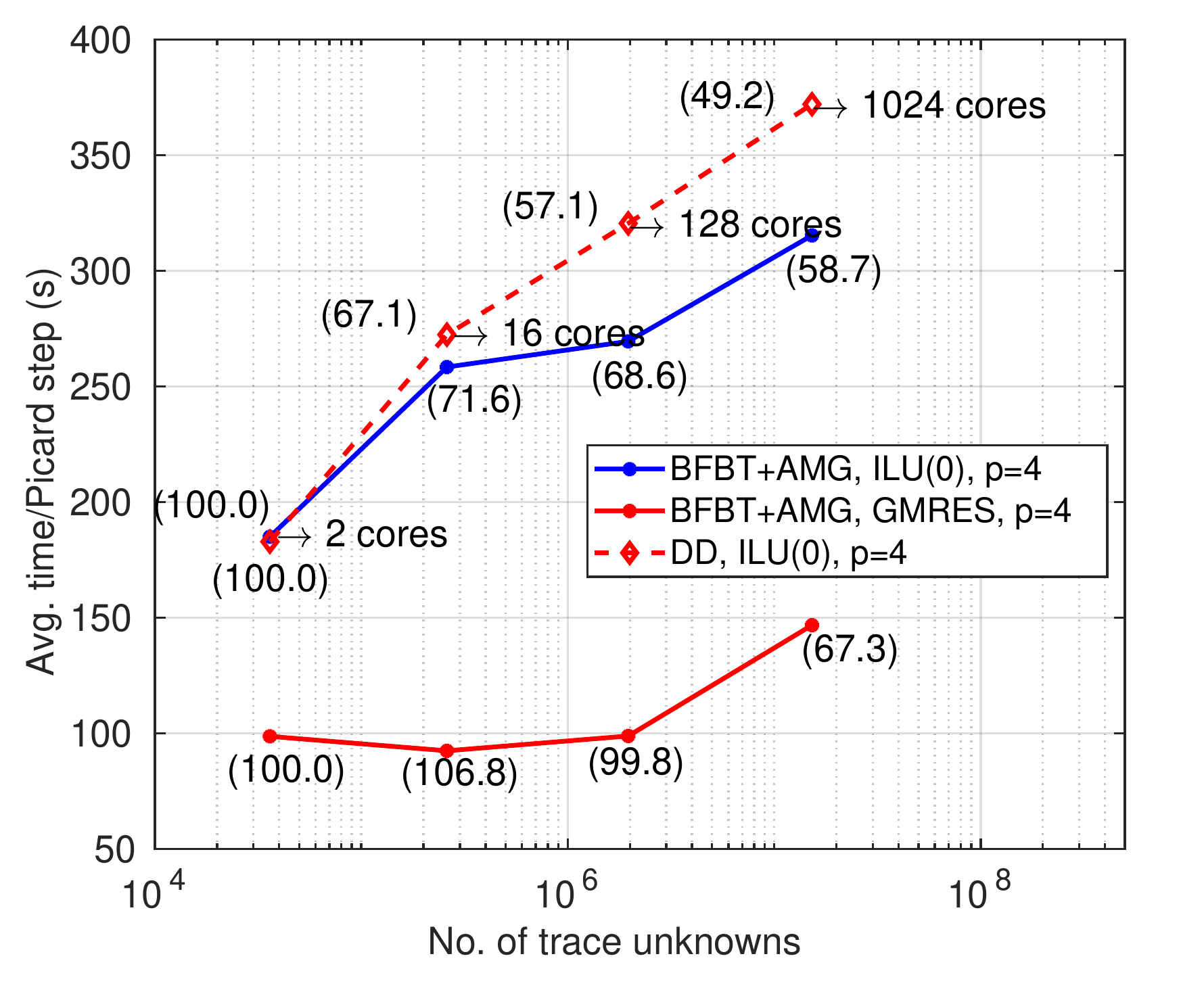}
  }

    \caption{3D island coalescence problem: weak scaling study of average iterations per Picard step (left) and average time per Picard step (right) for three preconditioners and solution order $\p=4$. The markers in BFBT+AMG with ILU(0) and 
    GMRES also represent the same number of processors as in DD with ILU(0). The values within parentheses in the right figure represent weak scaling parallel efficiencies.}
    \figlab{Island_weak_p4}
\end{figure}

We compare the performance in terms of average iterations per Picard step and average time per Picard step in Figure \figref{Island_weak_p4}. 
From Figure \figref{iter_island_p4} for solution order $\p=4$, we can see that the iterations of DD with ILU(0) preconditioner increases much faster than 
the other two block preconditioners and this is due to the lack of coarse solvers. Similar results are 
also observed in \cite{shadid2016scalable} for the stabilized FEM discretization. The average linear iterations 
per Picard step for the two block preconditioners are less than the one level domain decomposition preconditioner
and the growth with number of unknowns/processors (here the number of unknowns refers to the number of trace unknowns only) is also very mild. Between ILU(0) and GMRES smoothers we can see that the GMRES smoother
takes slightly more iterations than the ILU(0). 

In Figure \figref{time_island_p4}, we compare the average time per
Picard step for the three preconditioners and again DD with ILU(0) preconditioner takes more time compared to the
block preconditioners. However, a fact to notice is that the large difference in iteration counts evidenced in Figure 
\figref{iter_island_p4} between DD with ILU(0) and BFBT+AMG with ILU(0) is not much reflected in the timing figure in 
\figref{time_island_p4}. This is because ILU(0) with overlap one spends most of the time in the
setup cost and very little time in the application of the preconditioners. Since the setup cost in both DD with ILU(0) 
and BFBT+AMG with ILU(0) are very similar the difference comes only from the application of the preconditioners
which is high for DD with ILU(0) due to larger number of iterations. The scenario for the GMRES smoother is
exactly opposite with very little time spent for the setup and most of the cost coming from the application or solving
part of the preconditioner. This reflects in the timing results in Figure \figref{time_island_p4} with BFBT+AMG with GMRES
smoother taking the least time (in spite of its iteration count bit higher than ILU(0) smoother) and more than 
two times faster than the other two preconditioners. We want to make a remark here that even though the number of volume
unknowns per core is same in this weak scaling study the number of trace unknowns per core is not exactly the same and the 
cases of 16 cores, 128 cores and 1024 cores has slightly less number of trace unknowns per core compared to the 2 core case.
This is the reason in Figure \figref{time_island_p4} the time per Picard step for GMRES smoother decreases slightly for 16 cores
case even though the iterations increase a bit compared to that for 2 cores. However, asymptotically this 
difference gets negligible and hence does not affect much our inferences.

\begin{table}[h!b!t!]
\centering
\begin{tabular}{|r|c|c||c|c||c|c||}
\hline
    & \multicolumn{2}{c||}{DD, ILU(0)} & \multicolumn{2}{c||}{BFBT+AMG, ILU(0)} & \multicolumn{2}{c||}{BFBT+AMG, GMRES}\\
\cline{2-7}
    \!\!\! $S$ \!\!\!\! &  {\!\!\scriptsize iterations \!\!} &  {\!\!\scriptsize time (s)\!\!} &  {\!\!\scriptsize iterations \!\!} &  {\!\!\scriptsize time (s)\!\!} &  {\!\!\scriptsize iterations \!\!} &  {\!\!\scriptsize time (s)\!\!}\\
\cline{2-7}
\hline
$10^3$ & 94 & 489.4 & 3.5 & 287.8 & 5.3 & 124.4 \\
$10^4$ & * & * & * & * & 6 & 127.6 \\
$10^5$ & * & * & * & * & 7.5 & 128.2 \\
$10^6$ & * & * & * & * & 7.6 & 131.2 \\
\hline
\end{tabular}
    \caption{\label{tab:robustness_block}3D island coalescence problem.\ Robustness study with respect to Lundquist number $S$: iterations and time represent the average number of iterations taken per Picard step and the average time per Picard step. ``*" indicates the preconditioners did not reach the specified tolerance for the maximum number of iterations which is taken as $1000$ for DD with ILU(0) and $200$ for BFBT+AMG with ILU(0). The time step in this study
    is $\Delta t=0.05$ and the mesh, solution order and number of cores are $32^3$, $\p=4$ and $1024$ respectively. The average number of Picard iterations taken is approximately 3.}
\end{table}

In Table \ref{tab:robustness_block} we compare the robustness of the preconditioners with respect to Lundquist number $S$. Clearly,
BFBT+AMG with GMRES is more robust than the other two preconditioners
and takes nearly constant number of iterations and time per Picard 
step for the parameters mentioned in Table \ref{tab:robustness_block}.

Now, in addition to computational time and robustness another important aspect to take into account is, ILU(0) smoothers with overlap one require lot of memory and this is much pronounced at high solution orders.
In our numerical experiments we found that for solution orders greater than four, in the KNL nodes of Stampede2, we
are able to use only 8 cores per node or less even though the system has 68 cores per node due to memory limitations.
The GMRES smoother (GMRES preconditioned by ILU(0) with overlap zero) on the other hand has less memory requirements than the ILU(0) smoother with overlap one and this is the 
reason we test only BFBT+AMG with GMRES smoother for solution orders $\p=5,6$ in Figure \figref{Island_weak_p56}.

\begin{figure}[h!b!t!]
  \subfigure[]{
    \includegraphics[trim=2cm 6cm 2cm 7cm,clip=true,width=0.48\columnwidth]{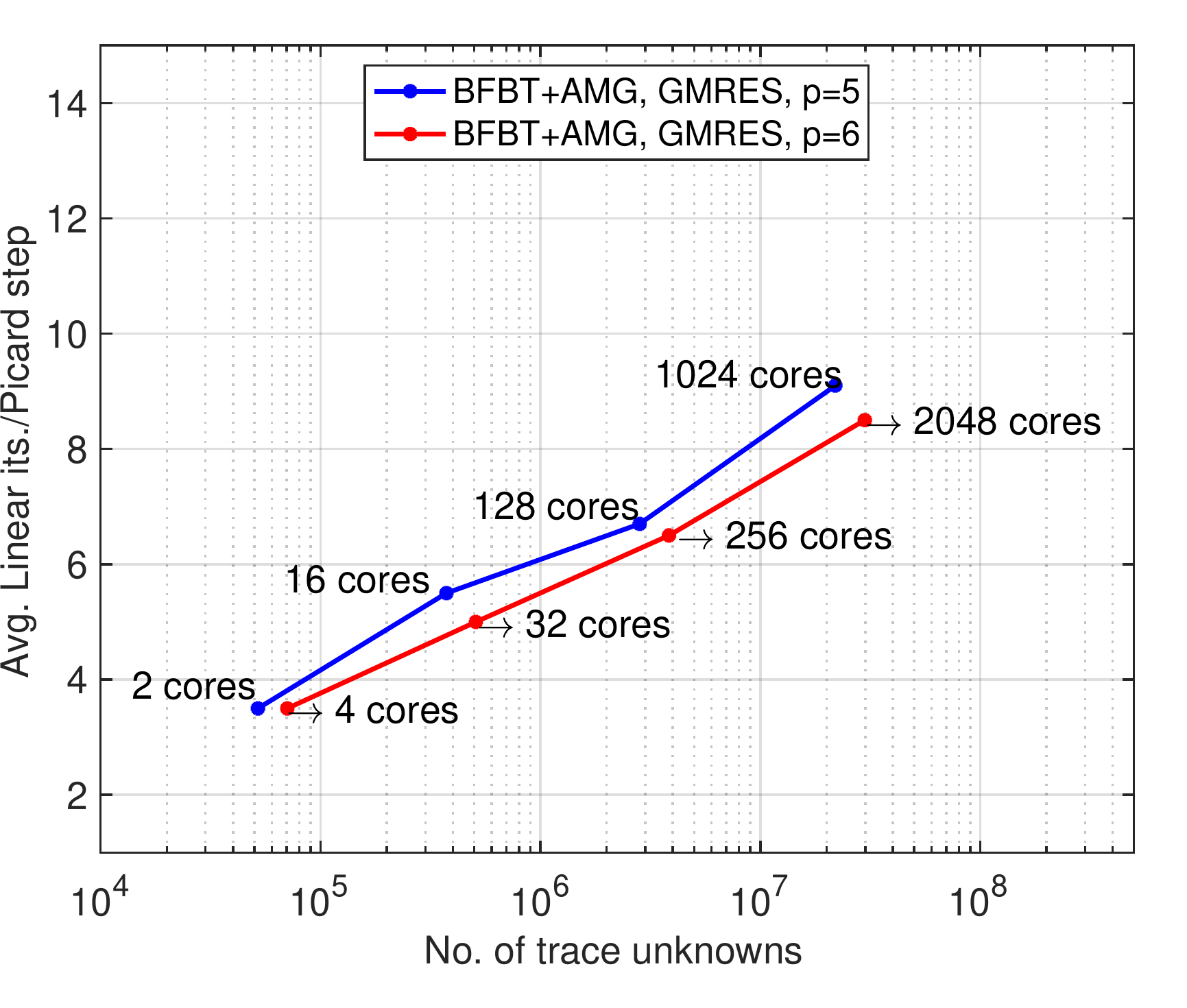}
  }
  \subfigure[]{
    \includegraphics[trim=1.5cm 6cm 2cm 7.5cm,clip=true,width=0.48\columnwidth]{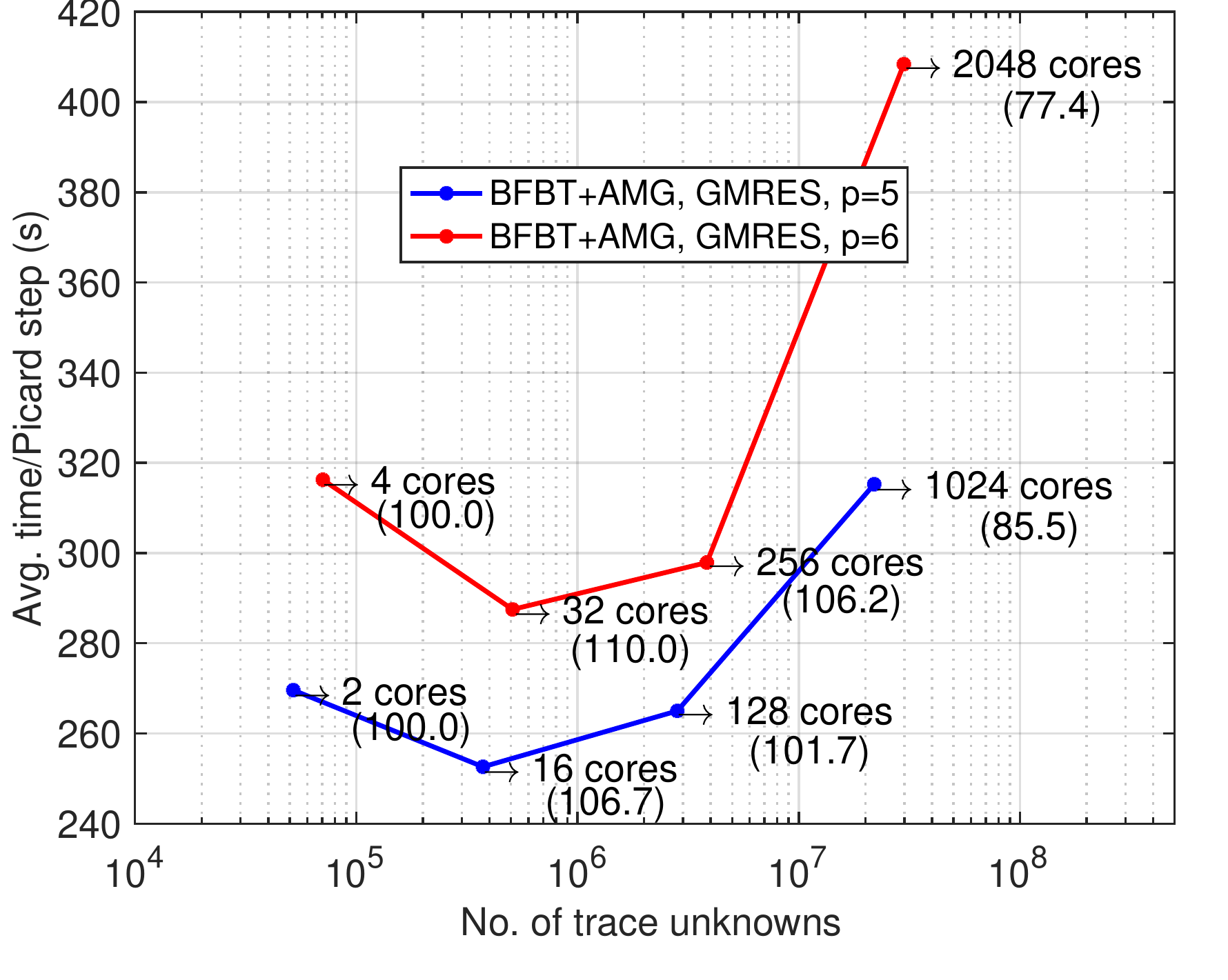}
  }

  \caption{3D island coalescence problem. BFBT+AMG with GMRES smoother: weak scaling study of average iterations per Picard step (left) and average time per Picard step (right) for solution orders $\p=5,6$. The values within parentheses in the right figure represent weak scaling parallel efficiencies.}
    \figlab{Island_weak_p56}
\end{figure}

In Figure \figref{Island_weak_p56}, the average number of iterations and time per Picard step are shown for the BFBT+AMG 
with GMRES smoother for solution orders $\p=5,6$. The iteration counts are very similar to those observed for $\p=4$
in Figure \figref{iter_island_p4} with a mild increase with respect to mesh refinements. The time per Picard step however, for $\p=4,5$ with 1024 cores
and $\p=6$ with 2048 cores show a significant increase compared to 128 and 256 cores respectively. This is a result of two things, (i) increase in iteration count; (ii) coarsening in AMG. We have not done repartitioning with the uncoupled 
aggregation performed in AMG and at high processor counts this resulted in fewer levels (3 or 4) in the AMG hierarchy with larger problem sizes on the coarsest level. Since, we use a serial sparse direct solver on the coarsest level
the timing increased. This trend is observed in the other two numerical experiments also and especially in Figure \figref{lid_weak_p456}
for the lid driven cavity problem. In that case we have almost constant iteration counts which is reflected in the timing till 256 
processors and after that both in 1024 and 2048 processors we see some increase in timing. In our ongoing work we are experimenting with different coarsening strategies in 
AMG and the initial studies have shown promise. We will report these findings as well as optimization of the
other components used in the block preconditioning strategy in the future.

\begin{table}[h!b!t!]
\begin{center}
\begin{tabular}{ | r | c | c | c | c | }
\hline
    \multicolumn{5}{|c|}{$\Nel=32,768$, $N_{trace} = 15.2M$, p=4, $\Delta t =0.1$, averaged over 6 time steps} \\
\hline
    {\#cores} & {time/Picard [s]} & {$\Nel$/core} & Trace DOF/core & {efficiency [\%]} \\
\hline
    128 & 211.6 & 256 & 118.4K & - \\
\hline    
    256 & 126.8 & 128 & 59.2K & 83.4 \\
\hline
    512 & 64.2 & 64 & 29.6K & 82.4 \\
\hline
    1024 & 36.3 & 32 & 14.8K & 73 \\
\hline
    2048 & 19.4 & 16 & 7.4K &  68 \\
\hline
    4096 & 16.6 & 8 & 3.7K &  40 \\
\hline
    6144 & 14.9 & 5.3 & 2.5K & 29.6  \\
\hline
\end{tabular}
\end{center}
\caption{3D island coalescence problem. Strong scaling study for BFBT+AMG with GMRES smoother for solution order $\p=4$. The simulation is carried out in the Skylake nodes of Stampede2 supercomputer.}
    \label{tab:strong_p4}
\end{table}
\begin{table}[h!b!t!]
\begin{center}
\begin{tabular}{ | r | c | c | c | c | }
\hline
    \multicolumn{5}{|c|}{$\Nel=32,768$, $N_{trace} = 29.8M$, p=6, $\Delta t =0.1$, averaged over 6 time steps} \\
\hline
    {\#cores} & {time/Picard [s]} & {$\Nel$/core} & Trace DOF/core & {efficiency [\%]} \\
\hline
    512 & 290 & 64 & 58.4K & - \\
\hline
    1024 & 129.5 & 32 & 29.2K & 112 \\
\hline
    2048 & 82.6 & 16 & 14.6K &  87.7 \\
\hline
    4096 & 46.6 & 8 & 7.3K & 77.7 \\
\hline
    8192 & 39.4 & 4 & 3.7K & 46 \\
\hline
\end{tabular}
\caption{3D island coalescence problem. Strong scaling study for BFBT+AMG with GMRES smoother for solution order $\p=6$. The simulation is carried out in the Skylake nodes of Stampede2 supercomputer.}
    \label{tab:strong_p6}
\end{center}
\end{table}

In Tables \ref{tab:strong_p4} and \ref{tab:strong_p6} we study the strong scalability of BFBT+AMG with GMRES smoother 
for solution orders $\p=4,6$. We see a good strong scaling performance and the efficiency improves with increase in order for similar number
of trace unknowns/core due to increased computation to communication ratio. For the strong scaling study we used the Skylake nodes of Stampede2 supercomputer and since it has a higher clock speed (2.1GHz) compared to the KNL nodes (1.4GHz) the average time 
per Picard step is $4-5$ times lesser than that for KNL nodes. Of the three preconditioners since the BFBT+AMG with GMRES smoother is faster, more robust and requires less memory we will consider only that in the subsequent numerical studies. Also, since the strong scaling performance is not very much problem dependent 
we only study the performance in terms of weak scaling.

We now use the BFBT+AMG with GMRES preconditioner to simulate a challenging problem which is 2D island coalescence at high Lundquist numbers. Here as an
example we consider a Lundquist number of $S=10^7$. The significance of 
this problem is that at high Lundquist numbers the thin current sheet which forms at the center of the domain
during reconnection breaks and gives rise to small structures called plasmoids \cite{bhattacharjee2009fast,huang2010scaling}. The dynamics of this problem is very different from the low Lundquist number cases
where the islands monotonically approach each other and coalesce to a single island. 

The settings of this experiment are as follows: While for the scaling studies we used a fixed time stepsize of $\Delta t=0.1$, 
in this case we use an adaptive time step where we initially start with a stepsize of 0.05 and if the Picard iteration fails to converge in 20 iterations
we cut the stepsize by half and try again. This is essential for capturing the highly nonlinear evolution of this problem and we observed a stepsize
of $0.0015625$ during the plasmoid formation and evolution stages. Once the time stepsize is cut by half we do not increase it later in the simulation, while increasing the stepsize may help in reducing the overall time of the simulation our focus here is to simulate the correct physics and also test the robustness of our preconditioner. The domain
is discretized by $128\times128$ uniform elements and solution order $\p=9$. For time stepping we use the five stage fourth order DIRK scheme. 
The problem is run on 512 cores in the Skylake nodes of Stampede2 supercomputer.


\begin{figure}[h!b!t!]
\centering
\includegraphics[trim=41cm 35cm 41cm 35cm,clip=true,width=\textwidth]{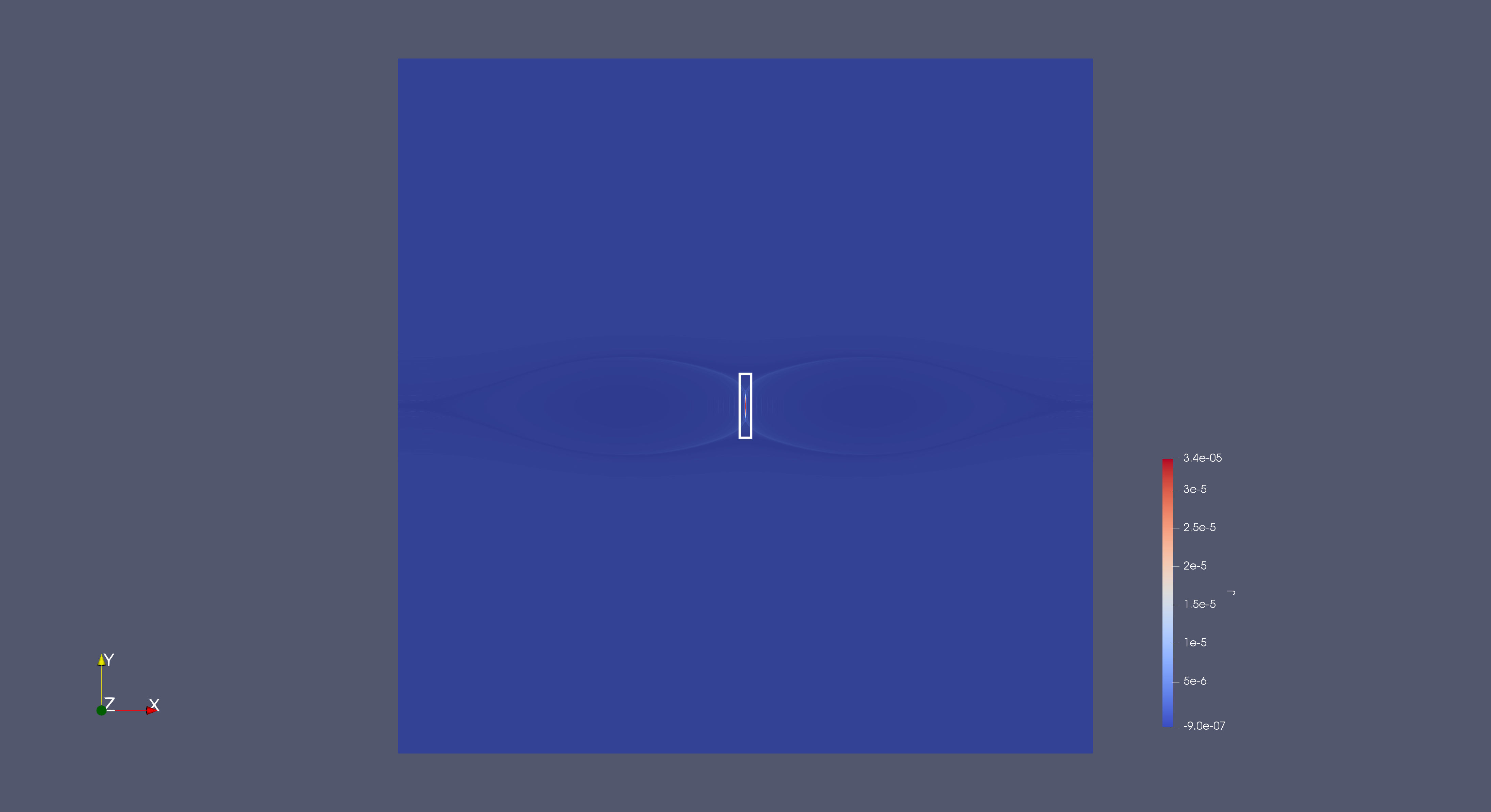}
\includegraphics[trim=41cm 35cm 41cm 35cm,clip=true,width=\textwidth]{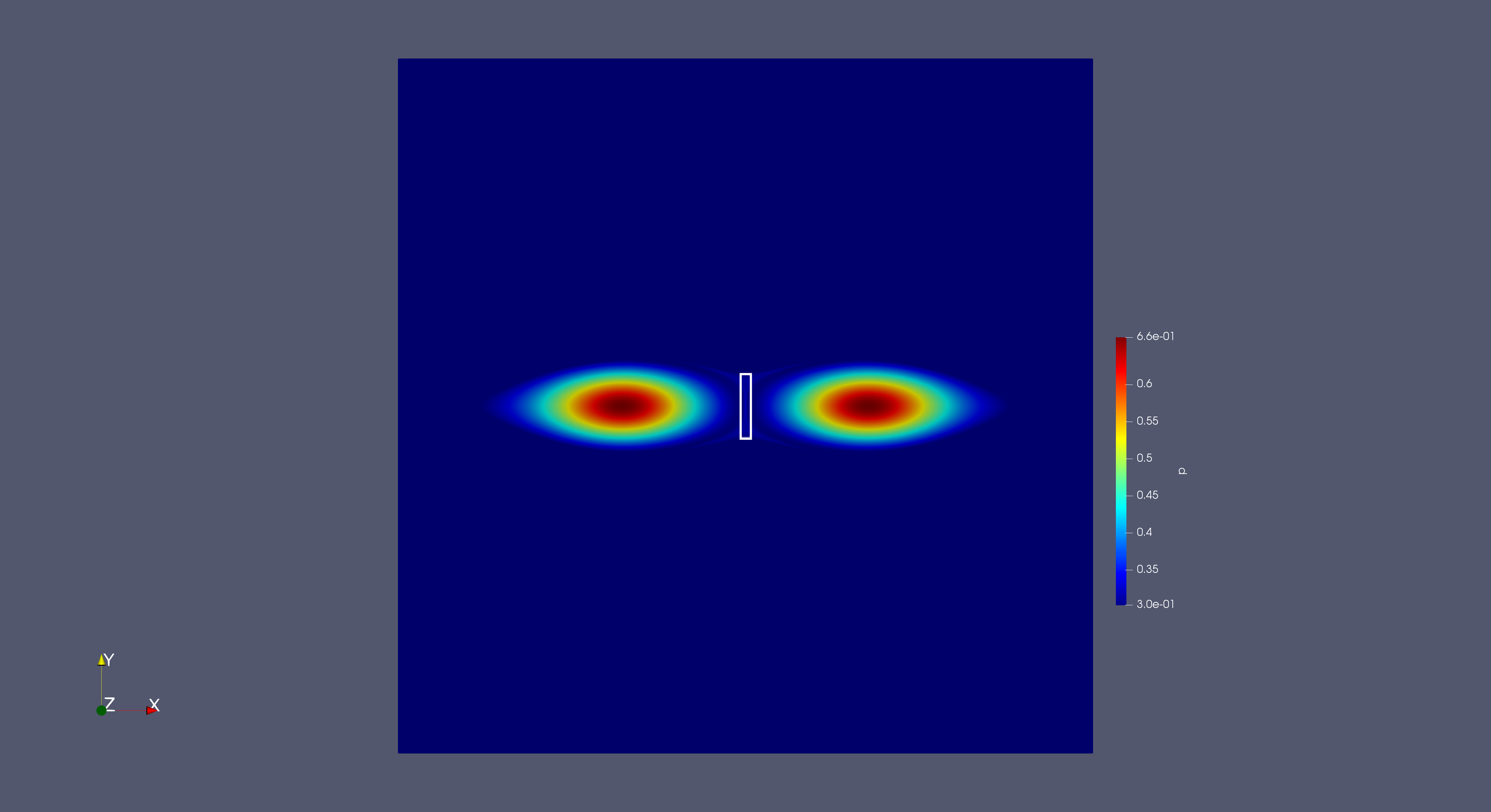}
\caption{2D island coalescence problem: Current (top) and pressure (bottom) at time $t=5.2$ showing the coalescence of the islands. In Figure \figref{J_p_plasmoid} we will focus on the
    box region marked in the central portion of the above figures.} 
\figlab{J_p_t_150}
\end{figure}
\begin{figure}[h!b!t!]
\subfigure[$t=5.85$]{
\includegraphics[trim=72cm 10cm 72cm 10cm,clip=true,width=0.14\columnwidth,height=0.865\columnwidth]{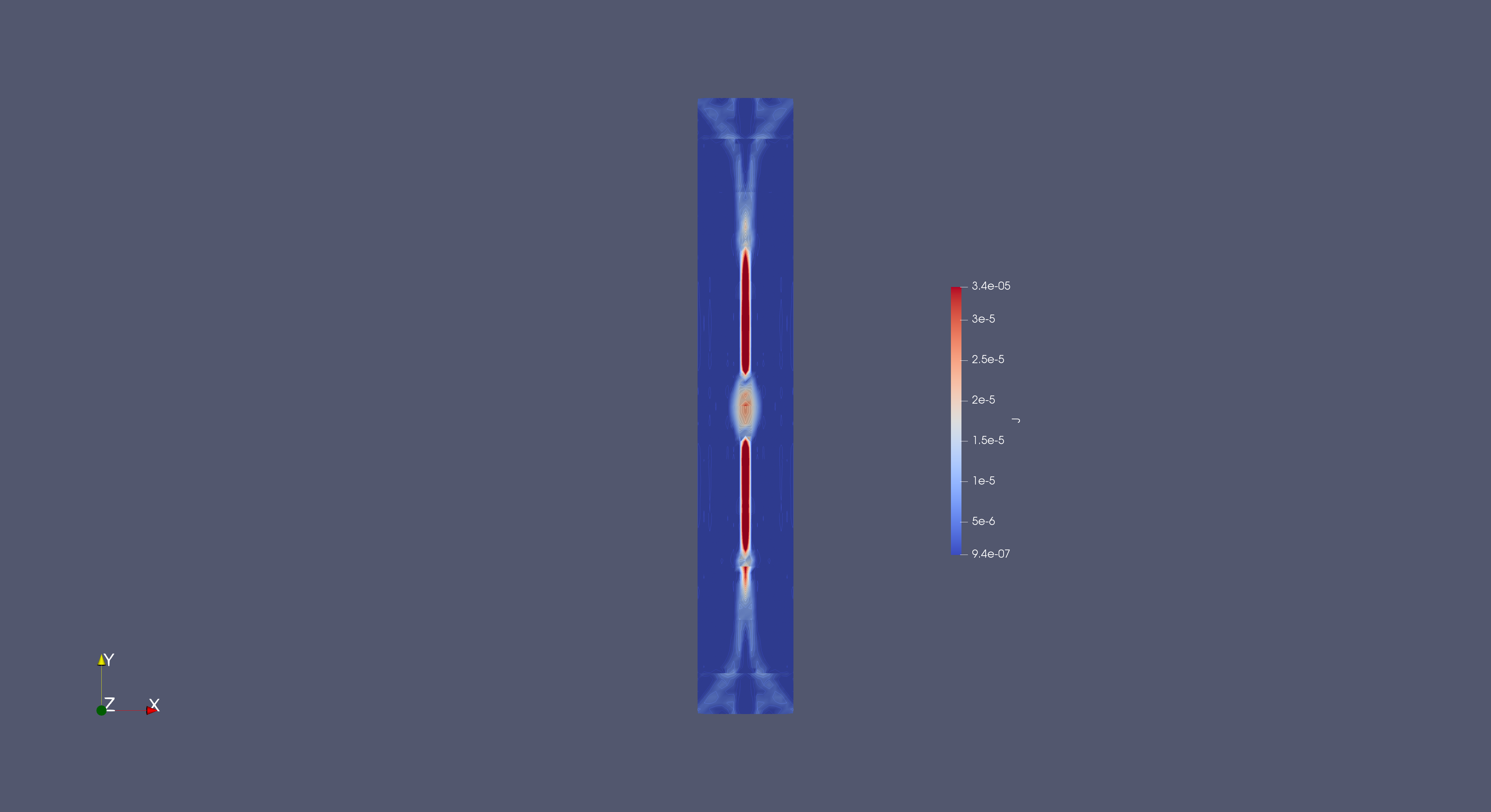}
}
\subfigure[$t=6.22$]{
\includegraphics[trim=72cm 10cm 72cm 10cm,clip=true,width=0.14\columnwidth,height=0.865\columnwidth]{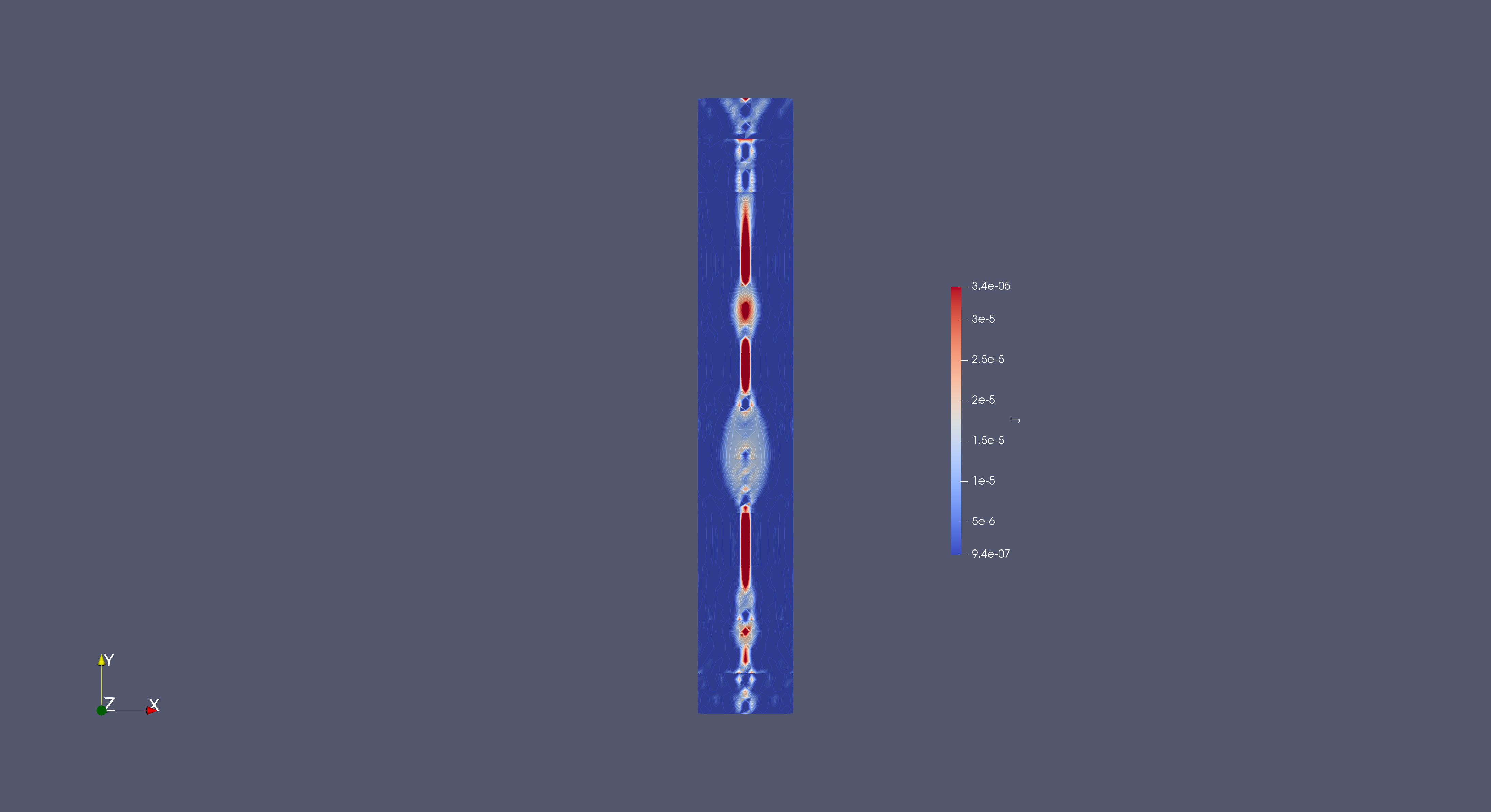}
}
\subfigure[$t=6.3$]{
\includegraphics[trim=72cm 10cm 72cm 10cm,clip=true,width=0.14\columnwidth,height=0.865\columnwidth]{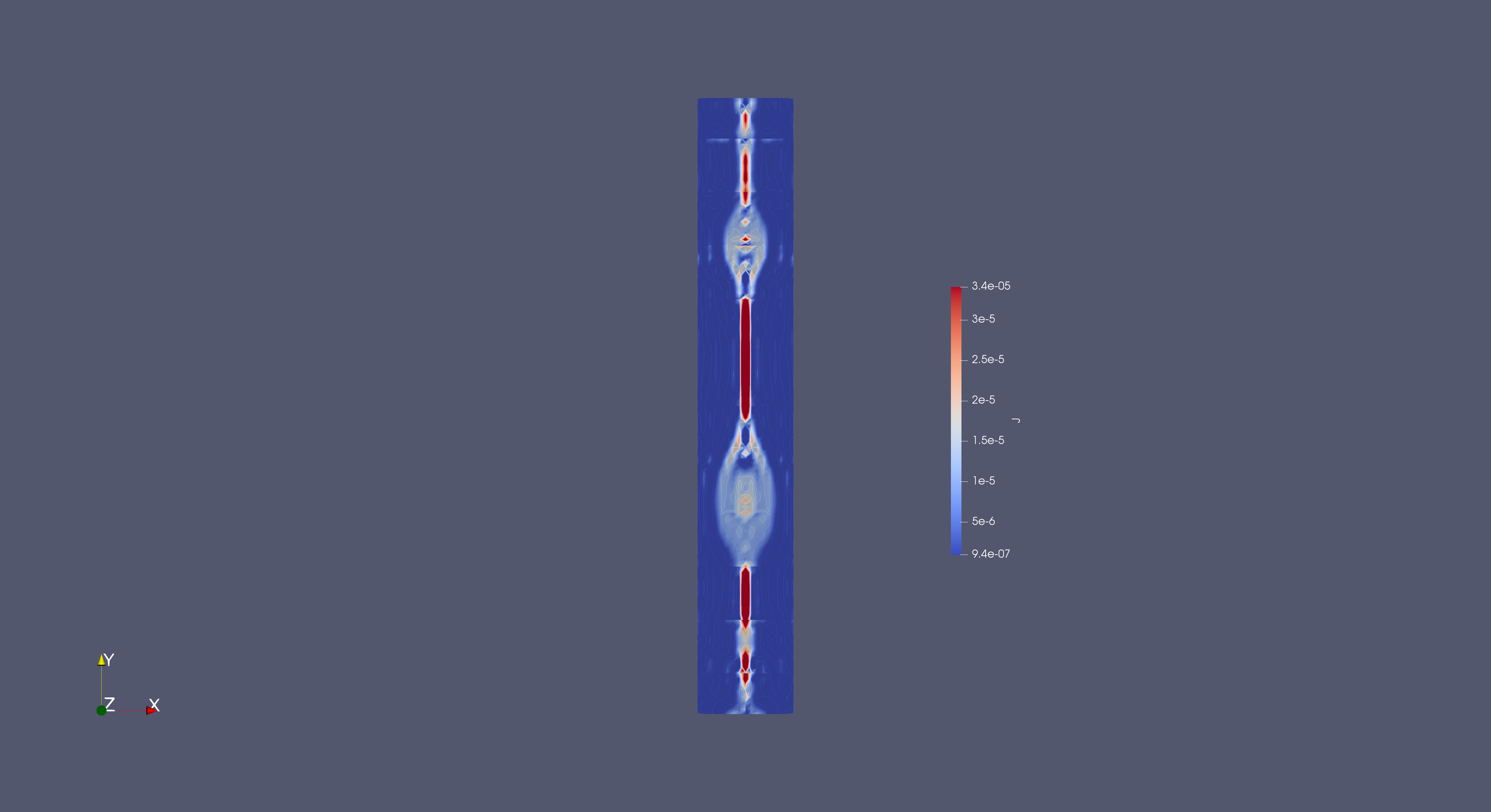}
}
\subfigure[$t=5.85$]{
\includegraphics[trim=71cm 3cm 71cm 3cm,clip=true,width=0.14\columnwidth]{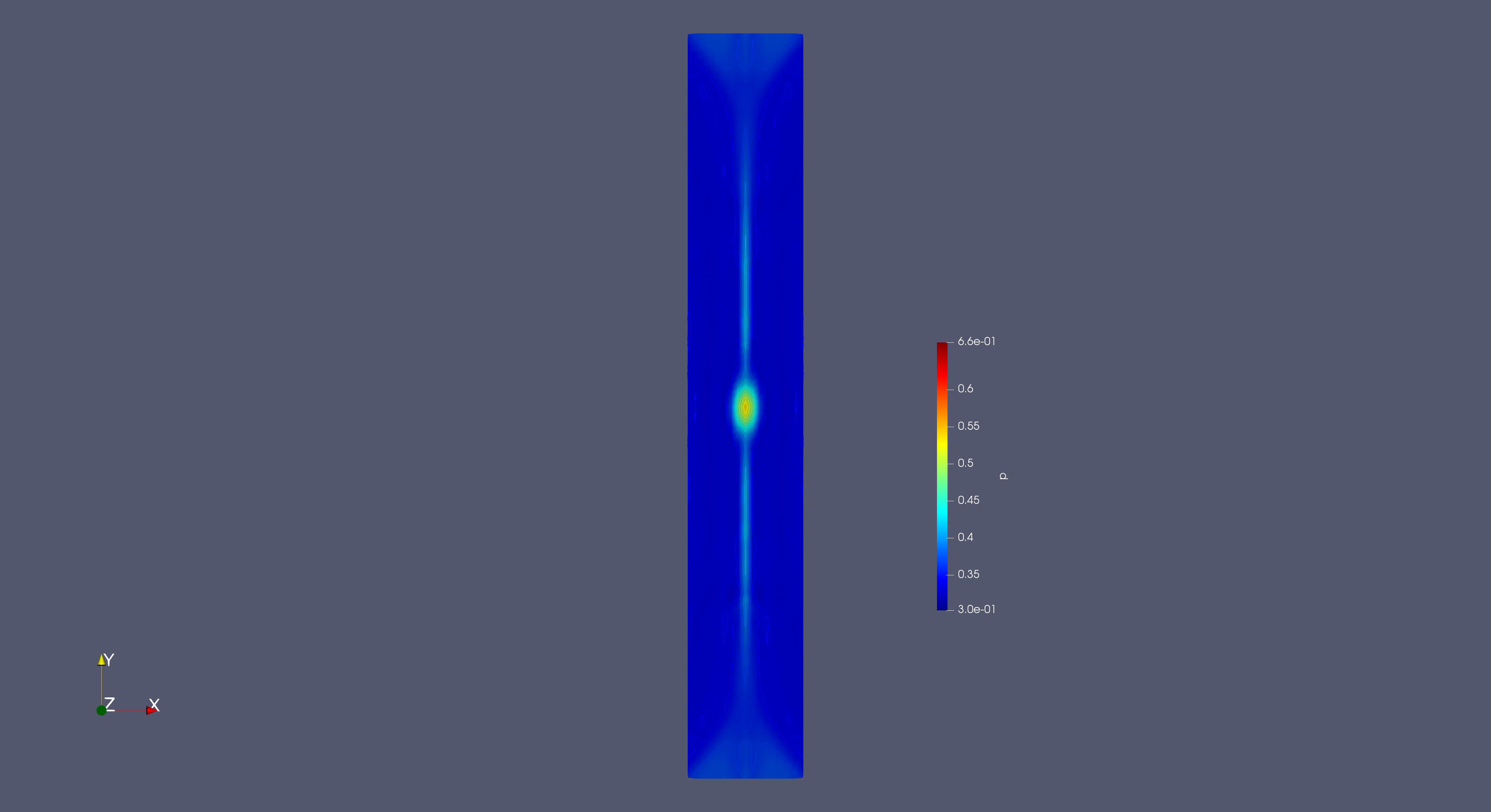}
}
\subfigure[$t=6.22$]{
\includegraphics[trim=71cm 3cm 71cm 3cm,clip=true,width=0.14\columnwidth]{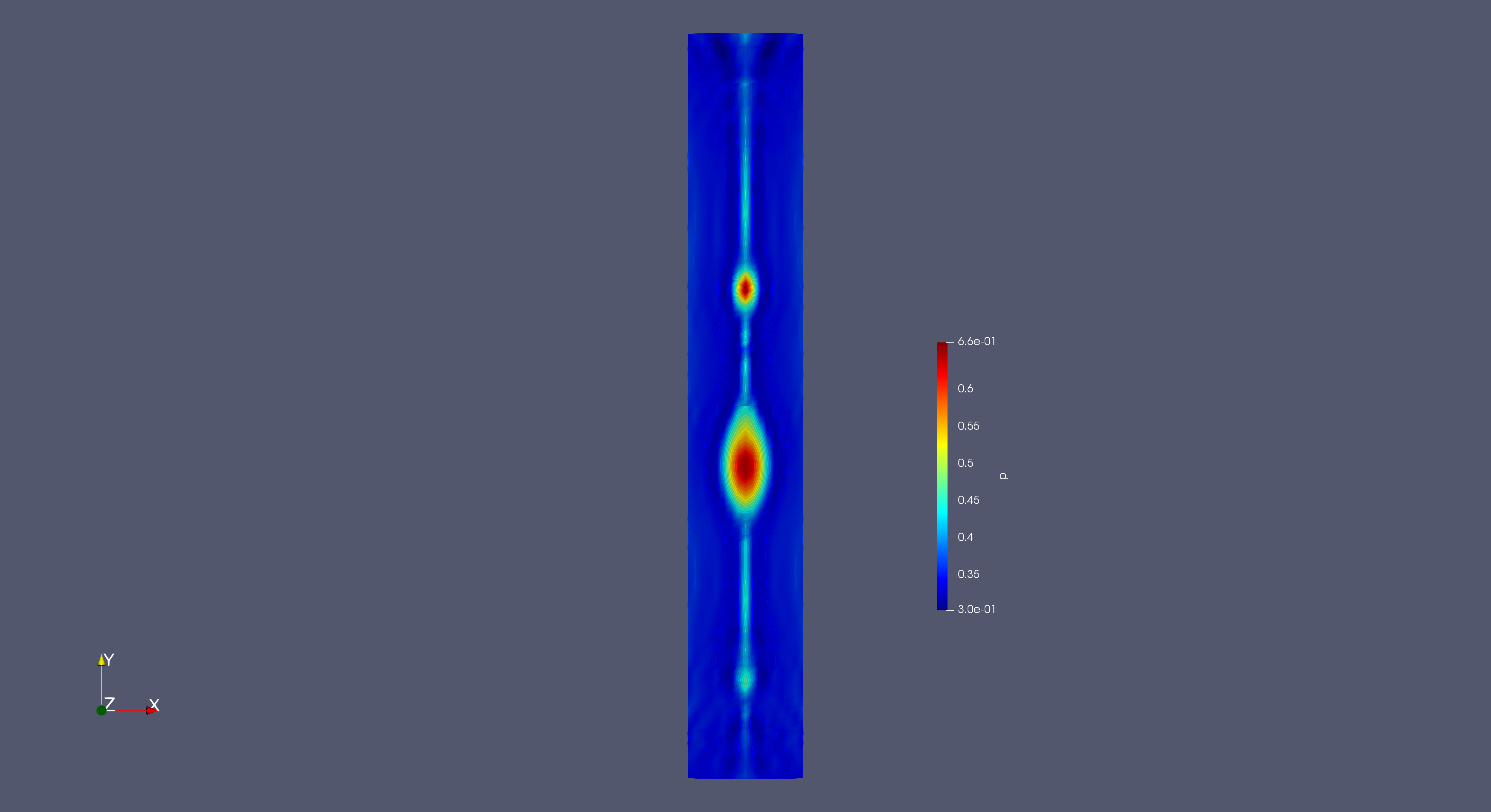}
}
\subfigure[$t=6.3$]{
\includegraphics[trim=71cm 3cm 71cm 3cm,clip=true,width=0.14\columnwidth]{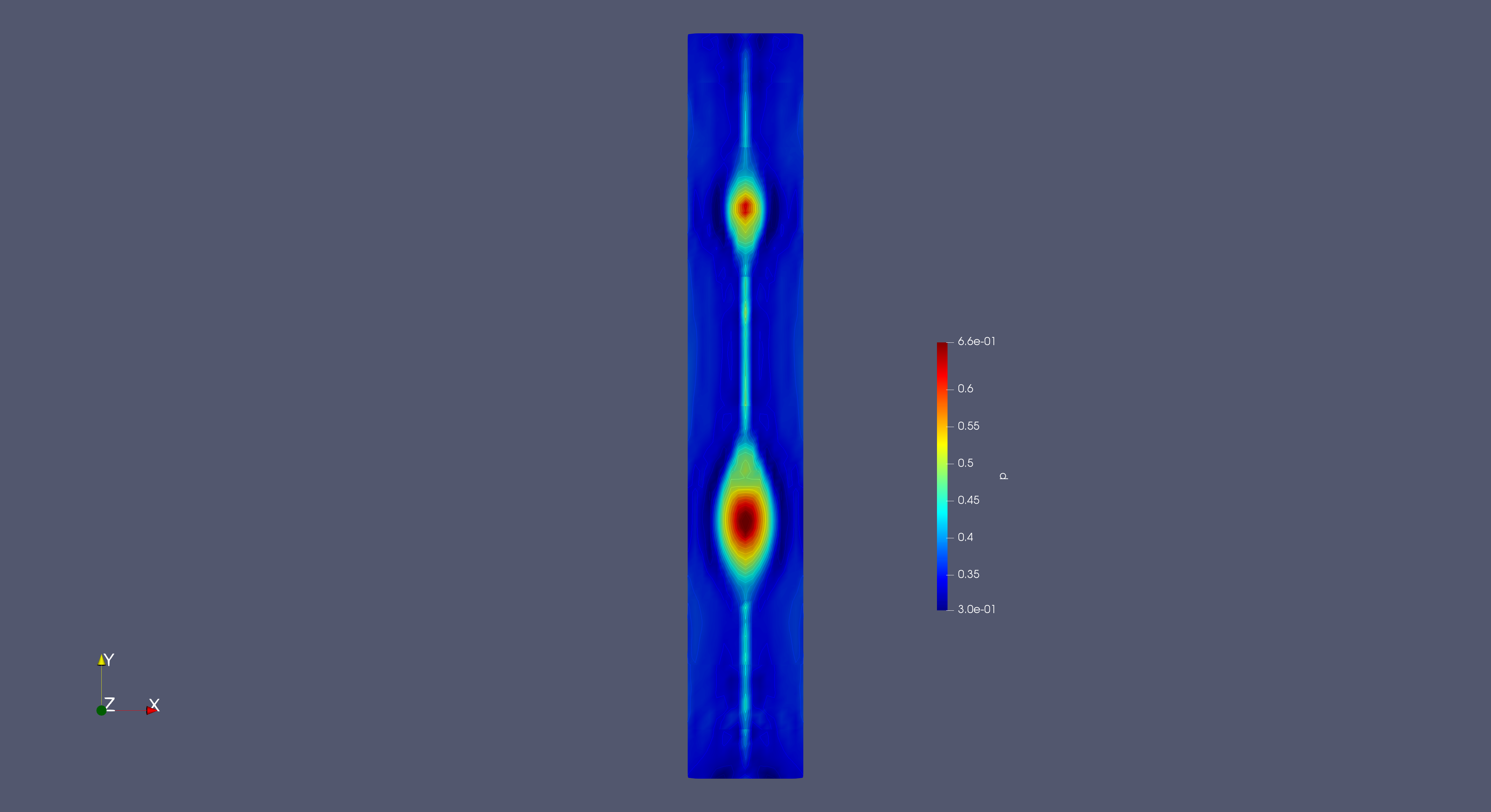}
}
    \caption{2D island coalescence problem: Current (first three from left to right) and pressure (last three) plots at the indicated times showing the formation and evolution of plasmoids with time.}
\figlab{J_p_plasmoid}
\end{figure}

Figure \figref{J_p_t_150} shows the formation of x-structure in the center of the
domain at $t=5.2$ similar to the 3D problem in Figure \figref{3D_x1}. We zoom in to the box marked at the center of the domain in Figure \figref{J_p_t_150} to see clearly the formation and evolution of 
plasmoids. Figure \figref{J_p_plasmoid} shows the current and pressure corresponding to the box region marked in Figure \figref{J_p_t_150}
at different times. From these two figures we can clearly see bubble like structures which are called plasmoids emanating from the breakdown of the thin current 
sheet and moving in the vertical direction. At least from Figure \figref{J_p_plasmoid} we can see two prominent plasmoids one moving 
downwards and the other upwards. However, there are lot more tiny plasmoids which continuously appear and merge as time proceeds. This problem clearly shows the
multiscale nature of the magnetic reconnection phenomenon in both space and time. The island widths are of $\mc{O}(0.1)$, whereas the size of plasmoids are of $\mc{O}(0.01)$.
Similarly the time scale in which the island moves is almost an order of magnitude larger than the plasmoid time scales. The linear solver iterations
for this problem is mostly less than 10 for majority of the Picard steps with some exceptions taking between $15-20$. Thus the BFBT+AMG preconditioner with GMRES smoother
seems to be fairly robust for this challenging problem which involves multiscale physics in space and time. The current plots in Figure
\figref{J_p_plasmoid} show some grid scale structures necessitating more elements in those areas or application of filters/artificial viscosity. We will investigate more into slope limiting, filtering and artificial viscosity strategies for our high-order scheme in the future works.

\subsection{Hydromagnetic Kelvin-Helmholtz instability}
\seclab{HMKH}

In this section we consider the 2D and 3D versions of hydromagnetic Kelvin-Helmholtz (HMKH) instability problem studied in 
\cite{shadid2016scalable,phillips2016block,cyr2013new}. This problem also involves magnetic reconnection similar to the island coalscence
problem in the previous section. The domain we consider is $[0, 4]\times[-2, 2]$ in 2D and 
$[0, 4]\times[-2, 2]\times[0, 2]$ in 3D. The initial conditions consists of two counter flowing conducting fluids of constant
velocities given by $\ub^0(x, y\geq0, z) = (1,0,0)$ and $\ub^0(x, y<0, z) = (-1,0,0)$ and a Harris sheet magnetic field defined by 
$\bb^0(x,y,z) = (B_0 \tanh(y/\delta), 0, 0)$. We choose a zero forcing for both fluid and magnetic equations. Similar to the island coalescence problem, the boundary conditions are periodic in
x- and z-directions. On the top and bottom faces the fluid boundary conditions are same as the island coalescence problem with zero 
normal velocity and zero shear stress and the magnetic field is defined by the Harris sheet in the initial condition. The Lagrange multiplier 
$r$ is set as zero on all the boundaries. We select the following parameters as per \cite{shadid2016scalable}: $\kappa=1$, $\Rey=\Rm=10^{4}$,
$B_0 = 0.3333$ and $\delta=0.1$. These values along with $\rho=\mu_0=1$ gives a super Alfv\'{e}nic 
Mach number of $M_A = u/u_A = 3$ as described in \cite{shadid2016scalable}. If $M_A>1$ then the magnetic field is not strong enough to suppress the instabilities and the shear layer is Kelvin-Helmholtz unstable. Thus the initial disturbances eventually grow to form vortices which roll up and merge as time proceeds.

\begin{figure}[h!b!t!]
\centering
\subfigure[$t=1.51$]{
\includegraphics[trim=0cm 37cm 0cm 23cm,clip=true,width=\textwidth]{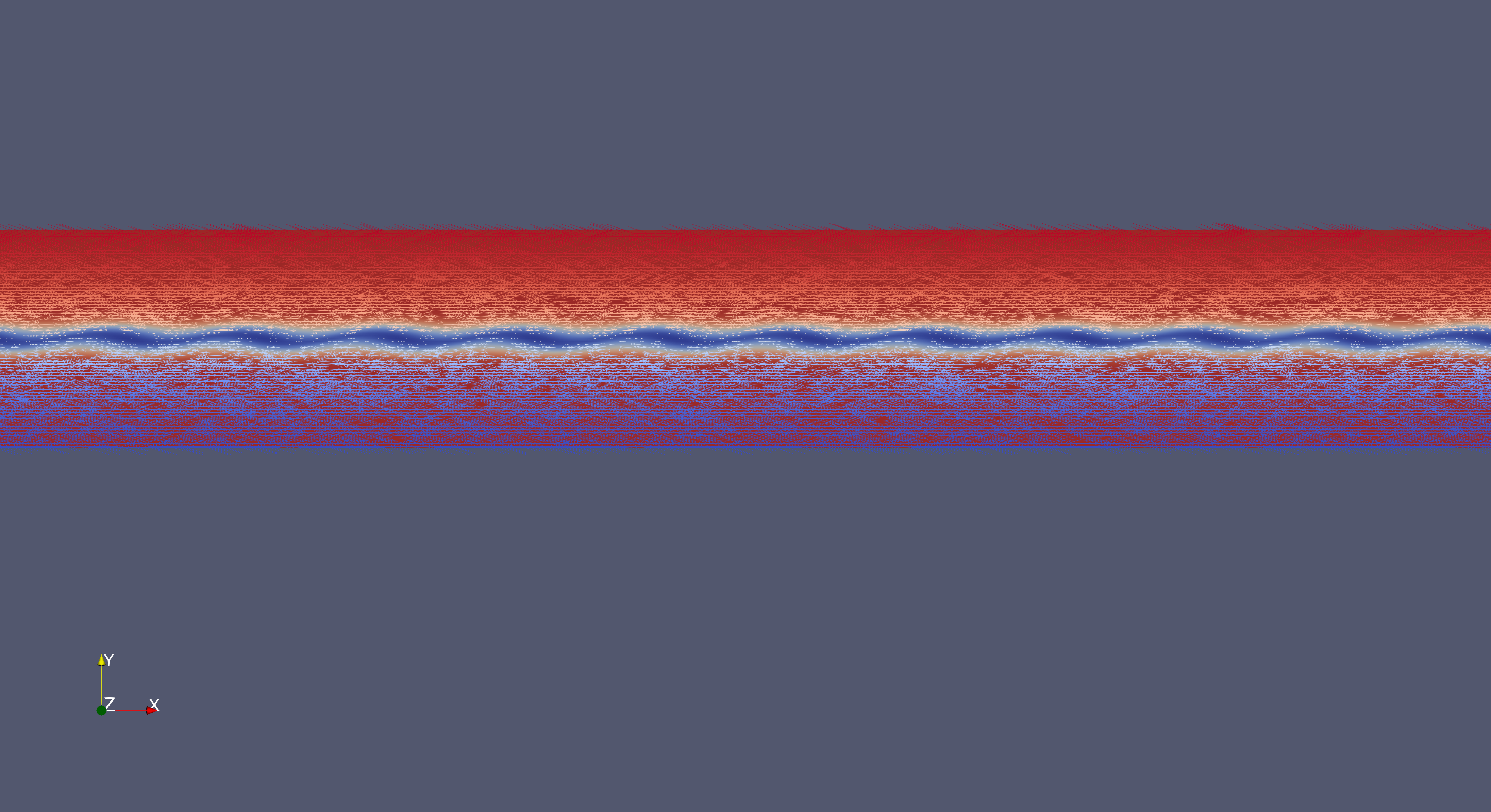}
}
\\
\subfigure[$t=1.73$]{
\includegraphics[trim=0cm 37cm 0cm 23cm,clip=true,width=\textwidth]{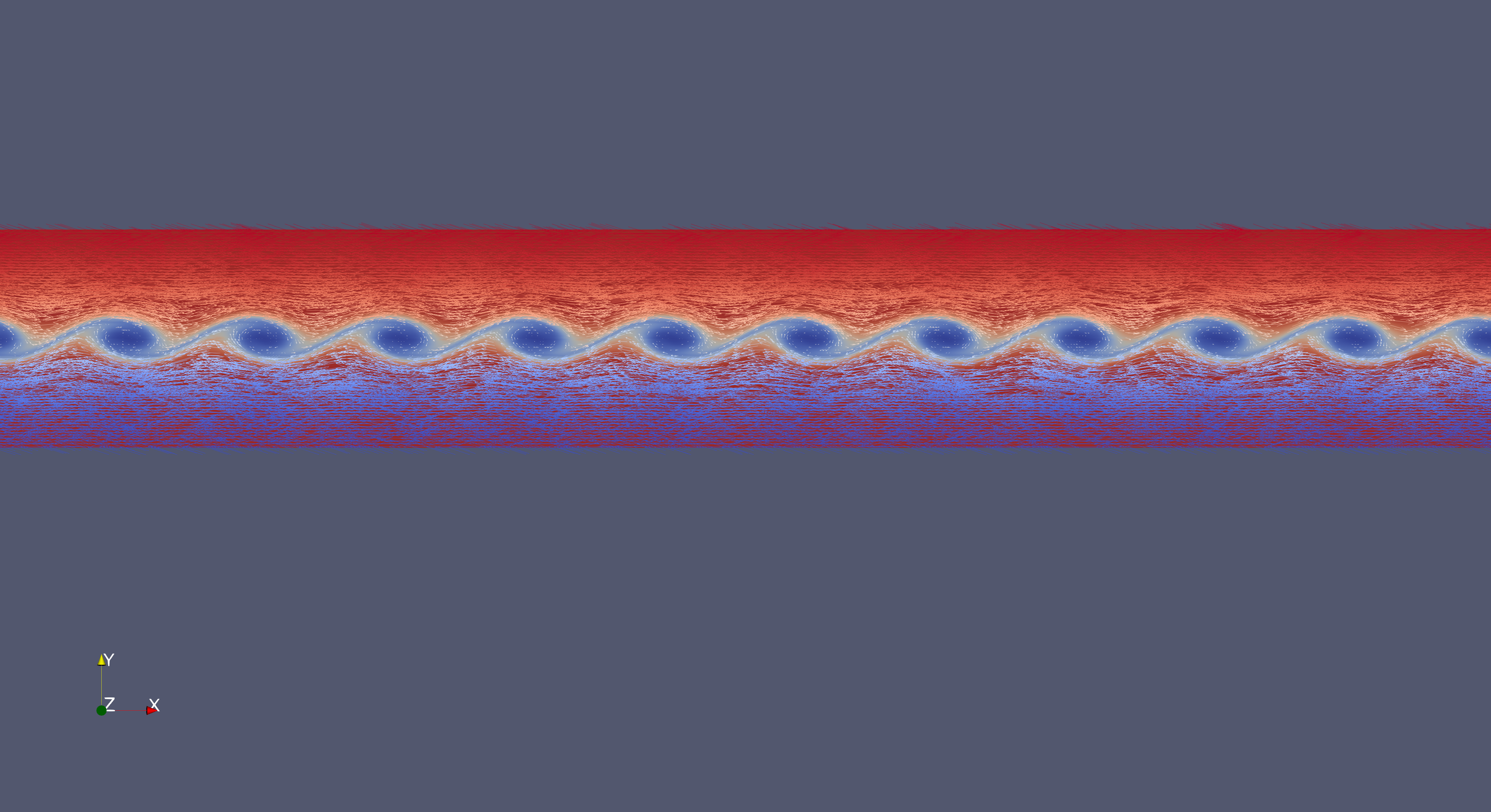}
}
\\
\subfigure[$t=2.61$]{
\includegraphics[trim=0cm 37cm 0cm 23cm,clip=true,width=\textwidth]{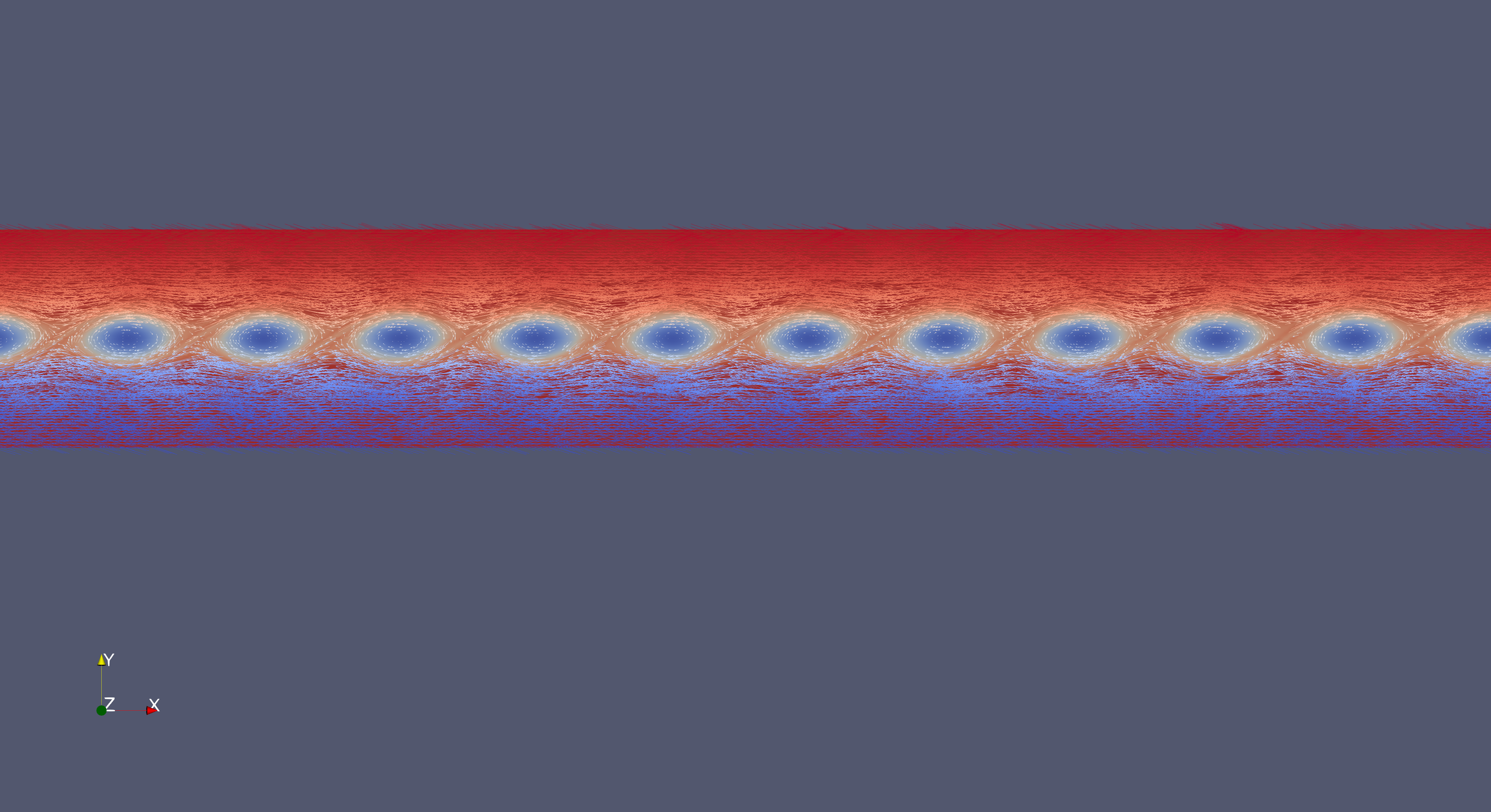}
}
    \caption{2D HMKH problem: Vorticity plots at the indicated times along with the magnetic vectors (marked as arrows). The magnetic vectors are scaled by their magnitude and are colored by the x-component of the magnetic field ($b_x$). The red arrows on the top represent the positive values of $b_x$ and blue arrows on the bottom represent the negative values. 
} 
\figlab{2d_hmkh}
\end{figure}

\begin{figure}[h!b!t!]
\centering
\includegraphics[trim=0cm 0cm 0cm 10cm, clip=true,width=0.8\textwidth]{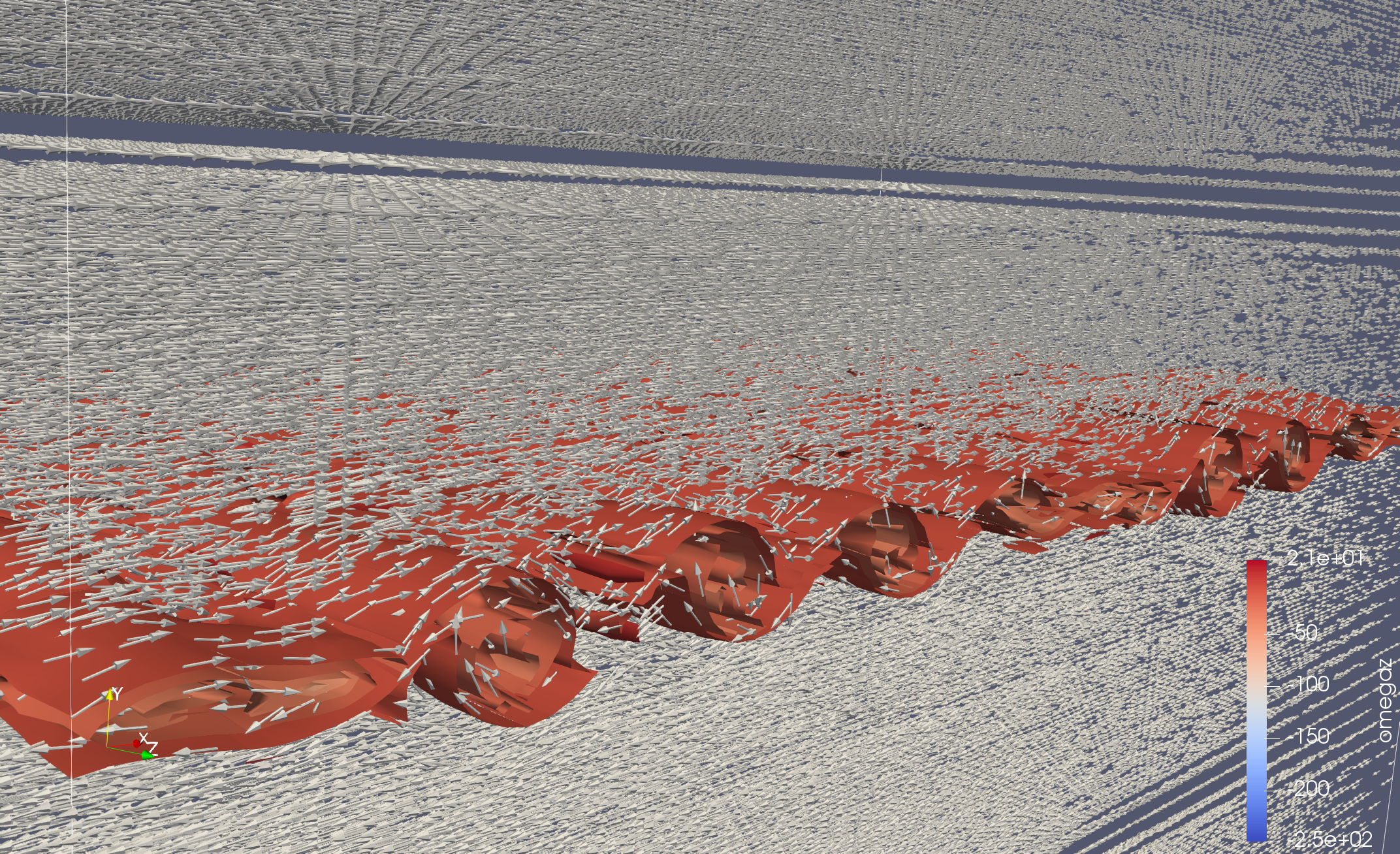}
    \caption{3D HMKH problem: $Z-$component of vorticity contours along with magnetic vectors at time $t=2.0625$. The magnetic vectors are not colored and scaled in this figure to improve visibility.}
\figlab{3d_hmkh}
\end{figure}

First, we consider the 2D HMKH problem discretized in a $128\times128$ mesh with solution order $\p=6$. We use a time stepsize of $\Delta t = 0.001$ in the five stage fourth order DIRK method. 
As mentioned in the beginning of the numerical section we needed a stricter tolerance of $10^{-9}$ in the linear solver to make the Picard iterations converge for this problem up to a 
relative tolerance of $10^{-4}$. The preconditioner is BFBT+AMG with GMRES smoother and the outer iterations are carried out by FGMRES. In Figure \figref{2d_hmkh}, the evolution of
vorticity with time is shown along with the magnetic vectors marked by arrows. The figure shows the roll up of vortices to form the familiar cat-eye pattern and the magnetic vectors bends and follows the 
fluid evolution as time proceeds. 
In Figure \figref{3d_hmkh}, we show the 3D HMKH problem discretized in a $20\times24\times7$ mesh clustered around the region of solution and solution order $\p=5$.
An initial time stepsize of $\Delta t=0.025$ is selected for the backward Euler time stepping and the adaptive time stepping procedure described in the previous section is employed. Here also as in 2D  we see the rollup of vortices and  the magnetic vectors following them.

\begin{figure}[h!t!b!]
  \subfigure[]{
    \includegraphics[trim=2cm 6cm 2cm 7cm,clip=true,width=0.48\columnwidth]{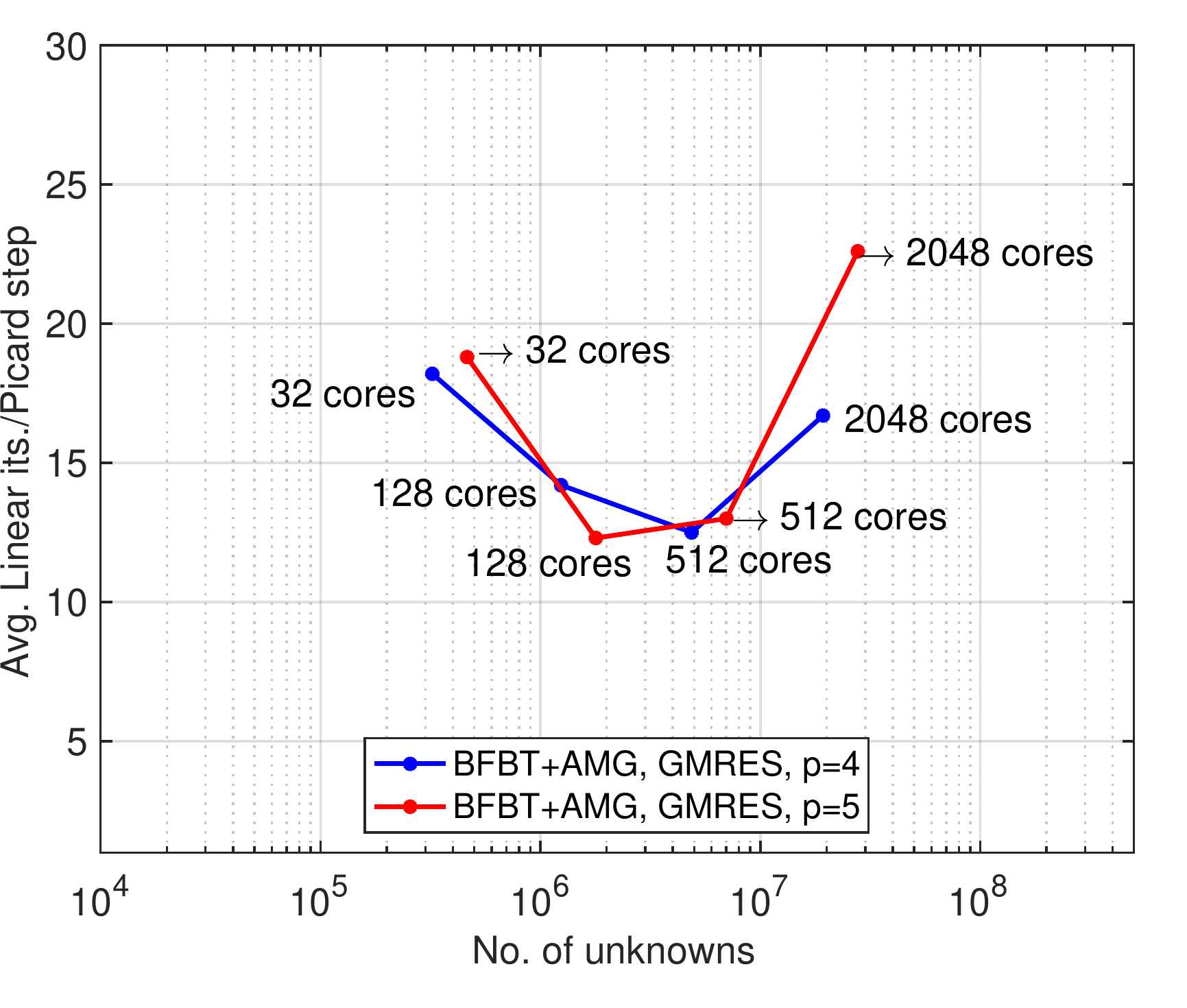}
  }
  \subfigure[]{
    \includegraphics[trim=1.5cm 6cm 2cm 7.5cm,clip=true,width=0.48\columnwidth]{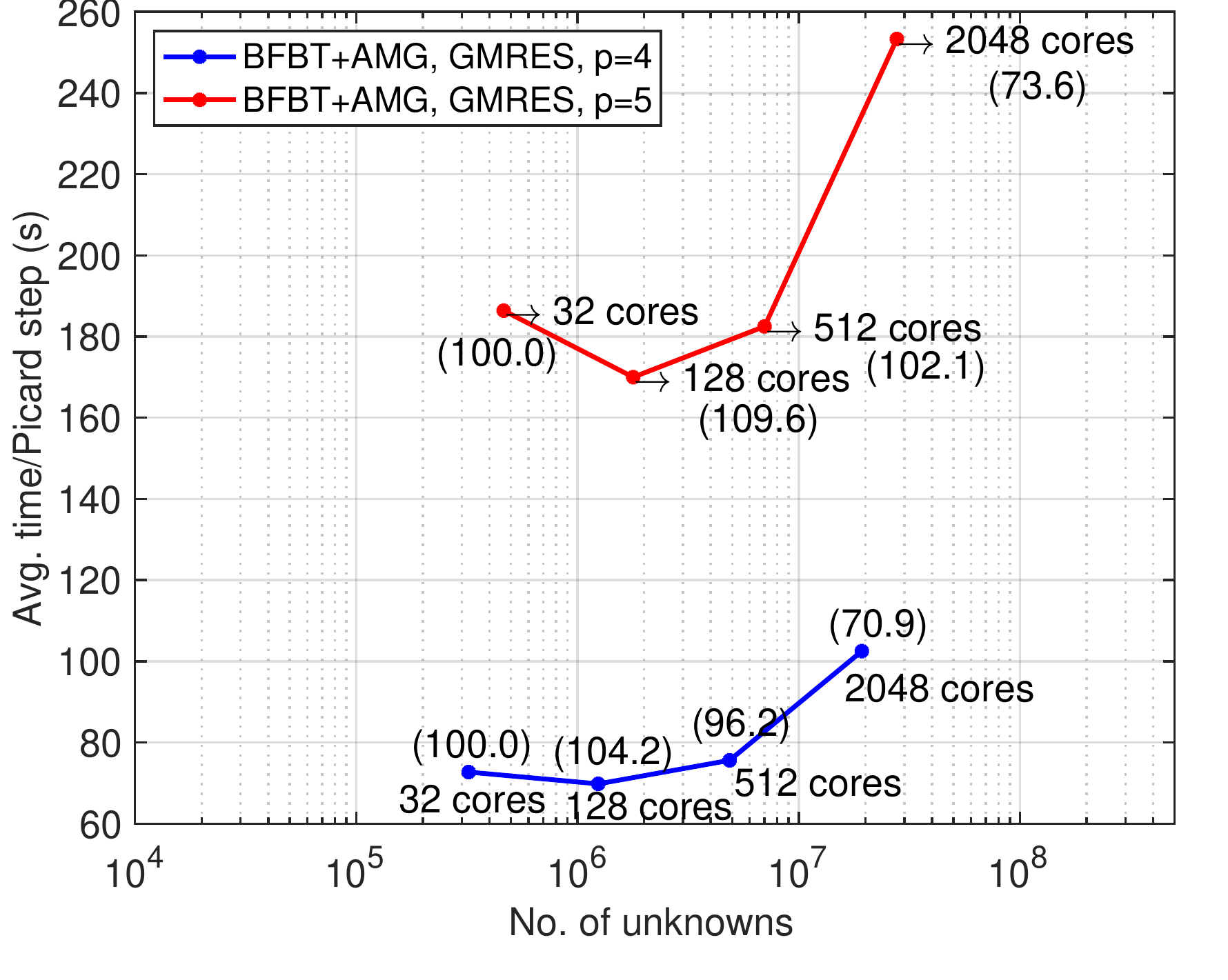}
  }

  \caption{3D HMKH problem. BFBT+AMG with GMRES smoother: weak scaling study of average iterations per Picard step (left) and average time per Picard step (right) for solution orders $\p=4,5$. The values within parentheses in the right figure represent weak scaling parallel efficiencies.}
    \figlab{HMKH_weak_p45}
\end{figure}

Next, we study the weak scaling performance of the BFBT+AMG with GMRES smoother for the 3D HMKH problem. For this we used a fixed time stepsize of $\Delta t=0.01$ in the backward Euler time stepping 
and the results are averaged over six time steps. We used $8\times8\times10$, $16\times16\times10$, $32\times32\times10$ and $64\times64\times10$ meshes corresponding to 32, 128, 512 and 2048 processors respectively.
Thus the number of elements per processor in this case is 20. The CFL ($\p^2\Delta t /h_{min}$) ranges from 0.8 to 2.56 for solution order $\p=4$ and 1.25 to 4 for $\p=5$ corresponding to the mesh sizes and time stepsize. 
The Picard solver took approximately 2 iterations in all the cases. Figure \figref{HMKH_weak_p45} shows the average iteration counts and time per Picard step, since the tolerance of the iterative solver in this
case ($10^{-9}$) is stricter than for the island coalescence problem ($10^{-6}$) we see an increase in overall iteration counts. However, the number of iterations still lies mostly between $10-20$ which is moderate considering the tight tolerance of $10^{-9}$ and the maximum CFL numbers of 2.56 and 4 for $\p=4$ and 5 respectively. The average time per Picard step reflects the trend in the iteration count together with the decrease in scalability 
in the 2048 processors regime as discussed in the previous section on island coalescence.

\subsection{Lid driven cavity}
\seclab{lid}

In this section we consider a hydromagnetic version of the classical lid driven cavity problem. The settings of this problem follow closely \cite{phillips2016block}. Even though we simulated 
the 2D version of this problem also, here we present the results only for the 3D problem for brevity. The domain is $[-0.5, 0.5]^3$, with no slip boundary conditions of $\ub={\bf 0}$
applied on all the walls except the top one where we apply a velocity of $\ub=(1,0,0)$ which drives the flow. For the magnetic field we set the tangential component on each wall as 
$\bb\times\n = (-1,0,0)\times\n$ which acts from right to left. The Lagrange multiplier 
$r$ is set as zero on all the boundaries. Both initial conditions and forcings are chosen as zero. We choose the following parameters: $\kappa=1$, $\Rey=\Rm=1000$ which
corresponds to a Hartmann number of $Ha=\sqrt{\kappa\Rey\Rm}=1000$.

\begin{figure}[h!b!t!]
\centering
\subfigure[$t=0.0125$]{
\includegraphics[trim=30cm 0cm 40cm 6cm,clip=true,width=0.31\columnwidth]{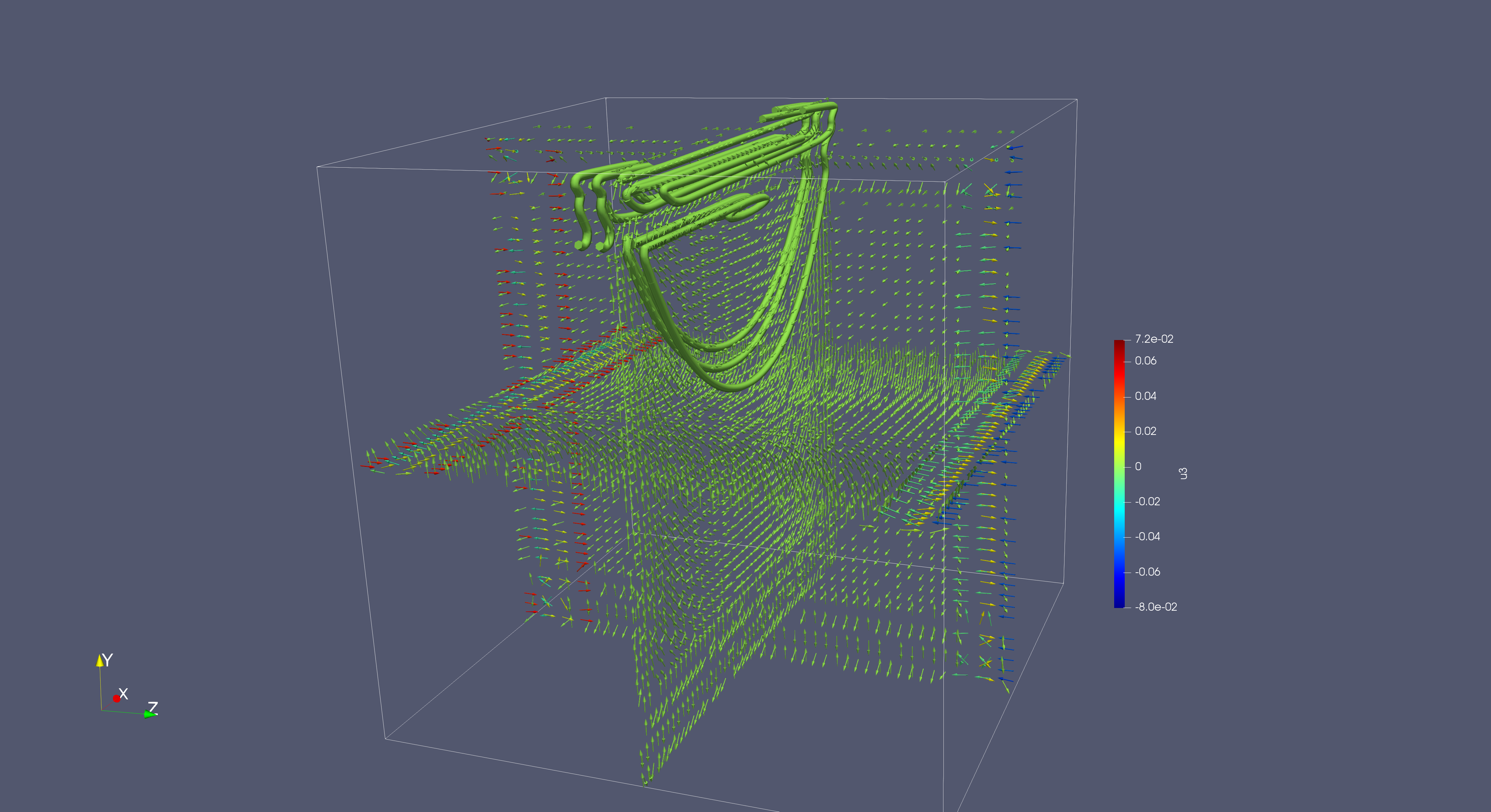}
}
\subfigure[$t=5.075$]{
    \figlab{v_lid1}
\includegraphics[trim=30cm 0cm 40cm 6cm,clip=true,width=0.31\columnwidth]{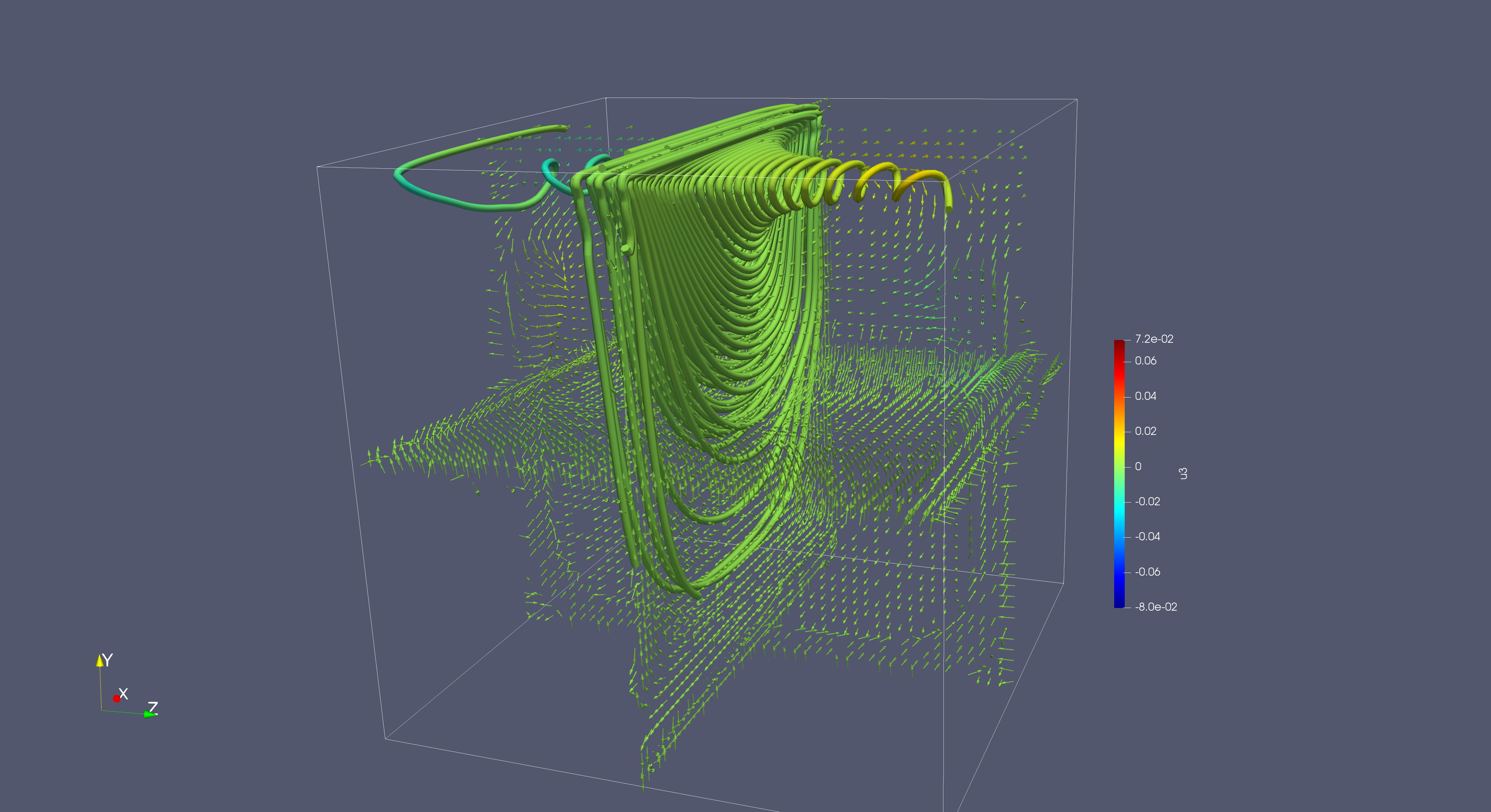}
}
\subfigure[$t=5.6625$]{
    \figlab{v_lid3}
\includegraphics[trim=30cm 0cm 40cm 6cm,clip=true,width=0.31\columnwidth]{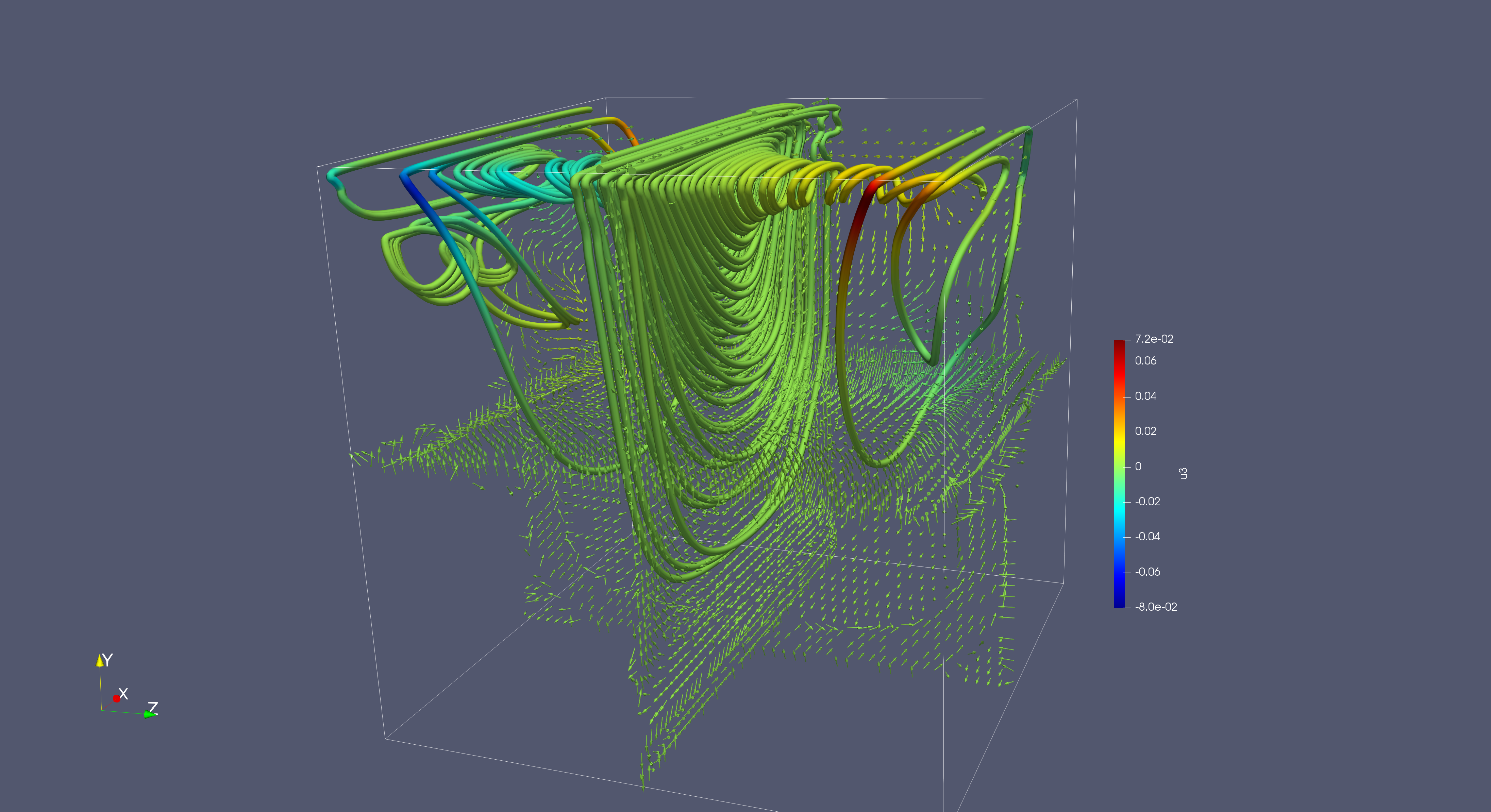}
}
    \caption{3D lid driven cavity problem: Evolution of the streamlines and velocity vectors with time. Both the streamlines and the velocity vectors are colored by the $z-$component of velocity and the velocity vectors are not scaled.}
\figlab{vel_lid}
\end{figure}

\begin{figure}[h!b!t!]
\centering
\subfigure[$t=0.0125$]{
\includegraphics[trim=38cm 0cm 34cm 10cm,clip=true,width=0.31\columnwidth]{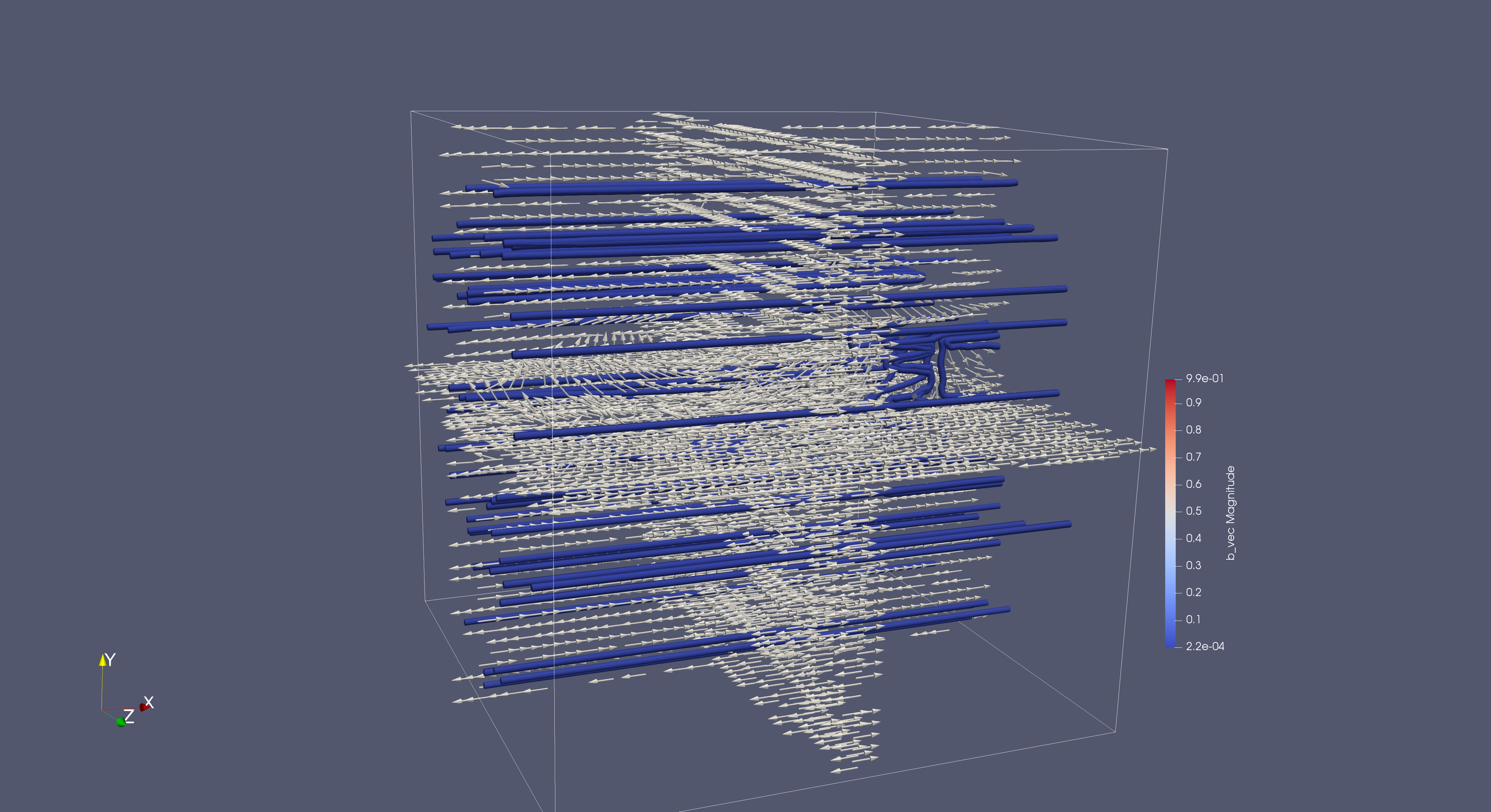}
}
\subfigure[$t=5.075$]{
\includegraphics[trim=38cm 0cm 34cm 10cm,clip=true,width=0.31\columnwidth]{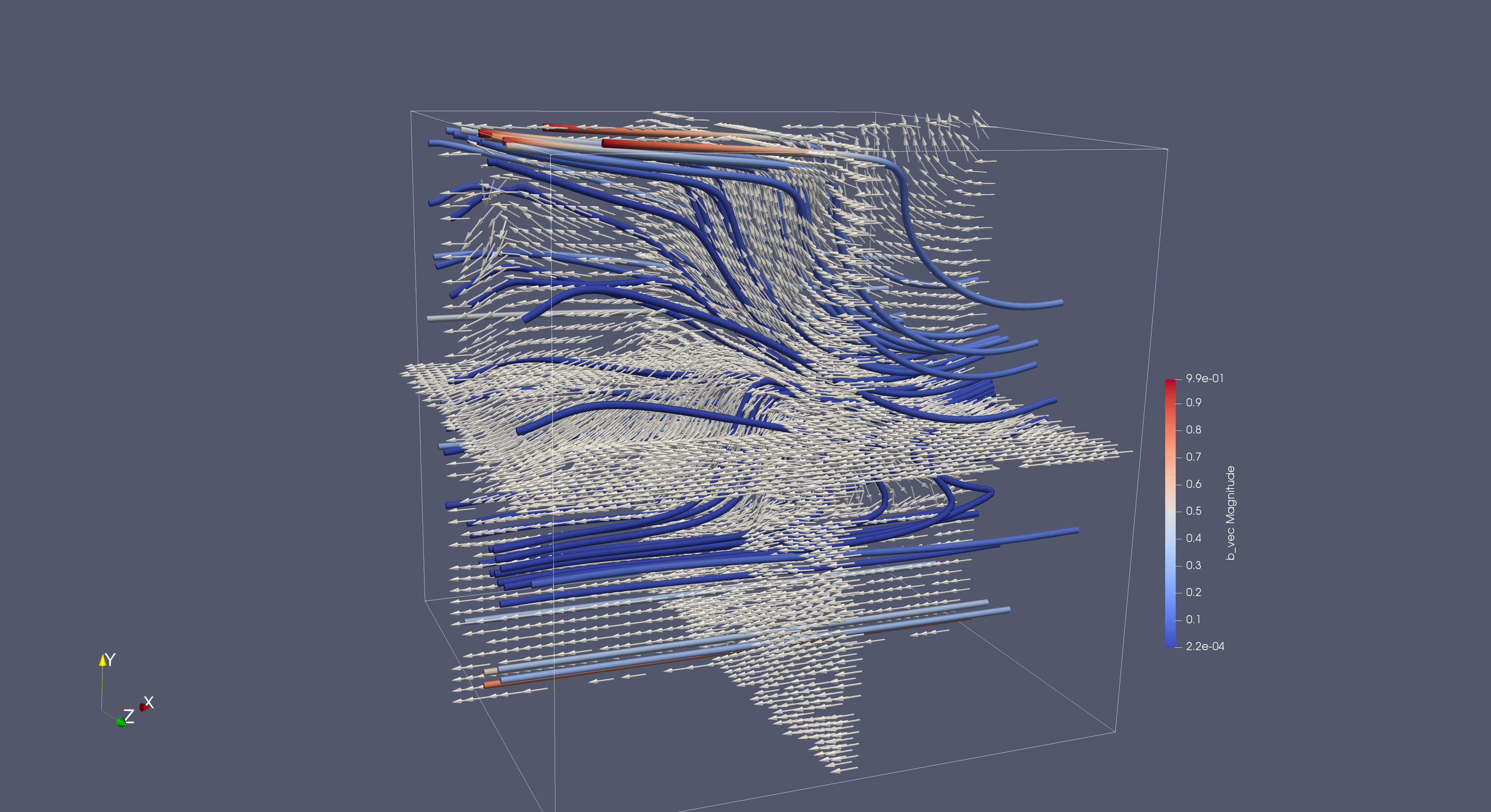}
}
\subfigure[$t=5.6625$]{
\includegraphics[trim=38cm 0cm 34cm 10cm,clip=true,width=0.31\columnwidth]{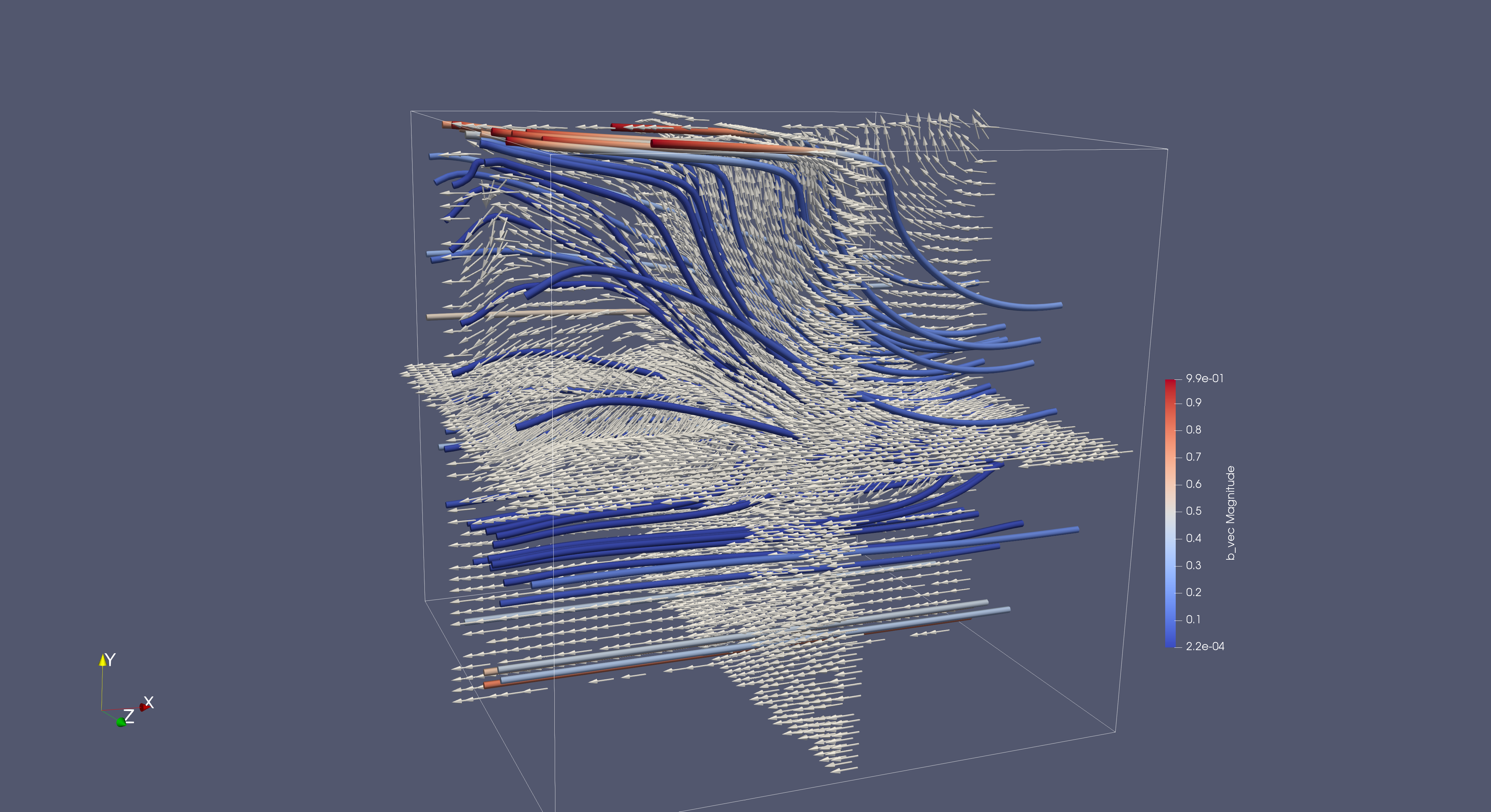}
}
    \caption{3D lid driven cavity problem: Evolution of the magnetic field lines and magnetic vectors with time. The magnetic vectors are not colored and scaled and the field lines are colored by $b_x$.}
\figlab{mag_lid}
\end{figure}

First, we consider a $8^3$ uniform mesh, solution order $\p=6$ and an initial time stepsize of $\Delta t=0.0125$ in the adaptive time stepping with backward Euler method.
The linear and nonlinear solver tolerances and the stopping criteria are same as the ones used for the island coalescence problem. Figure \figref{vel_lid} shows the evolution of the streamlines together with the velocity
vectors in the $x=0$, $y=0$ and $z=0$ planes. The presence of the third dimension allows the streamlines
to curl in the $z-$direction in Figures \figref{v_lid1} and \figref{v_lid3} which gives rise to complex velocity patterns. In Figure \figref{mag_lid} the corresponding magnetic field lines
along with the magnetic vectors are shown. Because of the applied tangential magnetic field from right to left we can see the magnetic lines and the vectors going from right to left with some bending caused by the 
interaction with the fluid components. 

\begin{figure}[h!t!b!]
  \subfigure[]{
      \figlab{iter_weak_lid}
    \includegraphics[trim=2cm 6cm 2cm 7cm,clip=true,width=0.48\columnwidth]{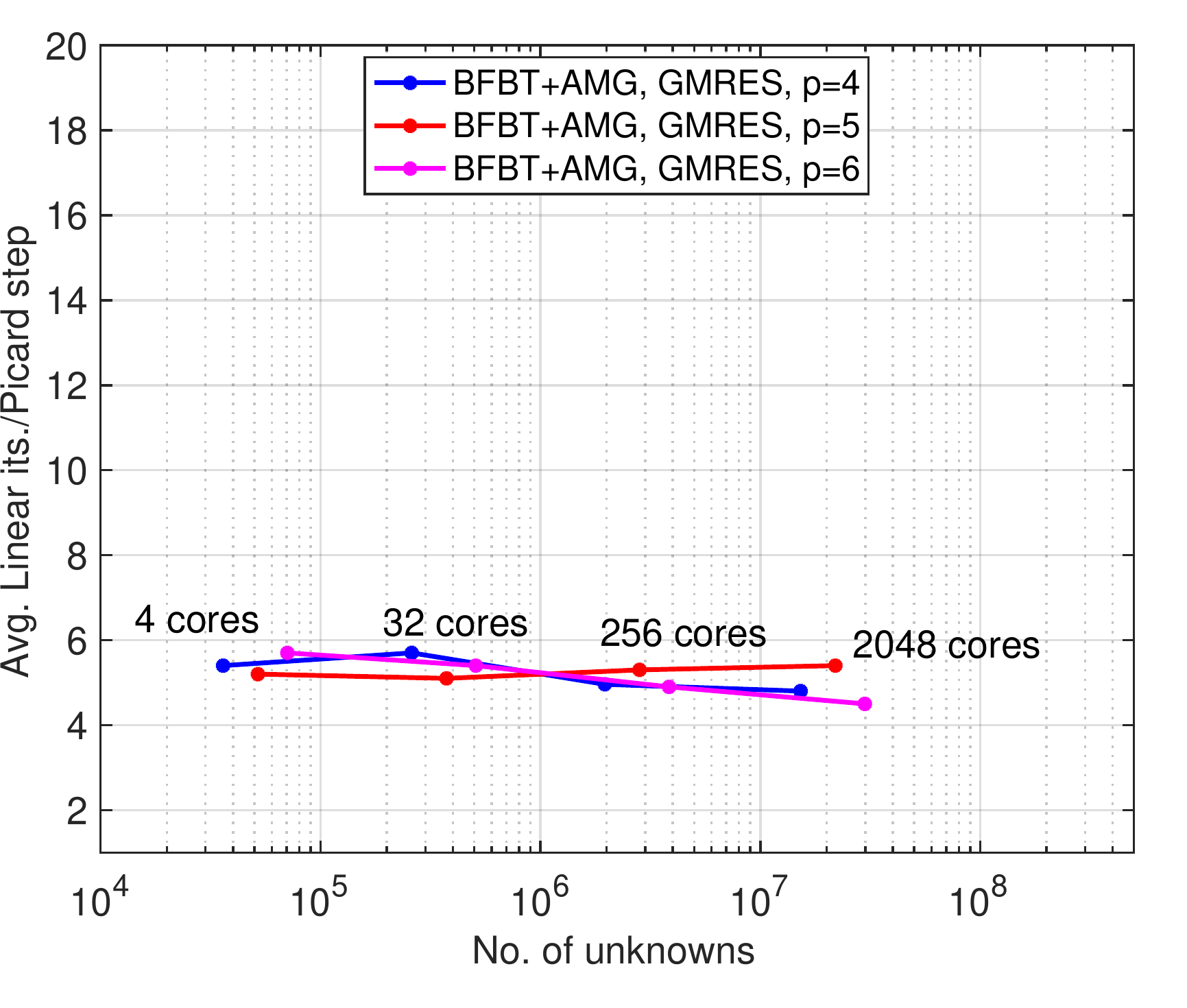}
  }
  \subfigure[]{
      \figlab{time_weak_lid}
    \includegraphics[trim=1.5cm 6cm 2cm 7.5cm,clip=true,width=0.48\columnwidth]{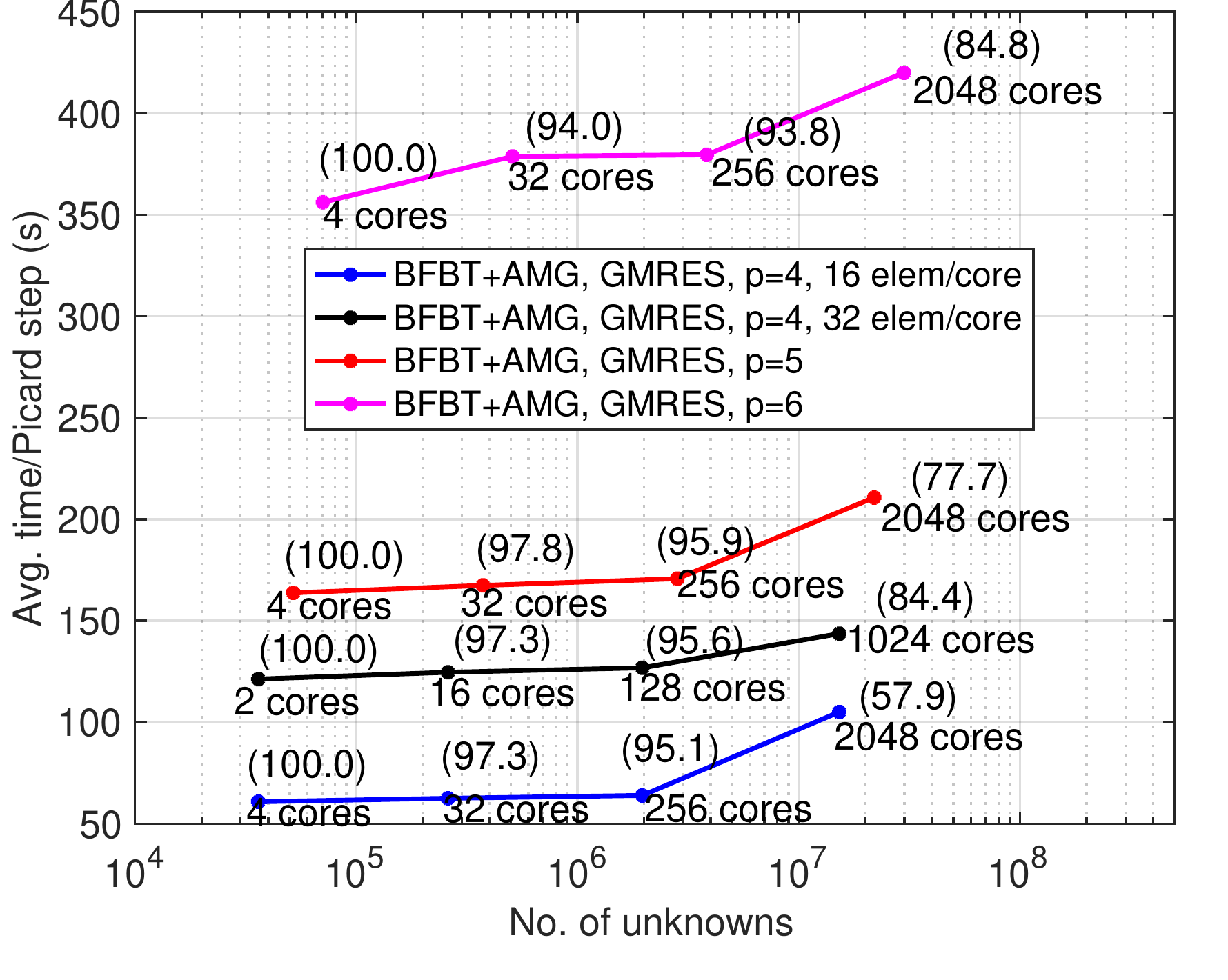}
  }

  \caption{3D lid driven cavity problem. BFBT+AMG with GMRES smoother: weak scaling study of average iterations per Picard step (left) and average time per Picard step (right) for solution orders $\p=4,5,6$. The values within parentheses in the right figure represent weak scaling parallel efficiencies.}
    \figlab{lid_weak_p456}
\end{figure}

For the weak scaling study we consider a fixed time stepsize of $\Delta t =0.0025$ for solution orders $\p=4,5$ whereas for $\p=6$ we consider $\Delta t=0.001$ in the backward Euler method. The meshes 
we consider are $4^3$, $8^3$, $16^3$ and $32^3$ corresponding to 4, 32, 256 and 2048 processors respectively with 16 elements per core. The CFL in this case ranges from 
$0.16-1.28$ for $\p=4$, $0.25-2$ for $\p=5$ and $0.14-1.15$ for $\p=6$. The Picard iterations in all the cases range from $3.5-4.7$. We also tested with $\Delta t=0.01$ for $\p=4$ and observed that the iterations
of the linear solver remains more or less the same as that for $\Delta t=0.0025$ whereas the average Picard iterations are $5.3, 6.5, 7.5$ and $9.3$ corresponding to the four mesh sizes. In Figure \figref{lid_weak_p456} we show the
average iterations and time per Picard step for the BFBT+AMG preconditioner with GMRES smoother and the results are averaged over six time steps. Compared to the previous two experiments, in this case we have a fairly
constant iteration count with mesh refinements for all the solution orders. In Figure \figref{time_weak_lid} we again see some increase in time per Picard step for 2048 cores which in this case comes only from the coarsening in AMG and other 
components of the block preconditioner as we have a flat iteration count. We consider one more case for solution order $\p=4$ with 2, 16, 128 and 1024 cores in Figure \figref{time_weak_lid},
    and this shows better weak scaling performance than the other case with 2048 cores for the finest mesh size. This once again highlights 
    the need for investigation of better coarsening strategies in the AMG preconditioner for better scalability at large number of cores.

\subsection{BFBT+Multilevel preconditioner}
\seclab{amg_vs_ml}
The AMG preconditioner used for the approximation of $\mc{F}^{-1}$ in the block preconditioner contributes to the dominant cost and 
in this section we will consider an alternative to it and compare the performance. To that extent  we will apply the multilevel preconditioner introduced in \cite{muralikrishnan2020multilevel} for $\tilde{\mc{F}}^{-1}$ in the block preconditioning \eqref{eq:BFBT_AMG}. The problem we consider is the 2D island coalescence studied in section \secref{island}. We will also compare the performance of the
iterative solver with respect to nested dissection direct solver in this section.

In \cite{muralikrishnan2020multilevel} the multilevel preconditioner is introduced for problems with scalar trace unknowns and here 
we will extend the idea in a natural manner and apply it to vector valued trace unknowns. Similar to AMG, 
we order the unknowns such that on each edge all the trace unknowns (except 
edge average pressure) corresponding to the first node are ordered first, 
followed by the unknowns in the second node and so on. The ordering of 
unknowns within each nodal point is $(\widehat{U},\widehat{B}^t,\widehat{R})$. 

We then apply one v-cycle of the iterative multilevel algorithm, Algorithm 1 in \cite{muralikrishnan2020multilevel}, for $\tilde{\mc{F}}^{-1}$. 
For the coarse-solver in the multilevel algorithm we use the enriched multilevel approach (EML) because of its robustness and better 
performance compared to the non-enriched version \cite{muralikrishnan2020multilevel}. By means of several numerical experiments we observed that the number of smoothing steps in the block-Jacobi part of the multilevel algorithm, $m_1=0$ and $m_2=1$ i.e., only one post-smoothing, gives the least number
of outer GMRES iterations and we use that in all the cases. We also observed that increasing the number of smoothing steps generally leads to more number of iterations in this case and post-smoothing performs better than pre-smoothing. Since the nodal block ($\mc{F}$) in the HDG 
discretization of MHD is a non-symmetric mixed parabolic-hyperbolic system and for the most part it is purely algebraic we cannot in general expect 
better performance by increasing the number of smoothing steps. The number of iterations also depends on how the fine scale
solver (block-Jacobi) in the multilevel algorithm interacts with the coarse-scale solver (EML) by means of capturing the overall spectrum and
we do not have a clear understanding of this yet. In our future work we will investigate this by means of Fourier analysis
and it can guide us to select the appropriate number of smoothing steps as well as the choice of the fine scale solver.

Having described the specifications of the multilevel algorithm for $\tilde{\mc{F}}^{-1}$ we will now compare BFBT+multilevel preconditioner (referred as BFBT+EML) to BFBT+AMG with GMRES smoother. The outer iterations are carried out with GMRES in the case of BFBT+multilevel preconditioner and FGMRES for
BFBT+AMG with GMRES smoother. For the time discretization we use the backward Euler time stepping with a fixed 
time stepsize of $\Delta t=0.1$. The results are averaged over six time steps for all the cases except $128\times128$ and
$256\times256$ meshes where we average over three time steps. We choose a Lundquist number of $S=10^3$. We run the AMG solver serially in this case to compare with the
serial implementation of the multilevel method. So this
study is mainly to assess the algorithmic scalability of multilevel and AMG preconditioners with respect to mesh refinement, solution
order and Lundquist number.  

\begin{figure}[h!t!b!]
  \subfigure[]{
      \figlab{iter_weak_EML_vs_AMG}
    \includegraphics[trim=3.3cm 8cm 4.5cm 9cm,clip=true,width=0.48\columnwidth]{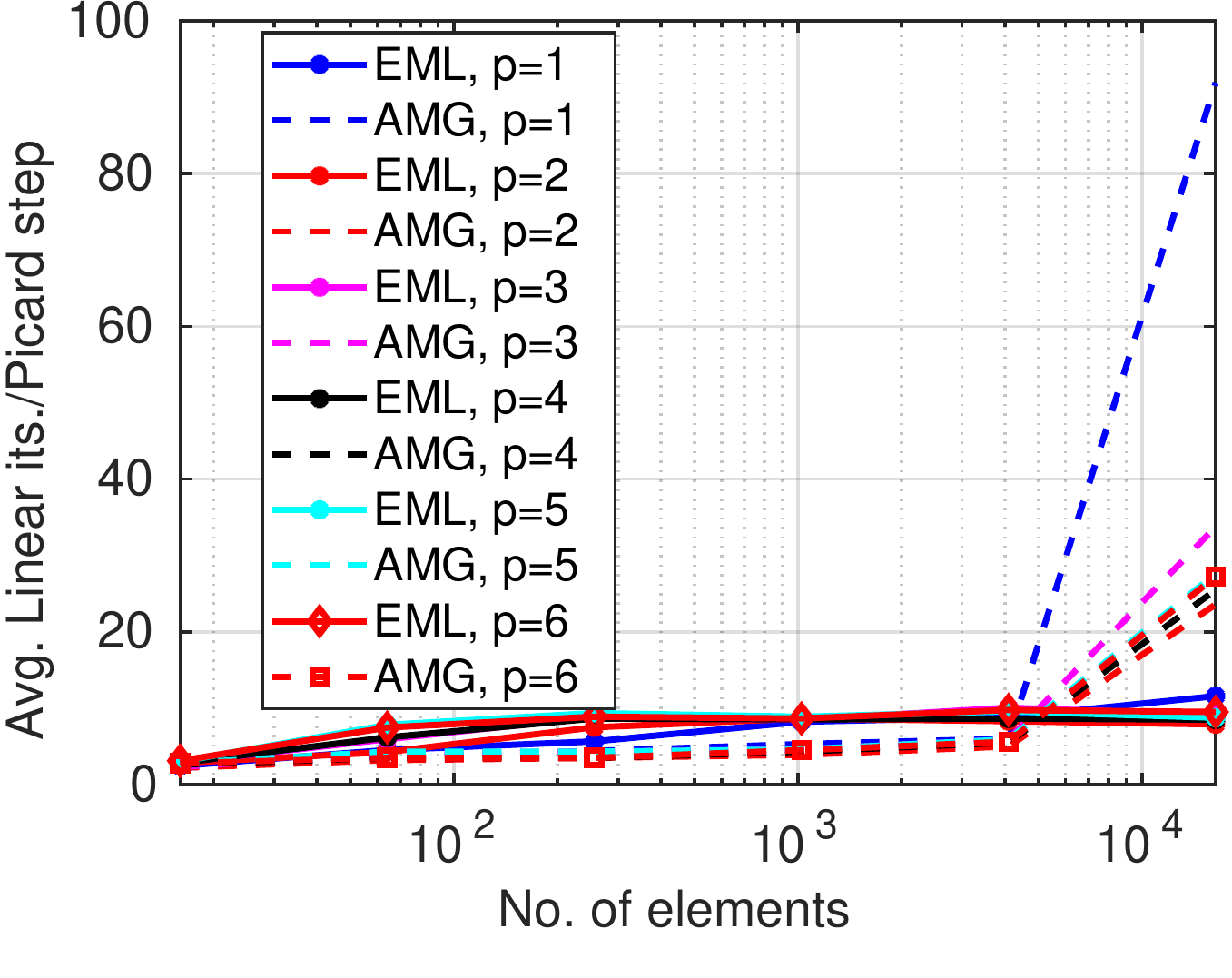}
  }
  \subfigure[]{
      \figlab{time_weak_EML_vs_AMG}
    \includegraphics[trim=3.5cm 8cm 4.5cm 9cm,clip=true,width=0.48\columnwidth]{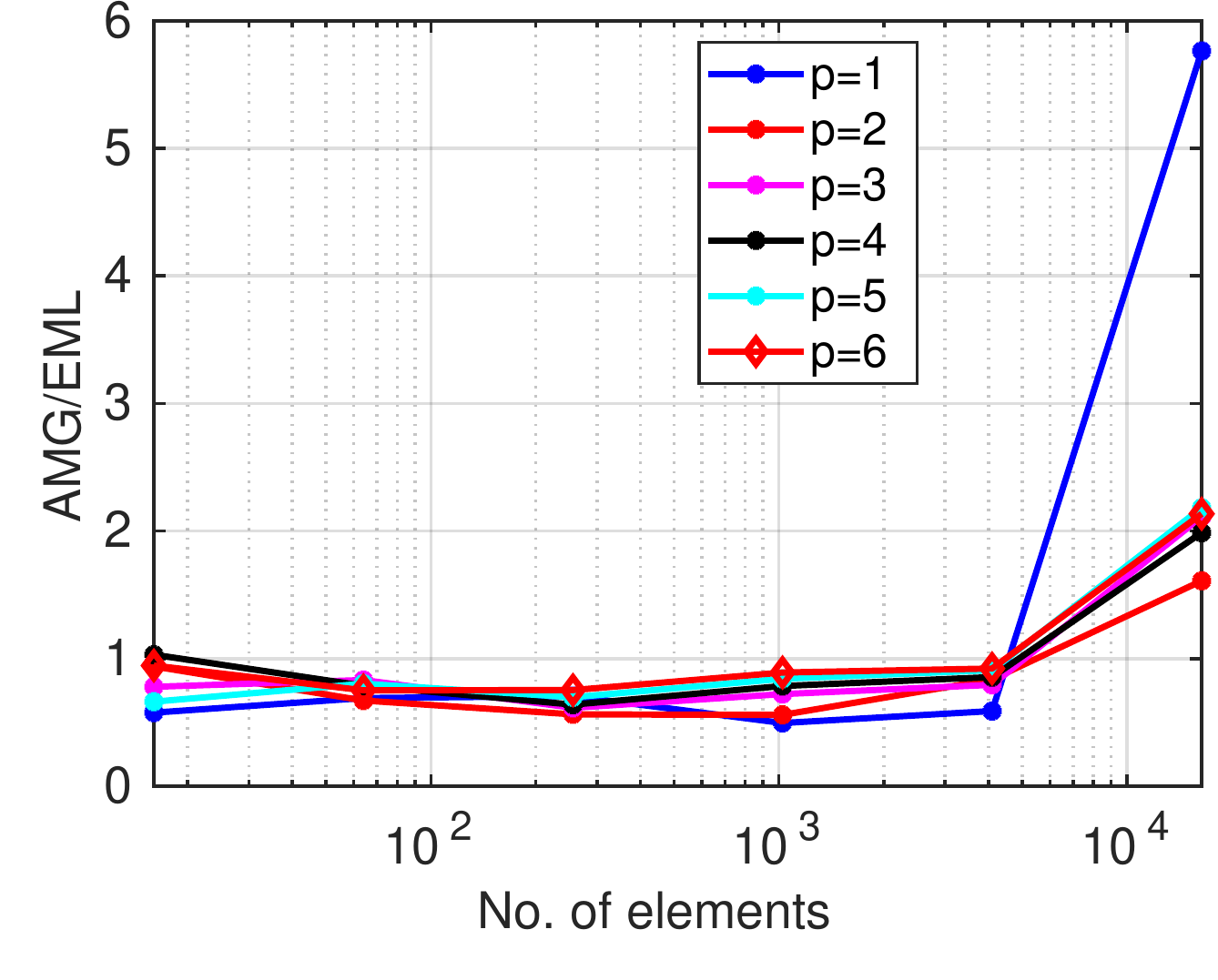}
  }

    \caption{2D island coalescence problem: Comparison of AMG and multilevel (EML) preconditioners as the mesh is refined in $h$ and $\p$. In both cases we use the BFBT approximation for the inverse of the Schur complement. On the left is the number of average iterations per Picard step 
    and on the right is the ratio of average time taken per Picard step for BFBT+AMG over BFBT+EML preconditioner.}
    \figlab{weak_EML_vs_AMG}
\end{figure}
  
  Figure \figref{iter_weak_EML_vs_AMG} shows the average number of iterations for both the AMG and the multilevel (EML) preconditioners together with the BFBT approximation for the inverse Schur complement in the block preconditioning \eqref{eq:BFBT_AMG} as the mesh and solution order are refined. In 
  calculating this average we have omitted the iteration counts for first Picard step in the first time step. This is because the iteration
  counts for both EML and AMG solvers are higher for this case than the rest of the steps as we start from a zero initial guess. Hence it is not a representative of the iteration counts taken in other time steps. Both the algorithms show almost flat iteration count until $64\times64$ mesh with AMG taking slightly less iteration counts
  than EML. However, for $128\times128$ mesh (last marker) and all solution orders AMG shows a sudden increase in iteration 
  count whereas the performance of EML almost remains the same maintaining the algorithmic scalability. Thus the
  EML solver is more scalable than AMG in terms of mesh refinements in the settings of this experiment. Figure
  \figref{time_weak_EML_vs_AMG} shows the ratio of average time taken per Picard step for the BFBT+AMG preconditioner over BFBT+EML and it reflects the trend observed in the iteration counts in Figure \figref{iter_weak_EML_vs_AMG}. We get approximately $2-6$ times speedup with 
  the BFBT+EML preconditioner over BFBT+AMG preconditioner for $128\times128$ mesh.

\begin{figure}[h!t!b!]
  \subfigure[]{
    \includegraphics[trim=1cm 6cm 2cm 5cm,clip=true,width=0.48\columnwidth]{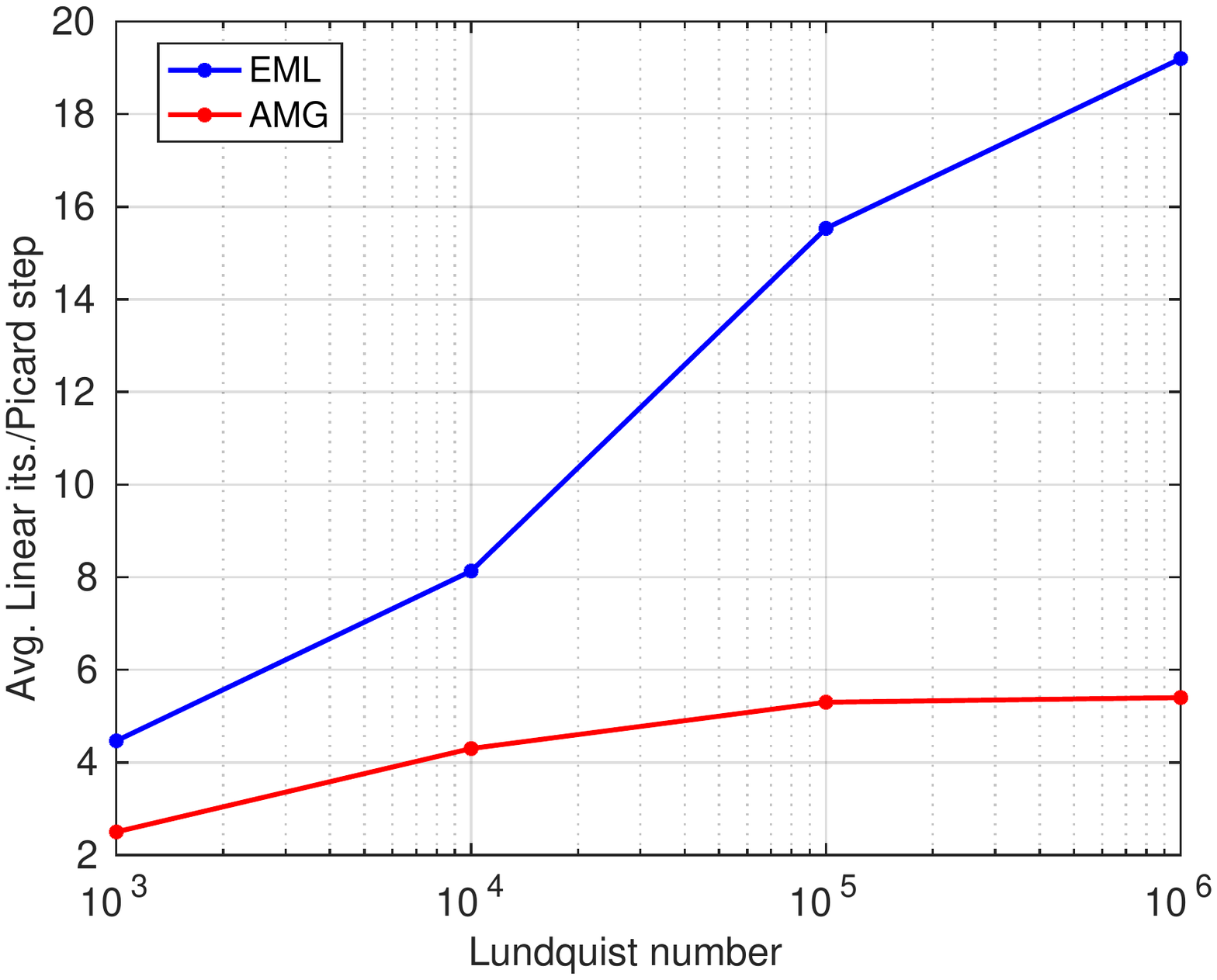}
  }
  \subfigure[]{
    \includegraphics[trim=1cm 6cm 2cm 5cm,clip=true,width=0.48\columnwidth]{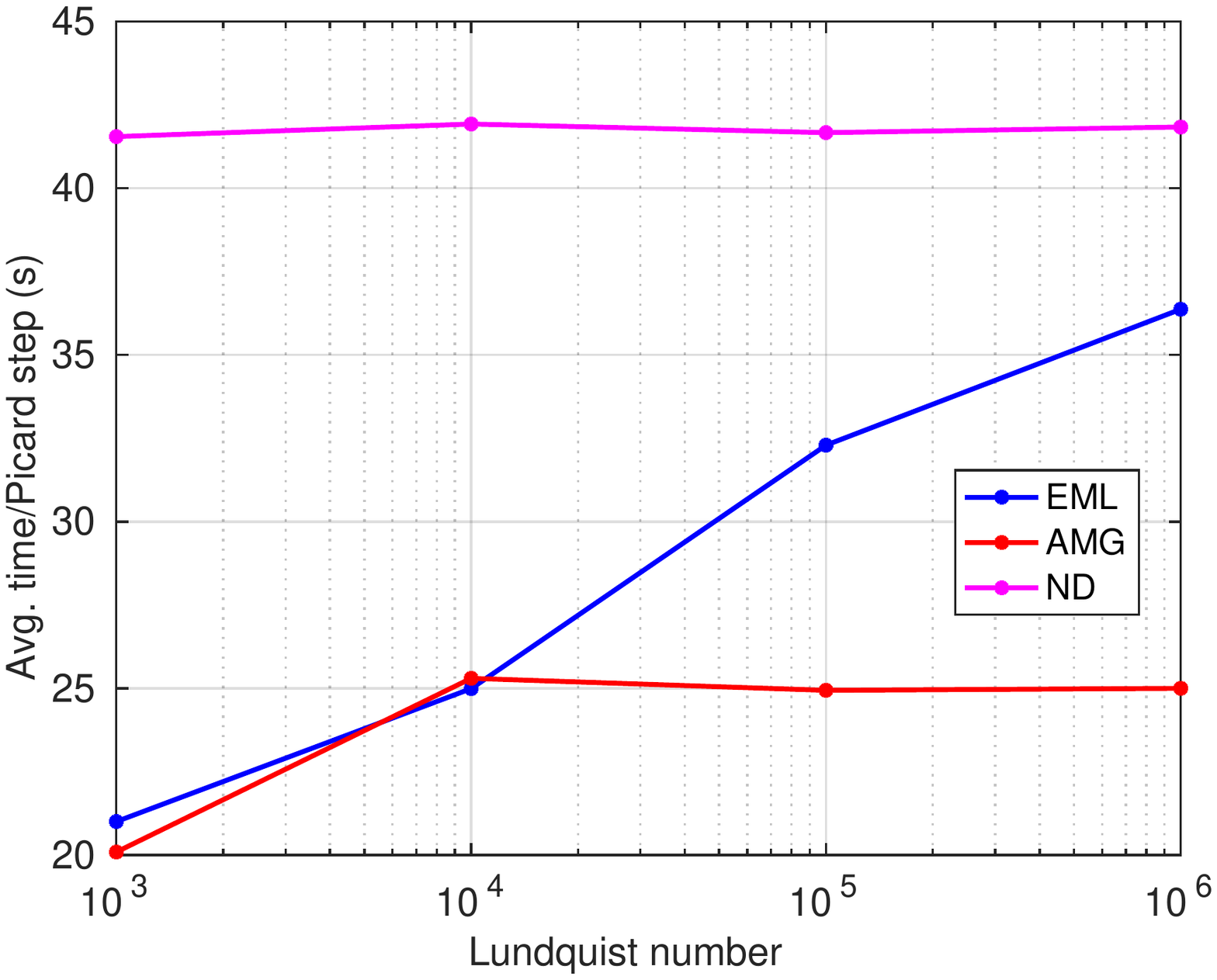}
  }

    \caption{2D island coalescence problem: Comparison of AMG and multilevel (EML) preconditioners with increase in Lundquist number $S$. In both cases we use the BFBT approximation for the inverse of the Schur complement. On the left is the average iterations per Picard step 
    and on the right is the average time taken per Picard step. We also show timings for the nested dissection (ND) direct solver on the right for comparison.}
    \figlab{robustness_EML_vs_AMG}
\end{figure}

Now we test the robustness of AMG and EML preconditioners with respect to Lundquist number for this problem. To that
extent we consider $64\times64$, $\p=6$ mesh and choose a time stepsize of $\Delta t=0.05$. The results are again 
averaged over six time steps. Figure \figref{robustness_EML_vs_AMG} shows the average iteration counts and time 
per Picard step for Lundquist numbers in the range $[10^3,10^6]$. We can see that the BFBT+AMG preconditioner is more robust with respect to increase in Lundquist numbers than the BFBT+EML preconditioner. Nevertheless, the growth in iterations 
for the BFBT+EML preconditioner is still moderate and in all the cases both the preconditioners take less time than the nested dissection (ND) direct solver. In our future work, we want to improve the robustness of the BFBT+EML preconditioner by exploring other fine scale solvers than block-Jacobi.

\begin{figure}[h!b!t!]
\centering
\includegraphics[trim=0cm 5cm 1cm 6cm, clip=true,width=0.7\textwidth]{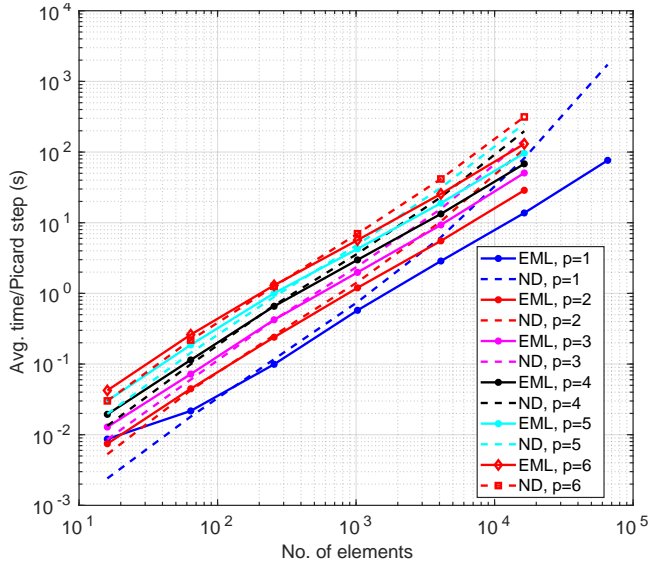}
    \caption{2D island coalescence problem: Comparison of scaling of BFBT+EML preconditioned GMRES (denoted as EML) and ND direct solver with number of elements for solution orders $\p={1-6}$.
}
\figlab{EML_ND_scaling}
\end{figure}


Finally, we will compare the BFBT+EML preconditioned GMRES with one of the fast direct solvers in 2D namely the nested dissection (ND). We use the 
UMFPACK library \cite{Davis04} for this purpose and order the trace
system matrix \eqref{eq:MHD_HDG_saddle} in nested dissection manner.
In Figure \figref{EML_ND_scaling} we compare the scaling of the average time taken per Picard step with mesh refinements for the BFBT+EML preconditioned GMRES and the nested dissection direct solver for different solution orders. The dominant cost in the BFBT+EML 
preconditioned solver comes from the factorization involved in EML and as per the theoretical complexities derived in \cite{muralikrishnan2020multilevel}, the iterative solver at all orders show close to linear scaling with the number of elements whereas the ND direct solver shows an asymptotic scaling of $\mc{O}(\Nel^{3/2})$.

\begin{figure}[h!t!b!]
  \subfigure[]{
    \figlab{EML_vs_ND_p_1}
    \includegraphics[trim=2cm 6cm 2.8cm 7cm,clip=true,width=0.48\columnwidth]{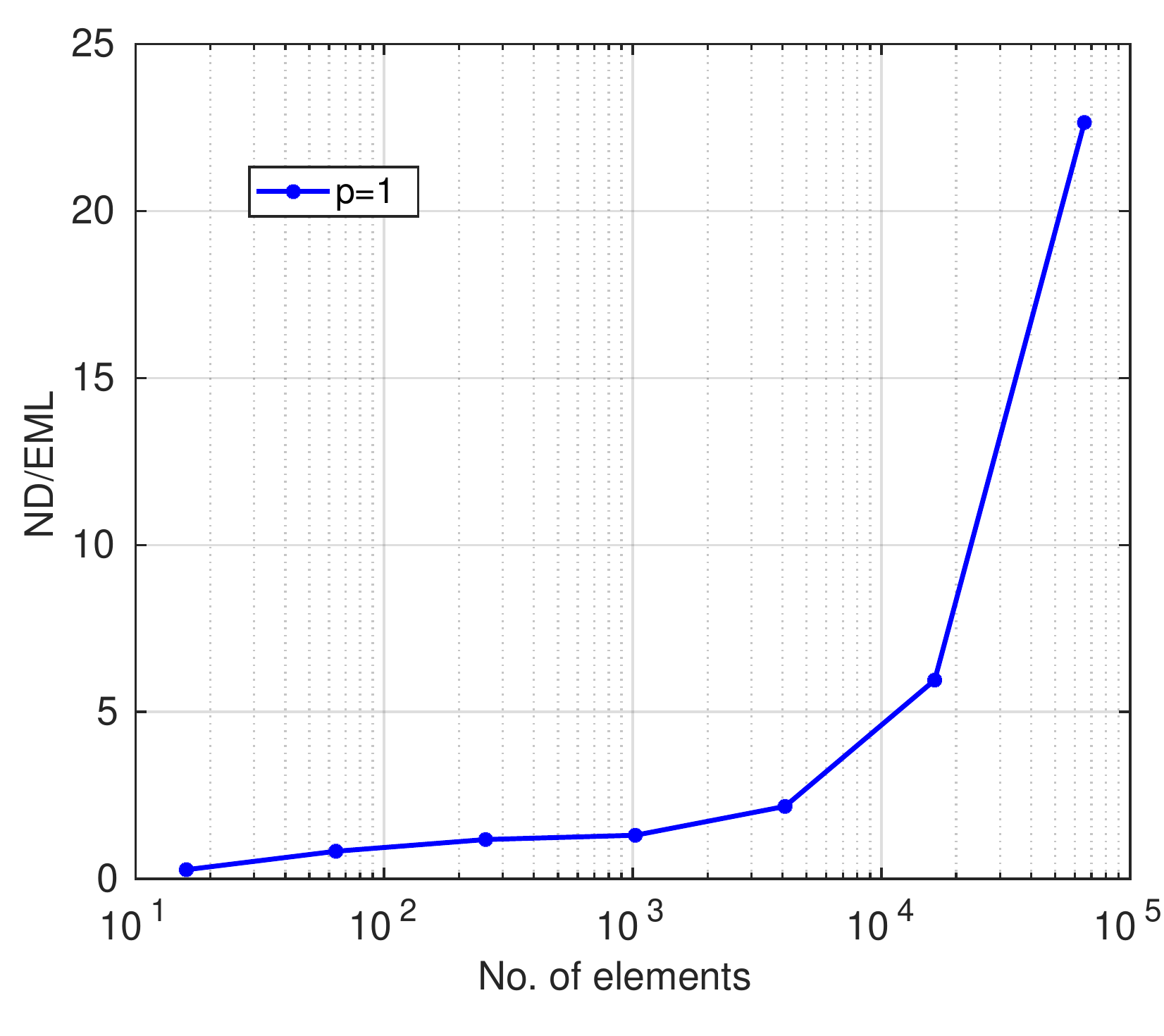}
  }
  \subfigure[]{
    \figlab{EML_vs_ND_p_rest}
    \includegraphics[trim=2cm 6cm 3cm 7cm,clip=true,width=0.48\columnwidth]{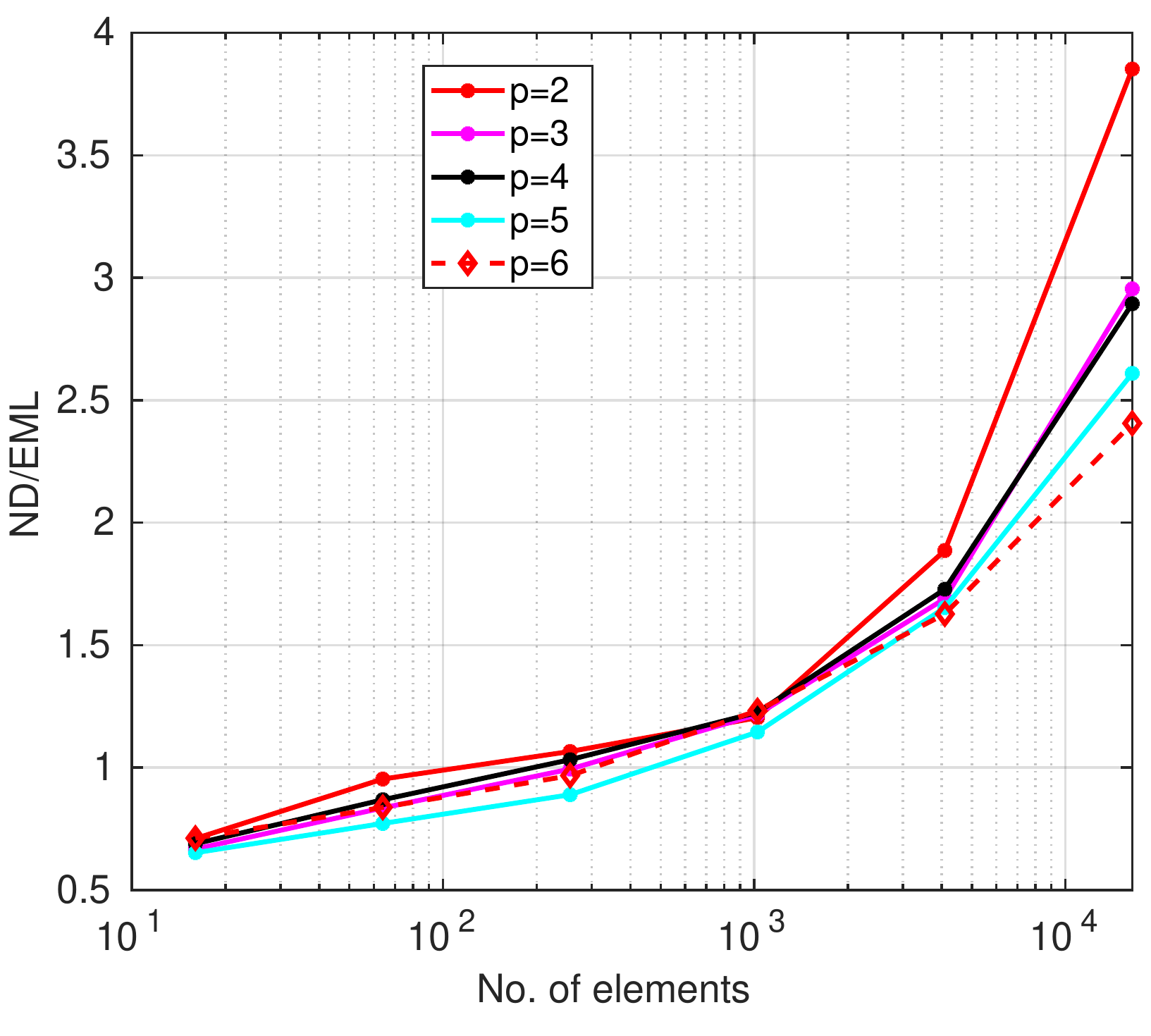}
  }
\caption{2D island coalescence problem: Ratio of time taken per Picard step for ND solver over BFBT+EML preconditioned GMRES for solution order $\p=1$ (left) and for orders $\p={2-6}$ (right). The ND solver ran into out of memory issues for orders $\p>1$ and $256^2$ elements.}
  \figlab{EML_vs_ND}
\end{figure}

In Figure \figref{EML_vs_ND} we show the speedup of BFBT+EML preconditioned GMRES compared to ND and we can see at all solution orders, when the 
number of elements is greater than $10^3$ (after 5 uniform refinements) the iterative solver is faster than ND.
We get a maximum speedup of approximately $23$ for $\p=1$ and $256\times256$ elements, and the ND solver ran into out of 
memory issues after this. In terms
of memory for $\p=1$ and $256\times256$ mesh, the ND solver needed $98.2$ GB for the L, U factors whereas the EML part in the iterative solver needed
only $2.4$ GB which is 41 times less memory compared to the direct solver. This speedup and memory reduction is significantly higher than the
one observed for the scalar problems studied in \cite{muralikrishnan2020multilevel}. Thus BFBT+EML preconditioned GMRES can deliver 
significant speedups and memory savings compared to the ND direct solver for vector valued problems even in 2D.


\section{Conclusion}
\seclab{conclusions}

In this work we present a block preconditioning strategy for the trace system coming from the HDG discretization of the incompressible resistive
MHD equations. In the block preconditioner, we use least squares commutator (BFBT) approximation for the inverse of the Schur complement
and a system projection AMG v-cycle for the approximate inverse of the nodal block. For the smoother inside AMG cycle, we compare preconditioned GMRES and 
ILU(0) smoother of overlap one and conclude that the GMRES smoother is faster, requires less memory and more robust compared to the ILU smoother. We test the performance of the
block preconditioner on several 2D and 3D transient test cases including, but not limited to, the island coalescence problem at high Lundquist numbers and demonstrate robustness and 
parallel scalability up to thousands of cores. We also show the application of the multilevel approximate nested disection preconditioner introduced in \cite{muralikrishnan2020multilevel} for the
approximate inverse of the nodal block and compare the performance with AMG and a full nested dissection direct solver. The BFBT+multilevel preconditioner shows better algorithmic
scalability compared to BFBT+AMG with respect to mesh refinements. In terms of robustness with respect to Lundquist numbers BFBT+AMG performs better
and strong smoothers are needed in the BFBT+multilevel preconditioner. In comparison with the nested dissection direct solver BFBT+multilevel
preconditioned GMRES is significantly faster and requires lot less memory up to an order of magnitude.

In terms of future works improving the parallel performance and scalability of the block preconditioner is one of the primary areas of focus.
This would enable large scale 3D simulation of realistic geometries that are of interest in fusion research as well as other areas of plasma physics research.
Since, the block preconditioner developed here is also applicable for Stokes and incompressible Navier Stokes equations studying the 
performance in those contexts is also of interest.

\section*{Acknowledgements}
The authors would like to thank Dr. Paul Lin for  help regarding various aspects of Trilinos. SM, SS and TBT are partially supported by the Department of Energy (grant DE-SC0018147), the National Science Foundation (grants NSF-DMS1620352, Early Career NSF-OAC1808576, and NSF-OAC1750863), and the Defense Threat Reduction Agency (grant HDTRA1-18-1-0020). We are grateful for the support. 
At Sandia this work was partially supported by the U.S. Department of Energy, Office of Science, Office of Advanced Scientific Computing Research, Applied Mathematics Program and by the U.S. Department of Energy, Office of Science, Office of Advanced Scientific Computing Research and Office of Fusion Energy Sciences, Scientific Discovery through the Advanced Computing (SciDAC) program.  Sandia National Laboratories is a multi-mission laboratory managed and operated by National Technology and Engineering Solutions of Sandia, 
LLC., a wholly owned subsidiary of Honeywell International, Inc., for the U.S. 
Department of Energy's National Nuclear Security Administration under contract 
DE-NA0003525. This paper describes objective technical results and analysis. Any subjective views or opinions that might be expressed in the paper do not necessarily represent the views of the U.S. Department of Energy or the United States Government.
The authors acknowledge the Texas Advanced Computing Center (TACC) at The University of Texas at Austin for providing HPC, visualization and storage resources that have contributed to the research results reported within this paper.  

\bibliography{references,ceo,DtNbib}

\end{document}